\documentclass[11pt,a4paper,usenames,dvipsnames]{article}

\usepackage{url,amsmath,amssymb,latexsym,mathrsfs,comment,amsthm,enumerate,geometry,calc,ifthen,mathdots,textcomp,multicol, adjustbox}

\usepackage{cases}

\usepackage[noadjust]{cite}

\usepackage[colorlinks]{hyperref}

\usepackage[margin=10pt,font=small,labelfont=bf, labelsep=period]{caption}

\usepackage{makecell}

\usepackage[all,cmtip,2cell]{xy}

\usepackage[x11names, svgnames, rgb,table]{xcolor}
\usepackage{hhline}
\usepackage{multirow}

\usepackage{pst-node}
\usepackage{tikz-cd} 
\tikzcdset{scale cd/.style={every label/.append style={scale=#1},
    cells={nodes={scale=#1}}}}

\usepackage{enumitem}

\usepackage[utf8]{inputenc}
\usepackage[T1]{fontenc}

\usepackage{tikz}
\usepgflibrary{snakes,arrows,shapes}
\pgfdeclarelayer{background layer}
\pgfsetlayers{background layer,main}

\usetikzlibrary{decorations.markings}
\usetikzlibrary{arrows,matrix}
\usepgflibrary{arrows}
  
\tikzset{->-/.style={decoration={
  markings,
  mark=at position #1 with {\arrow{{latex}}}},postaction={decorate}}}

\newcommand\nc\newcommand
\nc\rnc\renewcommand

\nc\trans[1]{\left(\begin{smallmatrix}#1\end{smallmatrix}\right)}
\nc\sm\setminus

\nc{\mtlab}[1]{\mapstochar \xrightarrow {\ #1\ }}

\nc{\ssim}{\mathrel{\raise0.25 ex\hbox{\oalign{$\approx$\crcr\noalign{\kern-0.84 ex}$\approx$}}}}

\nc{\ldb}{[\hspace{-0.55truemm}[}
\nc{\rdb}{]\hspace{-0.55truemm}]}

\nc\larr{\mathrel{\begin{tikzpicture}\draw[>-](0,0)--(.5,0);\end{tikzpicture}}}
\nc\rarr{\mathrel{\begin{tikzpicture}\draw[->](0,0)--(.5,0);\end{tikzpicture}}}
\nc\llarr{\mathrel{\begin{tikzpicture}\draw[>-<](0,0)--(.5,0);\end{tikzpicture}}}
\nc\rrarr{\mathrel{\begin{tikzpicture}\draw[<->](0,0)--(.5,0);\end{tikzpicture}}}
\nc\lrarr{\mathrel{\begin{tikzpicture}\draw[>->](0,0)--(.5,0);\end{tikzpicture}}}
\nc\arrl{\mathrel{\begin{tikzpicture}\draw[-<](0,0)--(.5,0);\end{tikzpicture}}}
\nc\arrr{\mathrel{\begin{tikzpicture}\draw[<-](0,0)--(.5,0);\end{tikzpicture}}}

\nc\Z{\mathbb Z}

\nc\dda{{\scalebox{0.5}{\begin{tikzpicture}\draw[thick,>->](0,.5)--(0,0);\end{tikzpicture}}}}

\nc\pre\preceq

\newcommand{\uv}[1]{\fill (#1,2)circle(.17);}
\newcommand{\lv}[1]{\fill (#1,0)circle(.17);}

\newcommand{\uvs}[1]{{\foreach \x in {#1} { \uv{\x}}}}
\newcommand{\lvs}[1]{{\foreach \x in {#1} { \lv{\x}}}}
\newcommand{\darcx}[3]{\draw(#1,0)arc(180:90:#3) (#1+#3,#3)--(#2-#3,#3) (#2-#3,#3) arc(90:0:#3);}
\newcommand{\darc}[2]{\darcx{#1}{#2}{.4}}
\newcommand{\uarcx}[3]{\draw(#1,2)arc(180:270:#3) (#1+#3,2-#3)--(#2-#3,2-#3) (#2-#3,2-#3) arc(270:360:#3);}
\newcommand{\uarc}[2]{\uarcx{#1}{#2}{.4}}
\newcommand{\stline}[2]{\draw(#1,2)--(#2,0);}
\nc{\uarcs}[1]{{\foreach \x/\y in {#1}{ \uarc{\x}{\y} }}}
\nc{\darcs}[1]{{\foreach \x/\y in {#1}{ \darc{\x}{\y} }}}

\nc\udotted[2]{\draw[dotted](#1+.5,2)--(#2-.5,2);}
\nc\ldotted[2]{\draw[dotted](#1+.5,0)--(#2-.5,0);}

\nc\mot[2]{{\lower1.0 ex\hbox{
\begin{tikzpicture}[scale=.23]
\uv1\uv2\uv3\lv1\lv2\lv3
\foreach \x in {#1} {\stline\x\x}
\foreach \x/\y in {#2} {\uarc\x\y \darc\x\y}
\end{tikzpicture}
}}}

\nc\mott[3]{{\lower1.0 ex\hbox{
\begin{tikzpicture}[scale=.23]
\uv1\uv2\uv3\lv1\lv2\lv3
\foreach \x/\y in {#1} {\stline\x\y}
\foreach \x/\y in {#2} {\uarc\x\y}
\foreach \x/\y in {#3} {\darc\x\y}
\end{tikzpicture}
}}}

\nc\moot[2]{{\lower1.4 ex\hbox{
\begin{tikzpicture}[scale=.3]
\uv1\uv2\uv3\uv4\lv1\lv2\lv3\lv4
\foreach \x in {#1} {\stline\x\x}
\foreach \x/\y in {#2} {\uarc\x\y \darc\x\y}
\end{tikzpicture}
}}}

\newcommand{\darcxcol}[4]{\draw[#4](#1,0)arc(180:90:#3) (#1+#3,#3)--(#2-#3,#3) (#2-#3,#3) arc(90:0:#3);}
\newcommand{\darccol}[3]{\darcxcol{#1}{#2}{.4}{#3}}
\nc{\darccols}[2]{{\foreach \x/\y in {#1}{ \darccol{\x}{\y}{#2} }}}
\newcommand{\uarcxcol}[4]{\draw[#4](#1,2)arc(180:270:#3) (#1+#3,2-#3)--(#2-#3,2-#3) (#2-#3,2-#3) arc(270:360:#3);}
\newcommand{\uarccol}[3]{\uarcxcol{#1}{#2}{.4}{#3}}
\nc{\uarccols}[2]{{\foreach \x/\y in {#1}{ \uarccol{\x}{\y}{#2} }}}

\newcounter{ncols}
\newcounter{incols}

\nc\la\langle
\nc\ra\rangle

\nc\ol\overline
\nc\ul\underline
\nc\da\downarrow
\nc\pr\bullet
\nc\ccirc\circledcirc
\nc\bc\odot

\nc\RSS{{\bf RSS}}
\nc\BS{{\bf BS}}
\nc\RBS{{\bf RBS}}
\nc\SBS{{\bf SBS}}
\nc\RSBS{{\bf RSBS}}
\nc\PA{{\bf PA}}
\nc\OG{{\bf OG}}
\nc\CPG{{\bf CPG}}
\nc\Set{{\bf Set}}
\nc\Sgp{{\bf Sgp}}
\nc\IS{{\bf IS}}
\nc\IG{\operatorname{\textup{\textsf{IG}}}}
\nc\RIG{\operatorname{\textup{\textsf{RIG}}}}
\nc\bG{{\bf G}}
\nc\bS{{\bf S}}
\nc\bE{{\bf E}}
\nc\bP{{\bf P}}
\nc\bEE{{\bf E}}
\nc\bPP{{\bf P}}
\nc\bB{{\bf B}}
\nc\bC{{\bf C}}
\nc\bD{{\bf D}}
\nc\bbG{\mathbb G}
\nc\bbS{\mathbb S}

\nc\PG{\operatorname{\textup{\textsf{PG}}}}
\nc\PGPz{\PG(P_0)}

\nc\JEnew[1]{{\color{WildStrawberry} #1}}
\nc\BOB[1]{{\color{RoyalPurple} #1}}

\nc\ds\displaystyle
\nc\ts\textstyle

\nc\bn{{\bf n}}

\nc\al\alpha
\nc\be\beta
\nc\ga\gamma
\nc\de\delta
\nc\ve\varepsilon
\nc\lam\lambda
\nc\om\omega
\nc\ze\zeta
\nc\si\sigma
\nc\Si\Sigma
\nc\Ga\Gamma
\nc\De\Delta
\nc\Om\Omega
\nc\io\iota

\nc\nab\nabla

\nc\vt\vartheta
\rnc\th\theta
\nc\Th\Theta

\nc\BY{\qquad\text{by}\qquad}
\nc\GIVENBY{\qquad\text{given by}\qquad}
\nc\ISGIVENBY{\qquad\text{is given by}\qquad}
\nc\OR{\qquad\text{or}\qquad}
\nc\Or{\quad\text{or}\quad}
\nc\AND{\qquad\text{and}\qquad}
\nc\ANDSIM{\qquad\text{and similarly}\qquad}
\nc\ANDSIm{\quad\text{and similarly}\quad}
\nc\ANDSO{\qquad\text{and so}\qquad}
\nc\ANd{\quad\text{and}\quad}
\nc\COMMA{,\qquad}
\nc\COMMa{,\quad}
\nc\WHERE{\qquad\text{where}\qquad}

\rnc\iff{\ \Leftrightarrow\ }
\nc\IFf{\quad \Leftrightarrow\quad }
\nc\Iff{\ \ \Leftrightarrow\ \ }
\nc\IFF{\qquad \Leftrightarrow\qquad }
\rnc\implies{\ \Rightarrow\ }
\nc\IMPLIES{\qquad \Rightarrow\qquad }

\nc\set[2]{\{#1:#2\}}
\nc\bigset[2]{\big\{#1:#2\big\}}
\nc\pres[2]{\la#1:#2\ra}

\nc\bit{\begin{itemize}[label=\textbullet, leftmargin=5mm]}
\nc\eit{\end{itemize}}
\nc\ben{\begin{enumerate}[label=\textup{(\roman*)},leftmargin=10mm]}
\nc\bena{\begin{enumerate}[label=\textup{(\alph*)},leftmargin=10mm]}
\nc\een{\end{enumerate}}
\nc\bmc{\begin{multicols}}
\nc\emc{\end{multicols}}

\nc\pf{\begin{proof}}
\nc\epf{\end{proof}}
\nc\pfclaim{\begin{quote}\begin{proof}}
\nc\epfclaim{\end{proof}\end{quote}}
\nc\epfres{\hfill\qed}
\nc\epfreseq{\tag*{\qed}}
\let\oldproofname=\proofname
\renewcommand{\proofname}{\rm\bf{\oldproofname}}
\nc{\pfitem}[1]{\medskip \noindent #1.}
\nc{\firstpfitem}[1]{#1.}
\nc{\pfcase}[1]{\medskip\noindent {\bf Case #1.}}
\nc\aftercases{\medskip\noindent}

\renewcommand{\H}{\mathrel{\mathscr H}}
\renewcommand{\L}{\mathrel{\mathscr L}}
\newcommand{\R}{\mathrel{\mathscr R}}
\newcommand{\D}{\mathrel{\mathscr D}}
\newcommand{\J}{\mathrel{\mathscr J}}
\newcommand{\K}{\mathbb K}
\nc\rH{\mathrel{\H}}
\nc\rL{\mathrel{\L}}
\nc\rR{\mathrel{\R}}
\nc\rD{\mathrel{\D}}
\nc\rJ{\mathrel{\J}}
\nc\rK{\mathrel{\K}}
\nc\rsi{\mathrel{\si}}

\nc\leqL{\leq_\L}
\nc\leqR{\leq_\R}
\nc\leqJ{\leq_\J}
\nc\leqH{\leq_\H}
\nc\leqF{\leq_{\F}}
\nc\geqF{\geq_{\F}}

\newcommand{\id}{\operatorname{id}}

\nc\mr\mathrel
\nc\pc[2]{(#1,#2)^\sharp}

\nc\U{{\bf U}}
\nc\bF{{\bf F}}
\nc\GG{{\bf G}}

\nc\V{\mathcal V}
\nc\G{\mathcal G}
\rnc\iff{\ \Leftrightarrow\ }
\rnc\implies{\ \Rightarrow\ }
\nc\Implies{\quad \Rightarrow\quad }
\nc\F{\mathrel{\mathscr F}}
\nc\C{\mathscr C}
\nc\M{\mathcal M}
\nc\CC{\mathcal C}
\nc\EP{\bE(P)}
\nc\EPd{\bE(P')}
\nc\PE{\bP(E)}
\nc\PEd{\bP(E')}
\nc\DD{\mathcal D}
\nc\I{\mathcal I}
\rnc\O{\mathcal O}
\rnc\P{\mathscr P}
\nc\TL{\mathcal T\!\mathcal L}
\nc\PP{\mathcal P}
\nc\T{\mathcal T}
\nc\p{\mathfrak p}
\nc\q{\mathfrak q}
\rnc\r{\mathfrak r}
\nc\s{\mathfrak s}
\rnc\t{\mathfrak t}
\nc\bd{{\bf d}}
\nc\br{{\bf r}}
\nc\op\oplus
\nc\lra{\mathrel\leftrightarrow}
\nc\es\varnothing
\nc\rev{\textup{rev}}
\nc\corestt{{\upharpoonleft}}
\nc\restt{{\upharpoonright}}
\nc\corest{{\downharpoonleft}}
\nc\rest{{\downharpoonright}}
\nc\WHERe{\quad\text{where}\quad}
\nc\ot\leftarrow
\rnc\a{\mathfrak a}
\rnc\b{\mathfrak b}
\rnc\c{\mathfrak c}
\rnc\d{\mathfrak d}
\nc\sub\subseteq
\nc\mt\mapsto
\nc\im{\operatorname{im}}
\nc\B{\mathcal B}
\nc\PB{\PP\B}
\nc\E{\mathbb E}
\nc\BP{\operatorname{\textup{\textsf{BP}}}}
\rnc\SS{\operatorname{\textup{\textsf{S}}}}
\nc\MM{\operatorname{\textup{\textsf{M}}}}

\numberwithin{equation}{section}

\newtheorem{thm}[equation]{Theorem}
\newtheorem{lemma}[equation]{Lemma}

\newtheorem{prop}[equation]{Proposition}

\theoremstyle{definition}

\newtheorem{defn}[equation]{Definition}
\newtheorem{rem}[equation]{Remark}
\newtheorem{eg}[equation]{Example}

\newcounter{caseco}

\newcounter{subcaseco}

\newcounter{stepco}

\newcounter{stageco}

\makeatletter
\def\blfootnote{\gdef\@thefnmark{}\@footnotetext}
\makeatother

\begin{document}

\title{Projection algebras and free projection- and idempotent-generated regular $*$-semigroups}

\date{}
\author{}

\maketitle

\vspace{-15mm}

\begin{center}
{\large 
James East,%
\hspace{-.25em}\footnote{\label{fn:JE}Centre for Research in Mathematics and Data Science, Western Sydney University, Locked Bag 1797, Penrith NSW 2751, Australia. {\it Emails:} {\tt J.East@westernsydney.edu.au}, {\tt A.ParayilAjmal@westernsydney.edu.au}.}
Robert D.~Gray,%
\hspace{-.25em}\footnote{School of Mathematics, University of East Anglia, Norwich NR4 7TJ, England, UK. {\it Email:} {\tt Robert.D.Gray@uea.ac.uk}.}
P.A.~Azeef Muhammed,%
\hspace{-.25em}\textsuperscript{\ref{fn:JE}}
Nik Ru\v{s}kuc%
\footnote{Mathematical Institute, School of Mathematics and Statistics, University of St Andrews, St Andrews, Fife KY16 9SS, UK. {\it Email:} {\tt Nik.Ruskuc@st-andrews.ac.uk}}
\blfootnote{This work was supported by the following grants:
Future Fellowship FT190100632 of the Australian Research Council;
EP/V032003/1, EP/S020616/1 and EP/V003224/1 of the Engineering and Physical Sciences Research Council.
The second author thanks the Sydney Mathematical Research Institute, the University of Sydney
and Western Sydney University for partially funding his visit to Sydney in 2023 during which part of
this research was undertaken.  We thank the referee for their careful reading of the paper, and for their helpful suggestions.
}
}
\end{center}

\maketitle

\begin{abstract}
\noindent
The purpose of this paper is to introduce a new family of semigroups---the free projection-generated regular $*$-semigroups---and initiate their systematic study.  
Such a semigroup~$\PG(P)$ is constructed from a projection algebra $P$, using the recent groupoid approach to regular $*$-semigroups.  
The assignment $P\mt\PG(P)$ is a left adjoint to the forgetful functor that maps a regular $*$-semigroup $S$ to its projection algebra $\bP(S)$.  
In fact, the category of projection algebras is coreflective in the category of regular $*$-semigroups.
The algebra $\bP(S)$ uniquely determines the biordered structure of the idempotents $\bE(S)$, up to isomorphism, and this leads to a category equivalence between projection algebras and regular $*$-biordered sets.  
As a consequence, $\PG(P)$ can be viewed as a quotient of the classical free idempotent-generated (regular) semigroups $\IG(E)$ and $\RIG(E)$, where $E=\bE(\PG(P))$; this is witnessed by a number of presentations in terms of generators and defining relations.  
%
%
The semigroup $\PG(P)$ can also be interpreted topologically, through a natural link to the fundamental groupoid of a simplicial complex explicitly constructed from~$P$.
The above theory is illustrated on a number of examples.
In one direction, the free construction applied to the projection algebras of adjacency semigroups yields a new family of graph-based path semigroups.
In another, it turns out that, remarkably, the Temperley--Lieb monoid $\TL_n$ is the free regular $*$-semigroup over its own projection algebra~$\bP(\TL_n)$.

\medskip

\noindent
\emph{Keywords}: 
Regular $*$-semigroup, projection algebra, chained projection groupoid, free projection-generated regular $*$-semigroup, 
forgetful functor, adjoint, coreflectivity,
regular semigroup, biordered set, free idempotent-generated semigroup, 
presentation, fundamental groupoid, Temperley--Lieb monoid.
\medskip

\noindent
MSC: 
20M50 (primary),  
18B40,  
20M05, 
20M10,  
20M17,  
20M20.  

\end{abstract}

\tableofcontents



\section{Introduction}\label{sect:intro}

Set-based free objects exist for many classes of algebras, such as groups, monoids, lattices, rings and modules.  These are typically built from a base set, and defined in terms of a universal mapping property.  Formally, the existence of such free algebras in a category $\bC$ amounts to the forgetful functor $\bC\to\Set$ (which maps an algebra to its underlying set) having a left adjoint $\Set\to\bC$; full definitions will be recalled later in the paper.  Such an adjoint exists for example if $\bC$ forms a variety \cite{Birkhoff1935,Bergman2015,BS1981}, but this is not the case in general, even for some `classical' algebras such as fields \cite[Exercise 2.3]{Bergman2015}.
From the early days of semigroup theory, it was recognised that the idempotents form a useful `skeleton' of a semigroup.  As the theory developed, it emerged that this phenomenon was governed by the existence of forgetful functors from various classes of semigroups into categories of idempotent-like structures such as semilattices, biordered sets and others.  Adjoints of some of these functors led to important classes of free semigroups, which will be discussed more below.

The current paper is concerned with \emph{regular $*$-semigroups}.  These were introduced in \cite{NS1978} as an intermediate class between inverse semigroups and regular semigroups.  They have attained prominence recently, as the so-called diagram monoids come equipped with a natural regular $*$-structure \cite{EG2017,Maz2002,EF2012,ACHLV2015,Auinger2014,BDP2002}.  These monoids are the building blocks of diagram algebras, such as the Brauer, Temperley--Lieb and partition algebras, among others, which in turn have important applications in theoretical physics, low-dimensional topology, representation theory, and many other parts of mathematics and science \cite{HR2005,Brauer1937,TL1971,Martin1994,Jones1994_2,Kauffman1987}.
Every regular $*$-semigroup~$S$ contains a set $\bP(S)$ of distinguished idempotents known as projections, which can be given the structure of a \emph{projection algebra} \cite{Imaoka1983}.  The category $\PA$ of such algebras played a key role in the groupoid representation of regular $*$-semigroups in the recent paper \cite{EPA2024}.  The assignment $S\mt\bP(S)$ is a forgetful functor from the category $\RSS$ of regular $*$-semigroups to $\PA$.
It turns out that this functor has a left adjoint $\PA\to\RSS$, and this leads to the notion of a \emph{free (projection-generated) regular $*$-semigroup} $\PG(P)$ over a projection algebra $P$.  The purpose of this paper is to show how to construct these free semigroups, and initiate their systematic study.

The link between regular $*$-semigroups and their projection algebras has many parallels with the link between arbitrary semigroups and their biordered sets of idempotents; the latter form the cornerstone of Nambooripad's theory of regular semigroups \cite{Nambooripad1979}.  In that situation we have a forgetful functor $\bE:\Sgp\to\BS$, which maps a semigroup $S$ to its biordered set
\[
E = \bE(S) = \set{e\in S}{e^2=e}.
\]
The latter has the structure of a partial algebra, where a product~$ef$ (for $e,f\in E$) is only defined if at least one of $ef$ or $fe$ is equal to $e$ or $f$.  
A deep result of Easdown and Nambooripad states that $\bE$ has a left adjoint $\BS\to\Sgp$.  This formulation can be found in \cite[Theorem 6.10]{Nambooripad_book}, but has its basis in \cite[Theorem 3.3]{Easdown1985}.
The adjoint of $\bE$ maps an (abstract) biordered set $E$ to the \emph{free (idempotent-generated) semigroup} $\IG(E)$, which is defined by the presentation
\[
\IG(E) = \pres{X_E}{x_ex_f=x_{ef}\ \text{if $ef$ is defined in $E$}},
\]
where here $X_E = \set{x_e}{e\in E}$ is an alphabet in one-one correspondence with $E$.  The key point is that the biordered set of the semigroup $\IG(E)$ is precisely $E$, when one identifies $e\in E$ with the equivalence class of the one-letter word $x_e$.

The biorder approach has additional power in the case of regular semigroups.  The restriction of the above forgetful functor $\bE:\Sgp\to\BS$ to regular semigroups maps into the category of regular biordered sets, which were axiomatised by Nambooripad \cite{Nambooripad1979}, but the adjoint $\BS\to\Sgp$ maps regular biordered sets to non-regular semigroups in general.  Instead, the restriction has a different adjoint, mapping $E$ to the \emph{free regular (idempotent-generated) semigroup} $\RIG(E)$; see~\cite{Nambooripad1979} for a combinatorial/topological definition, and~\cite{Pastijn1980} for a presentation.

The free semigroups $\IG(E)$ and $\RIG(E)$ turn out to have very intricate structure, and have therefore become a subject of broad interest in their own right \cite{
Nambooripad1979,
NP1980,
Easdown1985,
GR2012,
GR2012b,
DGR2017,
DG2014,
DG2016_2,
DGY2015,
BMM2009, 
Pastijn1980,
McElwee2002,
Dolinka2013,
Dandan2016,
Dolinka2012,
GY2014,
DR2013
}.  
As one strand of research, it was shown in \cite{BMM2009} that maximal subgroups of $\IG(E)$ and $\RIG(E)$ are (isomorphic to) fundamental groups of certain natural complexes associated to $E$.  The main motivation for this result was to address a folklore conjecture that such maximal subgroups were always free, and the topological viewpoint led to the discovery of a biordered set inducing the non-free subgroup $\Z\times\Z$; see \cite[Section 5]{BMM2009}.  Soon after, the conjecture was turned upside down, when it was shown in \cite{GR2012} that \emph{every} group is isomorphic to a maximal subgroup of some $\IG(E)$.  Since then, many substantial studies have emerged exploring the structure of $\IG(E)$ and $\RIG(E)$ in the case that $E=\bE(S)$ is the biordered set of some important semigroup $S$ \cite{GR2012b,DG2014,DGY2015,Dolinka2013}.

We note in passing that biordered sets are not the only structures that have been used to model the idempotent `skeleton' of semigroups.  For example, early work focussed on the so-called \emph{warp} of a semigroup $S$ \cite{Pastijn1980,Clifford1975,Clifford1976,Yamada1981}, which was again a partial algebra with underlying set $E = \bE(S)$, but which retained \emph{all} products $ef$ for which $e,f,ef\in E$.  Similarly, the biordered set can sometimes be given additional structure, as is indeed the case with regular $*$-semigroups, whose biordered sets have involutions \cite{NP1985}.

We now return to our main topic, regular $*$-semigroups and projection algebras, linked by the forgetful functor $\bP:\RSS\to\PA$.  The latter  maps a regular $*$-semigroup $S$ to its set
\[
P=\bP(S) = \set{p\in S}{p^2=p=p^*}
\]
of \emph{projections}, which is then given the structure of a \emph{projection algebra}, as originally defined (under a different name) by Imaoka \cite{Imaoka1983}.  This algebra has a unary operation $\th_p$ for each $p\in P$, which is defined by $q\th_p = pqp$ for $q\in P$.  Such algebras were axiomatised in \cite{Imaoka1983}, where it was shown that they are the appropriate vehicle for transformation representations of \emph{fundamental} regular $*$-semigroups, building on work of Munn in the inverse case \cite{Munn1970}; see also \cite{Jones2012}.  

Projection algebras took on a new level of significance in~\cite{EPA2024}, where they became the (structured) object sets of so-called \emph{chained projection groupoids}.  The main result \cite[Theorem 8.1]{EPA2024} is a category isomorphism $\RSS\cong\CPG$, where the latter is the category of all such groupoids.  These groupoids are in fact triples $(P,\G,\ve)$, where $P$ is an (abstract) projection algebra, $\G$ is an ordered groupoid whose structure has tight algebraic and order-theoretic links to $P$, and $\ve:\C\to\G$ is a functor from a natural \emph{chain groupoid} $\C=\C(P)$ built from $P$.  Seen through the groupoid lens, the forgetful functor $\CPG\to\PA$ maps $(P,\G,\ve)\mt P$.  The key construction in the current paper is a left adjoint $\PA\to\CPG$.  This maps $P\mt(P,\ol\C,\nu)$, where $\ol\C$ is an appropriate quotient of $\C$, and $\nu:\C\to\ol\C$ is the quotient map.  Topologically, $\ol\C$ is the fundamental groupoid of a natural complex built from~$P$.  
Applying the isomorphism $\bS:\CPG\to\RSS$ from~\cite{EPA2024} yields an adjoint $\PA\to\RSS$ to the forgetful functor $\bP:\RSS\to\PA$, and hence establishes the existence of the \emph{free (projection-generated) regular $*$-semigroups} $\PG(P) = \bS(P,\ol\C,\nu)$.  
In fact, we show that the adjoint is a right inverse of $\bP$, and this has the consequence that $\PA$ is \emph{coreflective} in $\RSS$ (and in $\CPG$).
We also remark that the existence of the adjoint answers a question not settled by the isomorphism $\RSS\cong\CPG$ from \cite{EPA2024}, in that it shows that every projection algebra $P$ can be realised as $P=\bP(S)$ for some regular $*$-semigroup $S$.  This fact was first established by Imaoka in his above-mentioned work on fundamental semigroups \cite{Imaoka1983}; see also~\cite{Jones2012}.

At this point, the semigroups $\PG(P)$ become the subject of the rest of the paper.  There are a number of ways to understand these semigroups, starting from their original definition as \emph{chain semigroups} $\PG(P) = \bS(P,\ol\C,\nu)$.  It is immediate from our construction that~$\PG(P)$ is generated by $P$.  Building on this, another of our main results establishes a presentation
\[
\PG(P) \cong \pres{X_P}{x_p^2=x_p,\ (x_px_q)^2=x_px_q,\ x_px_qx_p = x_{q\th_p}\ \text{for }p,q\in P},
\]
where $X_P = \set{x_p}{p\in P}$ is an alphabet in one-one correspondence with $P$.  Two further presentations 
allow us to understand the explicit relationship between the new regular $*$-semigroup $\PG(P)$ and the classical free idempotent-generated semigroups $\IG(E)$ and $\RIG(E)$, where $E=\bE(\PG(P))$.

As was the case with these classical semigroups, the free regular $*$-semigroups $\PG(P)$ are very interesting in their own right.  For example, when $P = \bP(A_\Ga)$ is the projection algebra of an adjacency semigroup (in the sense of Jackson and Volkov \cite{JV2010}), the free semigroup $\PG(P)$ is an apparently-new graph-based path semigroup.  As another example, it turns out that a finite \emph{Temperley--Lieb monoid} $\TL_n$ is a free regular $*$-semigroup over its own projection algebra~$\bP(\TL_n)$, giving yet another fundamental way to understand this important diagram monoid.  The situation for other diagram monoids is more delicate, and is taken up for the partition monoid in the forthcoming paper \cite{Paper4}.

We now give a brief overview of the structure of the paper; the introduction to each section contains a fuller summary of its contents.
%
%
\bit
\item
Sections \ref{sect:RSS} and \ref{sect:CPG} contain the preliminary material we need.  The former covers the basics on regular $*$-semigroups, and provides some key examples, namely adjacency semigroups and diagram monoids.  
The latter gives an overview of the relevant constructions and results from~\cite{EPA2024}, concerning regular $*$-semigroups, projection algebras and chained projection groupoids, leading up to the isomorphism $\RSS\cong\CPG$.  
\item
Sections \ref{sect:CP} and \ref{sect:F} are central for this paper: they introduce the semigroups $\PG(P)$, and demonstrate their freeness.  The former is the content of Definition \ref{defn:C_P} and Theorem~\ref{thm:C_P}.  The latter is achieved in Theorem \ref{thm:CP}, which shows that the functor ${\PA\to\RSS:P\mt\PG(P)}$ is indeed a left adjoint to the forgetful functor ${\RSS\to\PA:S\mt\bP(S)}$, and also establishes the coreflectivity of $\PA$ in $\RSS$.  A number of semigroup-theoretic consequences are given in Theorems \ref{thm:free} and \ref{thm:free2}.  
\item
Section \ref{sect:E} explores the connection between projection algebras and regular $*$-biordered sets, as defined in \cite{NP1985}.  The main result here is Theorem \ref{thm:equiv}, which establishes an equivalence between the categories $\PA$ and $\RSBS$ of all such structures.  As a consequence, Theorem~\ref{thm:CPE} shows that the semigroups $\PG(P)$ are also free with respect to the forgetful functor ${\RSS\to\RSBS:S\mt\bE(S)}$.  En route, we show in Proposition~\ref{prop:E} that the projection algebra of a regular $*$-semigroup uniquely determines its biordered set, up to isomorphism.
\item
In Section \ref{sect:pres} we give three presentations (by generators and defining relations) for the semigroups $\PG(P)$.  The first, in Theorem \ref{thm:pres}, involves $P$ as a generating set.  The other two, in Theorems \ref{thm:presE} and \ref{thm:presE2}, utilise the generating set $E=\bE(\PG(P))$, and exhibit $\PG(P)$ as an explicit quotient of $\IG(E)$ and of $\RIG(E)$, respectively.  
\item
Sections \ref{sect:AGa} and \ref{sect:TLn} illustrate our theory on some important examples.
In the former, we will see that the free construction applied to an adjacency semigroup results in an apparently-new graph-based \emph{bridging path semigroup}, which we believe is worthy of further study.
In the latter, we return to diagram monoids, showing in Theorem \ref{thm:TLnPGP} that the free regular $*$-semigroup associated to the projection algebra of a finite Temperley--Lieb monoid~$\TL_n$ is, somewhat remarkably, isomorphic to $\TL_n$ itself.
\item
Finally, Section \ref{sect:top} provides a topological interpretation of the free regular $*$-semigroups, with Theorems \ref{thm:fundgpd} and \ref{thm:fundgp} establishing a link with the fundamental group(oid)s of certain natural graphs and complexes associated to projection algebras.  We conclude the paper by recasting our earlier examples using this topological framework, and exploring some further ones.  As an intriguing consequence, the semigroup-theoretic structure of the Temperley--Lieb monoid $\TL_n$ allows us to immediately deduce that the components of its associated complex are simply connected; this is not obvious, a priori, as these complexes are highly intricate.
\eit
A feature of the work presented here is that it establishes tight connections between different types of mathematical objects.
This has certainly caused some presentational and notational challenges for the authors.
In the hope of somewhat easing such challenges for the reader, we have collected the key notation in Tables \ref{tab:cats}--\ref{tab:sscP}. The information presented there might not make much sense at this stage of reading the paper, but we hope that the reader can use the tables as a ready reference throughout.

\begin{table}[p]
\begin{center}
\begin{tabular}{|l|l|l|l|}
\hline
\textbf{Notation} & \textbf{Specification} & \textbf{Notes} & \textbf{Reference} \\
\hline\hline
$\BS$ & biordered sets (bosets) & & Subsection \ref{subsect:Eprelim}\\
\hline
$\CPG$ & chained projection groupoids & \parbox[t]{35mm}{isomorphic to $\RSS$\\ via $\bG$ and $\bS$} & 
\parbox[t]{37mm}{Subsection \ref{subsect:CPG},\\ Theorem \ref{thm:iso} \vspace{1mm}}\\
\hline
$\OG$ & ordered groupoids & & Subsection \ref{subsect:CP}\\
\hline
$\PA$ & projection algebras &\parbox[t]{35mm}{coreflective in $\RSS$\\ and in $\CPG$} & 
\parbox[t]{37mm}{Subsection \ref{subsect:PA},\\ Theorems \ref{thm:CP} and \ref{thm:PCnu} \vspace{1mm}}\\
\hline
$\RBS$ & regular bosets & & Subsection \ref{subsect:Eprelim}\\
\hline
$\RSBS$ & regular $*$-bosets & \parbox[t]{35mm}{equivalent to $\PA$\\ via $\bE$ and $\bP$} & 
\parbox[t]{37mm}{Subsection \ref{subsect:Eprelim},\\ Theorem \ref{thm:equiv} \vspace{1mm}}\\
\hline
$\RSS$ & regular $*$-semigroups & \parbox[t]{35mm}{isomorphic to $\CPG$\\ via $\bG$ and $\bS$} & 
\parbox[t]{37mm}{Subsection \ref{subsect:dbp}, \\ Theorem \ref{thm:iso} \vspace{1mm} } \\
\hline
$\Set$ & sets &  & 
\parbox[t]{37mm}{Section \ref{sect:intro}, \\ Remark \ref{rem:coref} \vspace{1mm} } \\
\hline
$\Sgp$ & semigroups &  & 
\parbox[t]{37mm}{Section \ref{sect:intro}, \\ Remark \ref{rem:coref} \vspace{1mm} } \\
\hline
\end{tabular}
\caption{Large categories.}
\label{tab:cats}
\end{center}
\end{table}

\begin{table}[p]
\begin{center}
\begin{tabular}{|l|l|l|l|}
\hline
\textbf{Notation} & \textbf{Meaning} & \textbf{Notes} &  \textbf{Reference} \\
\hline\hline
$\C:\PA\rightarrow\OG$ & \parbox[t]{60mm}{$\C(P)$ -- chain groupoid associated\\ with projection algebra $P$ \vspace{1mm}} & & Subsection \ref{subsect:CP}\\
\hline
$\bE:\RSS\rightarrow\RSBS$ & \parbox[t]{60mm}{$\bE(S)$ -- boset associated with\\ regular $*$-semigroup $S$ \vspace{1mm}} & forgetful & Subsection \ref{subsect:Eprelim}\\
\hline
$\bE:\PA\rightarrow\RSBS$ & \parbox[t]{60mm}{$\bE(P)$ -- boset associated with\\ projection algebra $P$ \vspace{1mm}} & equivalence & Subsection \ref{subsect:EP}\\
\hline
$\bF:\PA\rightarrow\CPG$ & \parbox[t]{70mm}{$\bF(P)$ -- free chained projection groupoid\\ associated with projection algebra $P$ \vspace{1mm}} & adjoint & Subsection \ref{subsect:olC}\\
\hline
$\bG:\RSS\rightarrow\CPG$ & \parbox[t]{65mm}{$\bG(S)$ -- chained projection groupoid\\ associated with regular $*$-semigroup $S$ \vspace{1mm}} & isomorphism & Subsection \ref{subsect:iso}\\
\hline
$\bP:\RSS\rightarrow\PA$ & \parbox[t]{50mm}{$\bP(S)$ -- projection algebra of\\ regular $*$-semigroup $S$ \vspace{1mm}} & forgetful & Subsection \ref{subsect:PA}\\
\hline
$\bP:\RSBS\rightarrow\PA$ & \parbox[t]{50mm}{$\bP(E)$ -- projection algebra of\\ regular $*$-boset $E$ \vspace{1mm}} & equivalence & Subsection \ref{subsect:PA}\\
\hline
$\bS:\CPG\rightarrow\RSS$ & \parbox[t]{71mm}{$\bS(P,\G,\varepsilon)$ -- regular $*$-semigroup associated\\ with chained projection groupoid $(P,\G,\varepsilon)$ \vspace{1mm}} & isomorphism & Subsection \ref{subsect:iso}\\
\hline
\end{tabular}
\caption{Functors.}
\label{tab:funct}
\end{center}
\end{table}

\begin{table}[p]
\begin{center}
\begin{tabular}{|l|l|l|}
\hline
\textbf{Notation} & \textbf{Meaning} &  \textbf{Reference} \\
\hline\hline
$\C(P)$ & chain groupoid of $P$ & Subsection  \ref{subsect:CP}\\
\hline
$\overline{\C}(P)$ & reduced chain groupoid of $P$ & Subsection  \ref{subsect:olC}\\
\hline
$G_P$ & graph associated with $P$ & Subsection  \ref{subsect:top}\\
\hline
$\IG(E)$ & free idempotent-generated semigroup over $E$ & Subsection  \ref{subsect:presE}\\
\hline
$K_P,K_P'$ & two complexes associated with $P$ & Subsection  \ref{subsect:top}\\
\hline
$\PG(P)$ & free projection-generated regular $*$-semigroup over $P$ & Subsection  \ref{subsect:C_P}\\
\hline
$\P(P)$ & path category of $P$ & Subsection  \ref{subsect:CP}\\
\hline
$\RIG(E)$ & free idempotent-generated regular semigroup over $E$ & Subsection  \ref{subsect:RIG}\\
\hline
\end{tabular}
\caption{Semigroups, small categories and other structures associated with a projection algebra $P$ or boset $E$.}
\label{tab:sscP}
\end{center}
\end{table}

\begin{table}[p]
\begin{center}
\begin{tabular}{|l|l|l|}
\hline
\textbf{Notation} & \textbf{Name} &  \textbf{Reference} \\
\hline\hline
$A_\Gamma$ & adjacency semigroup of digraph $\Gamma$ & Subsection \ref{subsect:AGa}\\
\hline
$B_\Gamma$ & bridging path semigroup of digraph $\Gamma$ & Section \ref{sect:AGa}\\
\hline
$\B_n$ & Brauer monoid & Subsection \ref{subsect:DM}\\
\hline
$\M_n$ & Motzkin monoid & Subsection \ref{subsect:DM}\\
\hline
$\PP_n$ & partition monoid & Subsection \ref{subsect:DM}\\
\hline
$\PB_n$ & partial Brauer monoid & Subsection \ref{subsect:DM}\\
\hline
$\TL_n$ & Temperley--Lieb monoid & Subsection \ref{subsect:DM}\\
\hline
\end{tabular}
\caption{Concrete semigroups.}
\label{tab:sscP}
\end{center}
\end{table}

\newpage

\section{\boldmath Regular $*$-semigroups}\label{sect:RSS}

\subsection{Definitions and basic properties}\label{subsect:dbp}

Here we gather the background on regular $*$-semigroups that we will need in the rest of the paper.  For proofs of the various assertions, see for example \cite{NS1978,Imaoka1983,EPA2024}.  For more on semigroups in general we refer to \cite{Howie1995,CPbook}.

A \emph{regular $*$-semigroup} is a semigroup $S$ with a unary operation ${{}^*:S\to S:a\mt a^*}$ satisfying
\[
(a^*)^* = a = aa^*a \AND (ab)^* = b^*a^* \qquad\text{for all $a,b\in S$.}
\]
From the identity $a=aa^*a$, it is clear that a regular $*$-semigroup is (von Neumann) regular.  The identities $(a^*)^*=a$ and $(ab)^*=b^*a^*$ say that ${}^*$ is an involution.

We write $\RSS$ for the category of regular $*$-semigroups with $*$-morphisms, i.e.~maps ${\phi:S\to S'}$ satisfying
\[
(ab)\phi = (a\phi)(b\phi) \AND (a^*)\phi = (a\phi)^* \qquad\text{for all $a,b\in S$.}
\]

Given a regular $*$-semigroup $S$, we write
\[
\bPP(S) = \set{p\in S}{p^2=p=p^*} \AND \bEE(S) = \set{e\in S}{e^2=e}
\]
for the sets of all \emph{projections} and \emph{idempotents} of $S$, respectively.  Important properties of these elements include the following:
\begin{enumerate}[label=\textup{\textsf{(RS\arabic*)}},leftmargin=13mm]
\item \label{RS1} The projections are precisely the elements of the form $aa^*$, for $a\in S$.
\item \label{RS2} The product of two projections is an idempotent, but need not be a projection.
\item \label{RS3} Any idempotent $e$ is the product of two projections, namely $e=(ee^*)(e^*e)$.
\item \label{RS4} The product of two idempotents need not be an idempotent.
\item \label{RS5} For all $p,q\in \bPP(S)$ we have $pqp\in \bPP(S)$.
\item \label{RS6} More generally, we have $aqa^*\in \bPP(S)$ for all $a\in S$ and $q\in \bPP(S)$.
\end{enumerate}
Our seventh item is an elaboration on \ref{RS3}.  For projections $p,q\in \bPP(S)$, we write $p\F q$ if $p=pqp$ and $q=qpq$, and say that $p$ and $q$ are \emph{friendly}.
 It is easy to check that $ee^* \F e^*e$ for any idempotent $e\in \bEE(S)$.  Thus, \ref{RS3} says that every idempotent is a product of a pair of $\F$-related projections.  As noted on \cite[p.~20]{EPA2024}, such expressions are unique:
\begin{enumerate}[label=\textup{\textsf{(RS\arabic*)}},leftmargin=13mm]\addtocounter{enumi}{6}
\item \label{RS7} $pq = rs \iff [p=r\text{ and }q=s]$ for all $(p,q),(r,s)\in{\F}$.
\end{enumerate}
More generally, any product of idempotents in a regular $*$-semigroup is equal to a product $p_1\cdots p_k$ of projections satisfying $p_1\F\cdots\F p_k$, though such expressions need not be unique; see \cite[Proposition 3.16]{EPA2024}.

Because of \ref{RS5}, each projection $p\in P=\bPP(S)$ induces a map
\begin{equation}\label{eq:thp}
\th_p:P\to P \GIVENBY q\th_p=pqp \qquad\text{for $q\in P$.}
\end{equation}
Taking these maps as unary operations gives $P$ the structure of a so-called \emph{projection algebra}; we will discuss these more formally in Section \ref{sect:CPG}, and will use them extensively in the rest of the paper.

We now describe some examples that we will use to illustrate the ideas developed in the paper.  For more examples, see \cite{EPA2024} or \cite{NP1985}, but note that regular $*$-semigroups were called `special $*$-semigroups' in the latter.

\subsection{Adjacency semigroups}\label{subsect:AGa}

Adjacency semigroups were introduced in \cite{JV2010}, and were discussed as a key example in \cite{EPA2024}.  They can be viewed in several ways: as combinatorial completely $0$-simple regular $*$-semigroups, or as symmetrical square $0$-bands.  Here we take the graph theoretic approach of \cite{JV2010}.  General completely $0$-simple regular $*$-semigroups were discussed in detail in \cite[Section 3.3]{EPA2024}.

Let $\Ga = (P,E)$ be a symmetric, reflexive digraph, with vertex set $P$, and edge set $E\sub P\times P$.  The \emph{adjacency semigroup} $A_\Ga$ is the regular $*$-semigroup with:
\bit
\item underlying set $A_\Ga = (P\times P)\cup\{0\}$,
\item  involution $0^*=0$ and $(p,q)^*=(q,p)$, and
\item product $0^2=0=0(p,q)=(p,q)0$ and $(p,q)(r,s) = \begin{cases}
(p,s) &\text{if $(q,r)\in E$}\\
0 &\text{otherwise.}
\end{cases}$
\eit
We identify a vertex $p\in P$ with the pair $(p,p)\in A_\Ga$.  In this way, the projections and idempotents of $A_\Ga$ are the sets $P_0=P\cup\{0\}$ and $E_0=E\cup\{0\}$, and we have $(p,q)=pq$ for all $(p,q)\in E$.  The projection algebra operations are given, for $p,q\in P_0$, by
\begin{equation}\label{eq:thpAGa}
q\th_p = \begin{cases}
p &\text{if $(p,q)\in E$}\\
0 &\text{otherwise.}
\end{cases}
\end{equation}
In particular, $\th_0$ is the constant map with image $0$.  It also follows that ${\F} = E\cup\{(0,0)\}$.
%

\subsection{Diagram monoids}\label{subsect:DM}

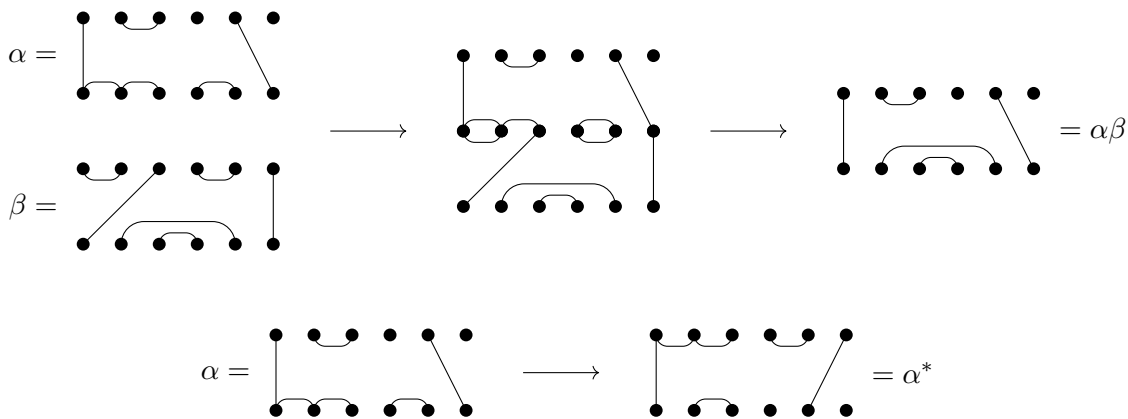
\begin{figure}[h]
\begin{center}
\begin{tikzpicture}[scale=.5]

\begin{scope}[shift={(0,0)}]	
\uvs{1,...,6}
\lvs{1,...,6}
\uarcx23{.3}
\darcx12{.3}
\darcx23{.3}
\darcx45{.3}
\stline11
\stline56
\draw(0.6,1)node[left]{$\alpha=$};
\draw[->](7.5,-1)--(9.5,-1);
\end{scope}

\begin{scope}[shift={(0,-4)}]	
\uvs{1,...,6}
\lvs{1,...,6}
\uarcx12{.3}
\uarcx45{.3}
\darcx34{.3}
\darcx25{.6}
\stline31
\stline66
\draw(0.6,1)node[left]{$\beta=$};
\end{scope}

\begin{scope}[shift={(10,-1)}]	
\uvs{1,...,6}
\lvs{1,...,6}
\uarcx23{.3}
\darcx12{.3}
\darcx23{.3}
\darcx45{.3}
\stline11
\stline56
\draw[->](7.5,0)--(9.5,0);
\end{scope}

\begin{scope}[shift={(10,-3)}]	
\uvs{1,...,6}
\lvs{1,...,6}
\uarcx12{.3}
\uarcx45{.3}
\darcx34{.3}
\darcx25{.6}
\stline31
\stline66
\end{scope}

\begin{scope}[shift={(20,-2)}]	
\uvs{1,...,6}
\lvs{1,...,6}
\uarcx23{.3}
\darcx34{.3}
\darcx25{.6}
\stline11
\stline56
\draw(6.4,1)node[right]{$=\alpha\beta$};
\end{scope}

\end{tikzpicture}\\[10mm]
%
%
%
%
\begin{tikzpicture}[scale=.5]

\begin{scope}[shift={(0,0)}]	
\uvs{1,...,6}
\lvs{1,...,6}
\uarcx23{.3}
\darcx12{.3}
\darcx23{.3}
\darcx45{.3}
\stline11
\stline56
\draw(0.6,1)node[left]{$\alpha=$};
\draw[->](7.5,1)--(9.5,1);
\end{scope}

\begin{scope}[shift={(10,0)}]	
\uvs{1,...,6}
\lvs{1,...,6}
\darcx23{.3}
\uarcx12{.3}
\uarcx23{.3}
\uarcx45{.3}
\stline11
\stline65
\draw(6.4,1)node[right]{$=\alpha^*$};
\end{scope}

\end{tikzpicture}
\caption{Diagrammatic representation, multiplication and involution in $\PP_6$.}
\label{fig:Pn}
\end{center}
\end{figure}

A key class of examples of regular $*$-semigroups comes from the so-called \emph{diagram monoids}.  This class includes important families such as partition monoids $\PP_n$, Brauer monoids $\B_n$, Temperley--Lieb monoids $\TL_n$, partial Brauer monoids $\PB_n$ and Motzkin monoids $\M_n$.  The elements of $\PP_n$ are the set partitions of $\{1,\ldots,n\}\cup\{1',\ldots,n'\}$.  These partitions are represented and multiplied diagrammatically, as in Figure \ref{fig:Pn}; for formal definitions and an extended discussion see \cite{EPA2024}.  The involution $\al\mt\al^*$ corresponds to a vertical reflection, as also shown in Figure \ref{fig:Pn}.  
The elements of $\B_n$ (resp.~$\PB_n$) are the partitions whose blocks have size $2$ (resp.~$\leq2$), while~$\TL_n$ (resp.~$\M_n$) consists of partitions from $\B_n$ (resp.~$\PB_n$) that can be drawn in planar fashion within the rectangle bounded by the vertices.  Thus, in Figure \ref{fig:Pn} we have $\be\in\TL_6$ and $\al\be\in\M_6$.  

The Temperley--Lieb monoid $\TL_n$ will be considered in Section \ref{sect:TLn} as a major example.  To study it we will require the following well-known presentation by generators and relations; see~\cite{BDP2002,East2021} for proofs.

\begin{thm}\label{thm:TLn}
The Temperley--Lieb monoid has monoid presentation
\[
\TL_n\cong\pres{X_T}{R_T},
\]
where $X_T = \{t_1,\ldots,t_{n-1}\}$, and where $R_T$ is the set of relations
\begin{align}
\tag*{\textsf{(T1)}} \label{T1} t_i^2 &= t_i &&\hspace{-2cm}\text{for all $i$}\\
\tag*{\textsf{(T2)}} \label{T2} t_it_j &= t_jt_i &&\hspace{-2cm}\text{if $|i-j|>1$}\\
\tag*{\textsf{(T3)}} \label{T3} t_it_jt_i &= t_i &&\hspace{-2cm}\text{if $|i-j|=1$.}\\
&&& \epfreseq
\end{align}
\end{thm}

In the above presentation, the generator $t_i$ corresponds to the diagram 
\begin{equation}\label{eq:taui}
\begin{tikzpicture}[scale=.53]
\uvs{1,3,4,5,6,8}
\lvs{1,3,4,5,6,8}
\uarc45
\darc45
\stline11
\stline33
\stline66
\stline88
\udotted13
\udotted68
\ldotted13
\ldotted68
\draw(0.6,1)node[left]{$\tau_i=$};
\node[above] ()at(1,2.1){\footnotesize$1$};
\node[above] ()at(4,2.1){\footnotesize$i$};
\node[above] ()at(5,2.05){\footnotesize$i\!\!+\!\!1$};
\node[above] ()at(8,2.1){\footnotesize$n$};
\node() at (8.75,.5){.};
\end{tikzpicture}
\end{equation}

\section{Projection algebras and chained projection groupoids}\label{sect:CPG}

In this section we give an overview of the constructions and results we will need from \cite{EPA2024}, building up to the isomorphism between the categories of regular $*$-semigroups and chained projection groupoids; see Theorem \ref{thm:iso}.  The presentation here is necessarily streamlined; for more details, and for proofs of the various assertions, see \cite{EPA2024}.

\subsection{Preliminaries on small categories}\label{subsect:C}

Categories considered in the paper come in two kinds: large categories whose objects and morphisms are algebraic structures and structure-preserving mappings; and small categories, which are thought of as algebraic objects in their own right.  The former are treated in the usual way; see for example \cite{Awodey2010,MacLane1998}.  The current section explains how we view the latter, following \cite{EPA2024}.

A small category $\CC$ will be identified with its morphism set.  The objects of $\CC$ are identified with the identities, the set of which is denoted $v\CC$.\footnote{The choice of notation alludes to the graph-theoretic interpretation of objects and morphisms as vertices and edges in a digraph. This viewpoint will become prominent in the final section of the paper.}
The domain and range maps are denoted $\bd,\br:\CC\to v\CC$, and we compose morphisms left-to-right, so that $a\circ b$ is defined when $\br(a)=\bd(b)$, and then $\bd(a\circ b)=\bd(a)$ and $\br(a\circ b)=\br(b)$.  For $p,q\in v\CC$ we write
\[
\CC(p,q)=\set{a\in\CC}{\bd(a)=p,\ \br(a)=q}
\]
for the set of all morphisms $p\to q$.

A \emph{$*$-category} is a small category $\CC$ with an involution, i.e.~a map ${\CC\to\CC:a\mt a^*}$ satisfying the following, for all $a,b\in\CC$:
\bit
\item $\bd(a^*)=\br(a)$, $\br(a^*)=\bd(a)$ and $(a^*)^*=a$.
\item If $\br(a)=\bd(b)$, then $(a\circ b)^*=b^*\circ a^*$.
\eit
A \emph{groupoid} is a $*$-category for which we additionally have $a\circ a^*=\bd(a)$ (and hence also $a^*\circ a=\br(a)$) for all $a\in\CC$.
In a groupoid, we typically write $a^*=a^{-1}$ for $a\in\CC$.

\newpage

An \emph{ordered $*$-category} (respectively, \emph{ordered groupoid}) is a $*$-category (respectively, groupoid)~$\CC$ equipped with a partial order $\leq$ satisfying the following, for all $a,b,c,d\in\CC$ and $p\in v\CC$: 
\bit
\item If $a\leq b$, then $\bd(a)\leq\bd(b)$, $\br(a)\leq\br(b)$ and $a^*\leq b^*$.
\item If $a\leq b$ and $c\leq d$, and if $\br(a)=\bd(c)$ and $\br(b)=\bd(d)$, then $a\circ c\leq b\circ d$.
\item For all $p\leq\bd(a)$ and $q\leq\br(a)$, there exist unique $u,v\leq a$ with $\bd(u)=p$ and $\br(v)=q$.  These elements are denoted $u={}_p\corest a$ and $v=a\rest_q$, and are called the \emph{left} and \emph{right restrictions} of $a$ to $p$ and $q$, respectively.
\eit

A \emph{congruence} on a small category $\CC$ is an equivalence relation $\approx$ on $\CC$ satisfying the following, for all $a,b,u,v\in\CC$:
\bit
\item $a\approx b \implies [\bd(a)=\bd(b)$ and $\br(a)=\br(b)]$,
\item $a\approx b \implies [u\circ a\approx u\circ b$ and $a\circ v\approx b\circ v]$, whenever the stated compositions are defined.
\eit
For a subset $\Om\sub\CC\times\CC$ with $\bd(a)=\bd(b)$ and $\br(a)=\br(b)$ for all $(a,b)\in\Om$, we write $\Om^\sharp$ for the congruence on $\CC$ generated by $\Om$.

An \emph{ordered $*$-congruence} on an ordered $*$-category $\CC$ is a congruence $\approx$ satisfying the following, for all $a,b\in\CC$  and $p\in v\CC$:
\bit
\item $a\approx b \implies a^*\approx b^*$,
\item $[a\approx b$ and $p\leq\bd(a)] \implies {}_p\corest a \approx {}_p\corest b$.
\eit
When ${\approx}=\Om^\sharp$, these two conditions can be verified by showing that they hold for all $(a,b)\in\Om$; see \cite[Lemma 2.7]{EPA2024}.

Given an (ordered $*$-) congruence $\approx$ on an (ordered $*$-) category $\CC$, the quotient $\CC/{\approx}$ is an (ordered $*$-) category.  Identifying an object $p\in v\CC$ with its $\approx$-class, we have ${v(\CC/{\approx}) = v\CC}$.

\subsection{Projection algebras}\label{subsect:PA}

The properties of projections of regular $*$-semigroups are formalised in what are now known as projection algebras, going back to Imaoka \cite{Imaoka1983}, who called them `$P$-groupoids'.  Here we recall the definition of these algebras, and list some known results that will be used later in the paper.

A \emph{projection algebra} is a unary algebra $P$, with set of operations $\set{\th_p}{p\in P}$ in one-one correspondence with $P$, satisfying the following axioms, for all $p,q\in P$:
\begin{enumerate}[label=\textup{\textsf{(P\arabic*)}},leftmargin=10mm]\bmc3
\item \label{P1} $p\theta_p = p$,
\item \label{P2} $\theta_p\theta_p=\theta_p$,
\item \label{P3} $p\theta_q\theta_p=q\theta_p$,
\item \label{P4} $\theta_p\theta_q\theta_p=\theta_{q\theta_p}$,
\item \label{P5} $\theta_p\theta_q\theta_p\theta_q=\theta_p\theta_q$.
\emc\end{enumerate}
The elements of a projection algebra are called \emph{projections}.  We write $\PA$ for the category of projection algebras with projection algebra morphisms, defined as maps $\phi:P\to P'$ satisfying
\begin{equation}
\label{eq:PAmorph}
(p\th_q)\phi = (p\phi) \th_{q\phi} \qquad\text{for all $p,q\in P$.}
\end{equation}
Projection algebras can also be thought of as binary algebras, with a single operation $\diamond$ given by $q\diamond p=q\th_p$.  In this formulation, projection algebra morphisms are simply maps $\phi:P\to P'$ satisfying $(q\diamond p)\phi=(q\phi)\diamond(p\phi)$ for $p,q\in P$.  The binary approach was used in \cite{Jones2012}, and compared in detail to the unary approach in \cite[Remark 4.2]{EPA2024}.

Given a regular $*$-semigroup $S$, the set of projections
\[
P = \bPP(S) = \set{p\in S}{p^2=p=p^*}
\]
becomes a projection algebra, with unary operations $\th_p$ as in \eqref{eq:thp}.  Conversely, any projection algebra is the projection algebra of some regular $*$-semigroup.  This latter fact was proved in \cite{Imaoka1983,Jones2012,EPA2024}, and it also follows from Theorem \ref{thm:C_P} below; see also \cite{NP1985}.  The assignment $S\mt \bPP(S)$ is the object part of a (forgetful) functor
\begin{equation}\label{eq:Pfunctor}
\bPP:\RSS\to\PA.
\end{equation}
Given a $*$-morphism $\phi:S\to S'$, the projection algebra morphism $\bPP(\phi):\bPP(S)\to \bPP(S')$ is simply the restriction $\bPP(\phi)=\phi|_{\bPP(S)}$.
It is easy to check that $\bPP(\phi)$ is indeed a morphism as defined in \eqref{eq:PAmorph}.

A projection algebra $P$ has three associated relations, $\leq$, $\leqF$ and $\F$, defined for $p,q\in P$ by
\begin{equation}\label{eq:rels}
p\leq q \iff p=p\th_q \COMMA p\leqF q \iff p=q\th_p \AND p\F q\iff[p\leqF q \text{ and } q\leqF p].
\end{equation}
The relation $\leq$ is a partial order, and we have $p\leq q \iff p=r\th_q$ for some $r\in P$.  The relation~$\leqF$ is reflexive, and $\F$ is reflexive and symmetric; neither~$\leqF$ nor~$\F$ is transitive in general.

We now list some important properties of projection algebras proved in \cite[Section 4]{EPA2024}.  Specifically, for any projection algebra $P$, and for any $p,q,r\in P$:  
\begin{enumerate}[label=\textup{\textsf{(PA\arabic*)}},leftmargin=13mm]
\item \label{PA1} $p\theta_q\F q\theta_p$,
\item \label{PA2} $[p\leq q\leqF r$ or $p\leqF q\leq r] \implies p\leqF r$,
\item \label{PA3} $p\leq q\implies p\leqF q$,
\item \label{PA4} $p\leq q\implies\th_p=\th_p\th_q = \th_q\th_p$,
\item \label{PA5} $p\leqF q \implies \th_p=\th_p\th_q\th_p$.
\end{enumerate}

We will also need the following simple result:

\begin{lemma}\label{lem:P}
For any $p_1,\ldots,p_k,q,r\in P$ we have
\[
\th_{q\th_{p_1}\cdots\th_{p_k}} = \th_{p_k}\cdots\th_{p_1}\th_q\th_{p_1}\cdots\th_{p_k} \AND r\th_{p_k}\cdots\th_{p_1}\th_q\th_{p_1}\cdots\th_{p_k}\th_r = q\th_{p_1}\cdots\th_{p_k}\th_r.
\]
\end{lemma}

\pf
The first claim follows by iterating \ref{P4}.  The second follows by applying the first, and then~\ref{P3}:
\[
r\th_{p_k}\cdots\th_{p_1}\th_q\th_{p_1}\cdots\th_{p_k}\th_r = r \th_{q\th_{p_1}\cdots\th_{p_k}}\th_r = q\th_{p_1}\cdots\th_{p_k}\th_r.  \qedhere
\]
\epf

\subsection{Path categories and chain groupoids}\label{subsect:CP}

Let $P$ be a projection algebra.  A \emph{($P$-)path}  is a path in the graph of the $\F$-relation, i.e.~a tuple
\[
\p = (p_1,p_2,\ldots,p_k)\in P^k \qquad\text{for some $k\geq1$, such that $p_1\F p_2\F\cdots\F p_k$.\footnotemark}
\]
\footnotetext{We allow repeated vertices in a path, in alignment with other graph-based algebraic structures such as Leavitt path algebras \cite{AAS2017} or graph inverse semigroups \cite{AH1975}.  Some authors would use the term `walk' for our paths.}%
We write $\bd(\p)=p_1$ and $\br(\p)=p_k$.  We identify each $p\in P$ with the path $(p)$ of length $1$.  
The \emph{path category} of $P$ is the ordered $*$-category $\P=\P(P)$ of all $P$-paths, with:
\newpage
\bit
\item object set $v\P=P$, 
\item composition $(p_1,\ldots,p_k)\circ(p_k,\ldots,p_l) = (p_1,\ldots,p_k,\ldots,p_l)$,
\item involution $(p_1,\ldots,p_k)^\rev = (p_k,\ldots,p_1)$,
\item restrictions ${}_q\corest(p_1,\ldots,p_k)=(q_1,\ldots,q_k)$ and $(p_1,\ldots,p_k)\rest_r=(r_1,\ldots,r_k)$, for $q\leq p_1$ and $r\leq p_k$, where 
\begin{equation}\label{eq:rest}
q_i = q\th_{p_2}\cdots\th_{p_i} = q\th_{p_1}\cdots\th_{p_i} \AND r_i = r\th_{p_{k-1}}\cdots\th_{p_i}= r\th_{p_k}\cdots\th_{p_i} \qquad\text{for $1\leq i\leq k$.}
\end{equation}
\eit

Given a projection algebra $P$, we write $\Om=\Om(P)$ for the set of all pairs $(\s,\t)\in\P\times\P$ of the following two forms:
\begin{enumerate}[label=\textup{\textsf{($\mathsf{\Om}$\arabic*)}},leftmargin=10mm]
\item \label{Om1} $\s=(p,p)$ and $\t=(p)$, for some $p\in P$,
\item \label{Om2} $\s=(p,q,p)$ and $\t=(p)$, for some $(p,q)\in {\F}$,
\end{enumerate}
and we write ${\approx}=\Om^\sharp$ for the congruence on $\P$ generated by $\Om$.  This is an ordered $*$-congruence, and the quotient is a groupoid, the \emph{chain groupoid} of $P$:
\[
\C = \C(P) = \P/{\approx}.
\]
The elements of $\C$ are called \emph{($P$-)chains}, and we denote by $[\p]$ the chain containing the path~${\p\in\P}$.
The order in $\C$ is determined by the restrictions, which are canonically inherited from $\P$.  Specifically, for $\p\in\P$ we have
\[
{}_q\corest[\p] = [{}_q\corest\p] \AND [\p]\rest_r = [\p\rest_r] \qquad\text{for $q\leq\bd(\p)=\bd[\p]$ and $r\leq\br(\p)=\br[\p]$.}
\]
These are well defined because $\approx$ is an \emph{ordered} congruence.

\begin{rem}\label{rem:conf}
Any chain $\c=[\p]\in\C$ can be uniquely represented in the form $[p_1,\ldots,p_k]$, where each $p_i$ is distinct from $p_{i+1}$ (if $i\leq k-1$) and from $p_{i+2}$ (if $i\leq k-2$).
This `reduced form' can be found by successively reducing $\p$, using the rules
\[
(p,p) \to (p) \ \ \text{for $p\in P$} \AND
(p,q,p) \to (p) \ \ \text{for $(p,q)\in{\F}$}.
\]
One way to establish uniqueness is to show that the rewriting system $(\P,{\to})$ is (locally) confluent and Noetherian, in the sense of \cite{Gerard1980}.  
\end{rem}

For any projection algebra morphism $\phi:P\to P'$, there is a well-defined ordered groupoid morphism
\begin{equation}\label{eq:Cphi}
\C(\phi):\C(P)\to\C(P') \GIVENBY [p_1,\ldots,p_k]\C(\phi) = [p_1\phi,\ldots,p_k\phi].
\end{equation}
In this way, $\C$ can be viewed as a functor from the category $\PA$ of projection algebras to the category $\OG$ of ordered groupoids.

\subsection{Chained projection groupoids}\label{subsect:CPG}

A \emph{weak projection groupoid} is a pair $(P,\G)$, consisting of an ordered groupoid~$\G$ whose object set~$P$ has the structure of a projection algebra, and for which the restriction to $P$ of the order on $\G$ 
coincides with the projection algebra order~$\leq$ from \eqref{eq:rels}.  For any $a\in\G$, we have a pair of maps
\[
\vt_a:\bd(a)^\da\to\br(a)^\da:p\mt \br({}_p\corest a) \AND \Th_a= \th_{\bd(a)}\vt_a:P\to\br(a)^\da,
\]
where here $q^\da=\set{p\in P}{p\leq q}$ is the down-set of $q\in P$.  
%
%
%
A \emph{projection groupoid} is a weak projection groupoid $(P,\G)$ for which:
\begin{enumerate}[label=\textup{\textsf{(G1)}},leftmargin=10mm]
\item \label{G1} $\th_{p\Th_a} = \Th_{a^{-1}}\th_p\Th_a$ for all $p\in P$ and $a\in\G$.
\end{enumerate}

Let $(P,\G)$ be a projection groupoid.  An \emph{evaluation map} is an ordered $v$-functor $\ve:\C(P)\to\G$, meaning that the following hold:
\[
\ve(p)=p \text{ for $p\in P$} \COMMa
\ve(\c\circ\d)=\ve(\c)\circ\ve(\d) \text{ if $\br(\c)=\bd(\d)$} \ANd
\ve({}_p\corest\c) = {}_p\corest(\ve(\c)) \text{ for $p\leq\bd(\c)$.}
\]
We note in passing that evaluation maps are written to the left of their arguments, as we feel they are easier to read this way (and never need to be composed).  
Since $\C(P)$ is generated by chains of length~$2$, the functor $\ve$ is completely determined by the elements $\ve[p,q]$ for $(p,q)\in{\F}$.  These elements feature in the remaining assumptions and constructions involving evaluation maps.

Given a morphism $b\in\G$, a pair of projections $(e,f)\in P\times P$ is said to be \emph{$b$-linked} if 
\[
f = e\Th_b\th_f \AND e = f\Th_{b^{-1}}\th_e.
\]
Given such a $b$-linked pair $(e,f)$, and writing $q=\bd(b)$ and $r=\br(b)$, we define further projections
\[
e_1 = e\th_q \COMMA e_2 = f\Th_{b^{-1}} \COMMA f_1 = e\Th_b \AND f_2 = f\th_r.
\]
The groupoid $\G$ then contains two well-defined morphisms:
\begin{align*}
\lam(e,b,f) &= \ve[e,e_1]\circ {}_{e_1}\corest b \circ \ve[f_1,f] &\text{and}&& \rho(e,b,f) &= \ve[e,e_2]\circ {}_{e_2}\corest b \circ \ve[f_2,f]\\
&= \ve[e,e_1]\circ b\rest_{f_1} \circ \ve[f_1,f] &&&  &= \ve[e,e_2]\circ b\rest_{f_2} \circ \ve[f_2,f].
\end{align*}

A \emph{chained projection groupoid} is a triple $(P,\G,\ve)$, where $(P,\G)$ is a projection groupoid, and $\ve:\C(P)\to\G$ is an evaluation map for which:
\begin{enumerate}[label=\textup{\textsf{(G2)}},leftmargin=10mm]
\item \label{G2} $\lam(e,b,f) = \rho(e,b,f)$ for every $b\in\G$, and every $b$-linked pair $(e,f)$.
\end{enumerate}
We write $\CPG$ for the category of all chained projection groupoids with \emph{chained projection functors} as morphisms.  A chained projection functor $(P,\G,\ve)\to(P',\G',\ve')$ is an ordered groupoid functor $\phi:\G\to\G'$ such that
\bit
\item $v\phi=\phi|_P:P\to P'$ is a projection algebra morphism, and
\item $\phi$ respects evaluation maps, in the sense that the following diagram commutes:
\[
\begin{tikzcd}[sep=scriptsize]
\C(P) \arrow{rr}{\C(v\phi)} \arrow[swap]{dd}{\ve} & ~ & \C(P') \arrow{dd}{\ve'} \\%
~&~&~\\
\G \arrow{rr}{\phi}& ~ & \G',
\end{tikzcd}
\]
where $\C(v\phi)$ is constructed from $v\phi$ as in \eqref{eq:Cphi}.  Explicitly, this is to say that
\[
(\ve[p_1,\ldots,p_k])\phi = \ve'[p_1\phi,\ldots,p_k\phi]
\]
whenever $p_1,\ldots,p_k\in P$ and $p_1\F\cdots\F p_k$.
\eit

\subsection{A category isomorphism}\label{subsect:iso}

The main result of \cite{EPA2024} is that the categories $\RSS$ and $\CPG$, of regular $*$-semigroups and chained projection groupoids, are isomorphic.  The proof involves two functors, $\bG$ and $\bS$, between the two categories, which operate as follows.

A regular $*$-semigroup $S$ determines a chained projection groupoid $\bG(S) = (P,\G,\ve)$, where:
\bit
\item $P$ is the projection algebra of $S$.
\item $\G$ is an ordered groupoid built from $S$ as follows.  The morphisms are the elements of $S$, and the objects/identities are the projections, with $\bd(a)=aa^*$ and $\br(a)=a^*a$ for $a\in S$.  The composition and involution are given by $a\circ b = ab$ when $\br(a)=\bd(b)$, and $a^{-1}=a^*$.  Restrictions are given by ${}_p\corest a=pa$ and $a\rest_q=aq$ for $p\leq\bd(a)$ and $q\leq\br(a)$.
\item $\ve:\C(P)\to\G$ is the evaluation map given by $\ve[p_1,\ldots,p_k]=p_1\cdots p_k$, where this product is taken in $S$.
\eit
Any $*$-morphism $S\to S'$ is also a chained projection functor $\bG(S)\to\bG(S')$.

Conversely, any chained projection groupoid $(P,\G,\ve)$ gives rise to a regular $*$-semigroup $\bS(P,\G,\ve)$, with underlying set $\G$, and:
\bit
\item involution given by $a^*=a^{-1}$, 
\item product defined, for $a,b\in\G$ with $\br(a)=p$ and $\bd(b)=q$, by
\begin{equation}\label{eq:pr}
a\pr b = a\rest_{p'} \circ \ve[p',q'] \circ {}_{q'}\corest b, \WHERE p'=q\th_p \ANd q'=p\th_q.
\end{equation}
\eit
Any chained projection functor $(P,\G,\ve)\to(P',\G',\ve')$ is also a $*$-morphism ${\bS(P,\G,\ve)\to\bS(P',\G',\ve')}$.

\begin{thm}[{see \cite[Theorem 8.1]{EPA2024}}]
\label{thm:iso}
$\bG$ and $\bS$ are mutually inverse isomorphisms between the categories $\RSS$ and $\CPG$.  \epfres
\end{thm}

\section{Construction of the chain semigroup}\label{sect:CP}

We now come to the main focus of our study:  the chain semigroup $\PG(P)$ associated to a projection algebra $P$.  This semigroup is built in Subsection \ref{subsect:C_P} by first constructing a chained projection groupoid $(P,\ol\C,\nu)$, and then applying the functor ${\bS:\CPG\to\RSS}$ from Theorem~\ref{thm:iso}.  Here~$\ol\C$ is a homomorphic image of the chain groupoid $\C$ (see Subsection~\ref{subsect:olC}), and $\nu$ is the quotient map $\C\to\ol\C$.  The definition of $\ol\C$ involves the notion of \emph{linked pairs} of projections, which are introduced in Subsection \ref{subsect:LP}.

\subsection{Linked pairs of projections}\label{subsect:LP}

\begin{defn}\label{defn:CP_P}
Let $P$ be a projection algebra, and let $p\in P$.  A pair of projections $(e,f)\in P\times P$ is said to be \emph{$p$-linked} if
\begin{equation}\label{eq:epf}
f = e\th_p\th_f \AND e = f\th_p\th_e.
\end{equation}
Associated to such a $p$-linked pair $(e,f)$ we define the tuples
\[
\lam(e,p,f) = (e,e\th_p,f) \AND \rho(e,p,f) = (e,f\th_p,f).
\]
\end{defn}

The next two results gather some important basic properties of $p$-linked pairs.  Recall that $\P=\P(P)$ denotes the path category of $P$.

\begin{lemma}\label{lem:epf}
If $(e,f)$ is $p$-linked, then
\ben
\item \label{epf1} $(f,e)$ is also $p$-linked, and we have $\lam(e,p,f)^\rev = \rho(f,p,e)$ and $\rho(e,p,f)^\rev = \lam(f,p,e)$,
\item \label{epf2} $e,f\leqF p$,
\item \label{epf3} Each of $e$ and $f$ is $\F$-related to each of $e\th_p$ and $f\th_p$; consequently both $\lam(e,p,f)$ and $\rho(e,p,f)$ belong to $\P(e,f)$.
\een
\end{lemma}

\pf
\firstpfitem{\ref{epf1}}  This follows directly by inspecting Definition~\ref{defn:CP_P}.

\pfitem{\ref{epf2}}  By the symmetry afforded by part \ref{epf1}, it suffices to show that $e\leqF p$.  For this we use~\eqref{eq:epf}, Lemma \ref{lem:P} and \ref{P3} to calculate
\[
p\th_e = p\th_{f\th_p\th_e} = p\th_e\th_p\th_f\th_p\th_e = e\th_p\th_f\th_p\th_e = f\th_p\th_e = e.
\]

\pfitem{\ref{epf3}}  By symmetry, it suffices to show that $e \F e\th_p \F f$.  Combining ${e\leqF p}$ with \ref{PA1}, it follows that $e=p\th_e\F e\th_p$.  We obtain $f \leqF e\th_p$ directly from \eqref{eq:epf}.  Using \ref{P4} and \eqref{eq:epf} we calculate $f\th_{e\th_p} = f\th_p\th_e\th_p = e\th_p$, so that $e\th_p \leqF f$.
\epf

\begin{rem}
Lemma \ref{lem:epf}\ref{epf3} says that $(e,f)$ being $p$-linked implies $e,f \F e\th_p,f\th_p$.  The converse of this holds as well, as \eqref{eq:epf} says that $f\leqF e\th_p$ and $e\leqF f\th_p$.  Thus, we could take ${e,f \F e\th_p,f\th_p}$ as an equivalent definition for $(e,f)$ to be $p$-linked.
\end{rem}

\begin{figure}[h]
\begin{center}
\begin{tikzpicture}[scale=0.5]
\tikzstyle{vertex}=[circle,draw=black, fill=white, inner sep = 0.07cm]
\nc\sss{4}
\nc\ttt{1}
%
\node (e) at (-1,0.5){$e$};
\node (p) at (4,7){$p$};
\node (e') at (3,2){$e\th_p$};
\node (f') at (5,-\ttt){$f\th_p$};
\begin{scope}[shift={(\sss,0)}]
\node (f) at (5,0.5){$f$};
\end{scope}
\draw[->-=0.6] (e)--(e');
\draw[->-=0.5] (e')--(f);
\draw[->-=0.55] (e)--(f');
\draw[->-=0.5] (f')--(f);
\draw[white,line width=2mm] (p)--(f');
\draw[dashed] (e')--(p)--(f');
\draw[dotted] (e)--(p) (f)--(p);
\end{tikzpicture}
\caption{A projection $p\in P$, and a $p$-linked pair $(e,f)$, as in Definition \ref{defn:CP_P}.  Dotted and dashed lines indicate $\leqF$ and~$\leq$ relationships, respectively.  
Lines with arrows indicate the paths $\lam(e,p,f)$ and~$\rho(e,p,f)$. Of course, by Lemma \ref{lem:epf}\ref{epf1}, the dual diagram in which all four arrows are reversed is also valid.
See Remark \ref{rem:LPP} for more details.}
\label{fig:plinked}
\end{center}
\end{figure}
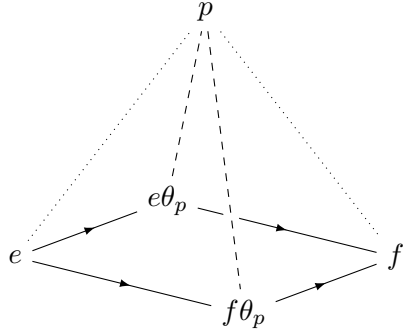

\begin{rem}\label{rem:LPP}
Consider a projection $p\in P$, and a $p$-linked pair $(e,f)$.  By Lemma \ref{lem:epf}\ref{epf2} we have $e,f\leqF p$, and of course we also have $e\th_p,f\th_p\leq p$.  These relationships are all shown in Figure~\ref{fig:plinked}.  In the diagram, each arrow $s\to t$ stands for the $P$-path $(s,t)\in\P$.  Thus, the upper and lower paths from $e$ to $f$ 
correspond to $\lam(e,p,f)$ and~$\rho(e,p,f)$, respectively.
\end{rem}

\begin{lemma}\label{lem:e'pf'}
If $(e,f)$ is $p$-linked, and if $e'\leq e$, then $(e',f')$ is $p$-linked, where $f'=e'\th_p\th_f$, and we have
\[
{}_{e'}\corest\lam(e,p,f) = \lam(e',p,f') \AND {}_{e'}\corest\rho(e,p,f) = \rho(e',p,f').
\]
\end{lemma}

\pf
To show that $(e',f')$ is $p$-linked, we must show that $f'=e'\th_p\th_{f'}$ and $e'=f'\th_p\th_{e'}$.  For the first we use the definition of $f'$ and Lemma \ref{lem:P} several times to calculate
\[
e'\th_p\th_{f'} = e'\th_p\th_{e'\th_p\th_f} 
= e'\th_p\th_f\th_p\th_{e'}\th_p\th_f 
= f\th_p\th_{e'}\th_p\th_f 
= {e'}\th_p\th_f  = f'.
\]
For the second we have
\begin{align*}
f'\th_p\th_{e'} = e'\th_p\th_f\th_p\th_{e'}
&= f\th_p\th_{e'} &&\text{by definition, and by Lemma \ref{lem:P}}\\
&= f\th_p\th_e\th_{e'} &&\text{by \ref{PA4}, as $e'\leq e$}\\
&= e\th_{e'} &&\text{by \eqref{eq:epf}}\\
&= e' &&\text{as $e'\leqF e$ by \ref{PA3}.}
\end{align*}
We now show that ${}_{e'}\corest\lam(e,p,f) = \lam(e',p,f')$, the proof that ${}_{e'}\corest\rho(e,p,f) = \rho(e',p,f')$ being analogous.  Using \eqref{eq:rest} and Definition \ref{defn:CP_P}, we have
\[
{}_{e'}\corest\lam(e,p,f) = (e',e'\th_{e\th_p},e'\th_{e\th_p}\th_f) \AND \lam(e',p,f') = (e',e'\th_p,f').
\]
Both have first component $e'$.  We deduce equality of the second and third components from the following calculations:
\bit
\item $e'\th_{e\th_p} = e'\th_p\th_e\th_p = e'\th_e\th_p\th_e\th_p = e'\th_e\th_p = e'\th_p$, using $e'\leq e$, \ref{P4} and \ref{P5},
\item $(e'\th_{e\th_p})\th_f = e'\th_p\th_f = f'$, using the previous calculation, and the definition of $f'$.  \qedhere
\eit
\epf

\subsection{The reduced chain groupoid}\label{subsect:olC}

We now use linked pairs to define the reduced chain groupoid of a projection algebra.

\begin{defn}\label{defn:olC}
Let $P$ be a projection algebra, and let $\Xi=\Xi(P)$ be the set of all pairs ${(\s,\t)\in\P\times\P}$ of the forms \ref{Om1} and \ref{Om2}, as well as:
\begin{enumerate}[label=\textup{\textsf{($\mathsf{\Om}$\arabic*)}},leftmargin=10mm]\addtocounter{enumi}{2}
\item \label{Om3} $\s=\lam(e,p,f)$ and $\t=\rho(e,p,f)$, for some $p\in P$, and some $p$-linked pair $(e,f)$.
\end{enumerate}
Let ${\ssim}=\Xi^\sharp$ be the congruence on $\P$ generated by $\Xi$, and define the quotient
\[
\ol\C = \ol\C(P) = \P/{\ssim}.
\]
As $\Om\sub\Xi$, we have ${\approx}\sub{\ssim}$, and so $\ol\C$ is a quotient of the chain groupoid $\C=\C(P)=\P/{\approx}$.  As such $\ol\C$ is itself a groupoid, which we call the \emph{reduced chain groupoid of~$P$}.  The elements of~$\ol\C$ are called \emph{reduced ($P$-)chains}, and we denote by $\ldb\p\rdb$ the reduced chain containing the path~${\p\in\P}$.  We denote the quotient map $\C\to\ol\C$ by
\[
\nu:\C\to\ol\C , \qquad\text{which is given by}\qquad \nu[\p]=\ldb\p\rdb \qquad\text{for $\p\in\P$.}
\]
\end{defn}

\begin{rem}
We will see in Subsection \ref{subsect:degen} that there exists a natural (and generally smaller) subset of $\Xi$ that also generates the congruence $\ssim$.  For now it is more convenient to use $\Xi$, as its definition is more symmetrical.
\end{rem}

In Proposition \ref{prop:PCnu} we show that $(P,\ol\C,\nu)$ is a chained projection groupoid.  We build towards this with a number of lemmas.

\begin{lemma}\label{lem:OmP}
$\ssim$ is an ordered $*$-congruence, and $\ol\C$ is an ordered groupoid.
\end{lemma}

\pf
We have already observed that $\ol\C$ is a groupoid.  
As explained in Section \ref{subsect:C}, we can show that ${\ssim}=\Xi^\sharp$ is an ordered $*$-congruence by showing that 
\[
\s^\rev\ssim\t^\rev \AND {}_{p}\corest\s \ssim {}_{p}\corest\t \qquad\text{for all $(\s,\t)\in\Xi$, and all $p\leq\bd(\s)$.}
\]
When $(\s,\t)\in\Xi$ has type \ref{Om1} or \ref{Om2}, this was done in \cite[Lemma 5.16]{EPA2024}.  For pairs of type~\ref{Om3} we apply Lemmas \ref{lem:epf} and~\ref{lem:e'pf'}.
\epf

The order in $\ol\C$ is determined by the restrictions, which are inherited from $\C$ (and hence ultimately from $\P$).  That is, for $\p\in\P$ they are given by
\[
{}_q\corest\ldb\p\rdb = \ldb{}_q\corest\p\rdb \AND \ldb\p\rdb\rest_r = \ldb\p\rest_r\rdb \qquad\text{for $q\leq\bd(\p)=\bd\ldb\p\rdb$ and $r\leq\br(\p)=\br\ldb\p\rdb$.}
\]
These are well defined by Lemma \ref{lem:OmP}

It is clear that $(P,\ol\C)$ is a weak projection groupoid.  To show that it is a projection groupoid we need to verify~\ref{G1}.  To do so, we need to understand the maps $\Th_\c=\th_{\bd(\c)}\vt_\c$.

\begin{lemma}\label{lem:Thc}
For $\c=\ldb p_1,\ldots,p_k\rdb\in\ol\C$, we have $\Th_\c = \th_{p_1}\cdots\th_{p_k}$.
\end{lemma}

\pf
Using \eqref{eq:rest}, and writing $\p=(p_1,\ldots,p_k)\in\P$, we first calculate
\[
\label{eq:vtc} q\vt_\c = \br({}_q\corest\c) = \br({}_q\corest\ldb\p\rdb) = \br\ldb{}_q\corest\p\rdb = \br({}_q\corest\p) = q\th_{p_2}\cdots\th_{p_k} \qquad\text{for $q\leq p_1=\bd(\c)$.}
\]
It follows from this that $t\Th_\c = t\th_{\bd(\c)}\vt_\c = t\th_{p_1}\th_{p_2}\cdots\th_{p_k}$ for arbitrary $t\in P$.
\epf


\begin{lemma}
If $P$ is a projection algebra, then $(P,\ol\C)$ is a projection groupoid.
\end{lemma}

\pf
To verify \ref{G1}, consider a reduced chain $\c=\ldb p_1,\ldots,p_k\rdb\in\ol\C$, and let $q\in P$.  
Then by Lemmas \ref{lem:Thc} and \ref{lem:P} we have
%
\[
\th_{q\Th_\c} = \th_{q\th_{p_1}\cdots\th_{p_k}} = \th_{p_k}\cdots\th_{p_1}\th_q\th_{p_1}\cdots\th_{p_k} = \Th_{\c^{-1}}\th_q\Th_\c . \qedhere
\]
\epf

We now bring in the quotient map $\nu:\C\to\ol\C$ from Definition \ref{defn:olC}.


\begin{lemma}
If $P$ is a projection algebra, then $\nu$ is an evaluation map.
\end{lemma}

\pf
Clearly $\nu$ is a $v$-functor.  To see that it is ordered, consider a path $\p=(p_1,\ldots,p_k)\in\P$, and let $q\leq p_1$.  
Then with the $q_i$ as in \eqref{eq:rest}, we have
\[
\nu({}_q\corest[\p]) = \nu[q_1,\ldots,q_k] = \ldb q_1,\ldots,q_k\rdb = {}_q\corest\ldb p_1,\ldots,p_k\rdb = {}_q\corest(\nu[\p]).  \qedhere
\]
\epf

\begin{prop}\label{prop:PCnu}
If $P$ is a projection algebra, then $(P,\ol\C,\nu)$ is a chained projection groupoid.
\end{prop}

\begin{figure}[h]
\begin{center}
\scalebox{0.85}{
\begin{tikzpicture}[xscale=0.6*0.9,yscale=0.8*0.9]
\tikzstyle{vertex}=[circle,draw=black, fill=white, inner sep = 0.07cm]
\nc\sss{4}
\nc\ttt{1}
\nc\uuu{6}
\node[left] () at (-1.2,2){${}_{\phantom{0}}e=$};
\node (e) at (-1,2){$u_0$};
\node (f) at (1+5*\sss,-\ttt){$v_5$};
\node[right] () at (1+5*\sss+.2,-\ttt){$=f_{\phantom{0}}$};
\node (p1) at (4,\uuu){$p_1$};
\node (u1) at (3,2){$u_1$};
\node (v1) at (5,-\ttt){$v_1$};
\begin{scope}[shift={(\sss,0)}]
\node (p2) at (4,\uuu){$p_2$};
\node (u2) at (3,2){$u_2$};
\node (v2) at (5,-\ttt){$v_2$};
\end{scope}
\begin{scope}[shift={(2*\sss,0)}]
\node (p3) at (4,\uuu){$p_3$};
\node (u3) at (3,2){$u_3$};
\node (v3) at (5,-\ttt){$v_3$};
\end{scope}
\begin{scope}[shift={(3*\sss,0)}]
\node (p4) at (4,\uuu){$p_4$};
\node (u4) at (3,2){$u_4$};
\node (v4) at (5,-\ttt){$v_4$};
\end{scope}
\draw[->-=0.5] (p1)--(p2);
\draw[->-=0.5] (p2)--(p3);
\draw[->-=0.5] (p3)--(p4);
\draw[->-=0.6] (e)--(u1);
\draw[->-=0.66] (u1)--(u2);
\draw[->-=0.66] (u2)--(u3);
\draw[->-=0.66] (u3)--(u4);
\draw[->-=0.5] (u4)--(f);
\draw[->-=0.5] (e)--(v1);
\draw[->-=0.5] (v1)--(v2);
\draw[->-=0.5] (v2)--(v3);
\draw[->-=0.5] (v3)--(v4);
\draw[->-=0.5] (v4)--(f);
\draw[->-=0.5] (u1)--(v2);
\draw[->-=0.5] (u2)--(v3);
\draw[->-=0.5] (u3)--(v4);
\draw[white,line width=2mm] (p1)--(v1) (p2)--(v2) (p3)--(v3) (p4)--(v4);
\draw[dashed] (u1)--(p1)--(v1) (u2)--(p2)--(v2) (u3)--(p3)--(v3) (u4)--(p4)--(v4);
\end{tikzpicture}
}
\caption{The projections $e,f,p_i,u_i,v_i$ from the proof of Proposition \ref{prop:PCnu}, shown here in the case $k=4$.  Dashed lines indicate $\leq$ relationships.  Each arrow $s\to t$ represents the $P$-path $(s,t)\in\P$, so the upper and lower paths $e\to f$ represent $\lam(e,\c,f) = (e,u_1,\ldots,u_k,f)$ and $\rho(e,\c,f) = (e,v_1,\ldots,v_k,f)$, respectively.}
\label{fig:clinked}
\end{center}
\end{figure}
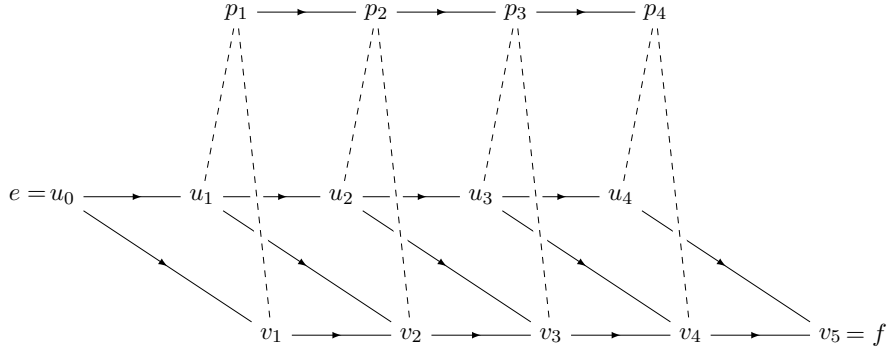

\pf
It remains to verify \ref{G2}.  To do so, fix a $\c$-linked pair $(e,f)$, where $\c=\ldb p_1,\ldots,p_k\rdb\in\ol\C$.  
So
\begin{equation}\label{eq:ecf}
f = e\Th_\c\th_f \AND e = f\Th_{\c^{-1}}\th_e,
\end{equation}
and we must show that $\lam(e,\c,f) = \rho(e,\c,f)$.  Keeping $\ve=\nu$ in mind, these morphisms are defined by
\begin{equation}\label{eq:lamecf}
\lam(e,\c,f) = \ldb e,e_1\rdb \circ {}_{e_1}\corest \c \circ \ldb f_1,f\rdb  \AND \rho(e,\c,f) = \ldb e,e_2\rdb \circ \c\rest_{f_2} \circ \ldb f_2,f\rdb ,
\end{equation}
in terms of the projections
\[
e_1 = e\th_{p_1} \COMMA e_2 = f\Th_{\c^{-1}} \COMMA f_1 = e\Th_\c \AND f_2 = f\th_{p_k}.
\]
For convenience, we will write
\[
{}_{e_1}\corest \c = \ldb  u_1,\ldots,u_k\rdb  \AND \c\rest_{f_2} = \ldb v_1,\ldots,v_k\rdb .
\]
Using \eqref{eq:rest}, and $e_1=e\th_{p_1}$ and $f_2 = f\th_{p_k}$, we have
\[
u_i = e_1\th_{p_2}\cdots\th_{p_i} = e\th_{p_1}\cdots\th_{p_i} \AND v_i = f_2\th_{p_{k-1}}\cdots\th_{p_i} = f\th_{p_k}\cdots\th_{p_i} ,
\]
for each $i$.  Keeping in mind that the compositions in \eqref{eq:lamecf} exist, we have
\[
\lam(e,\c,f) = \ldb e,u_1,\ldots,u_k,f\rdb  \AND \rho(e,\c,f) = \ldb e,v_1,\ldots,v_k,f\rdb ,
\]
and we must show that these are equal.  It will also be convenient to additionally write $u_0=e$ and $v_{k+1}=f$.  To assist with understanding the coming arguments, these projections are shown in Figure \ref{fig:clinked} (in the case $k=4$).  Our task is essentially to show that the large `rectangle' at the bottom of the diagram commutes, modulo $\ssim$.

We recall the notion of $p$-linked pairs (Definition \ref{defn:CP_P}) and claim that:
\begin{equation}
\label{eq:uvlinked}
(u_{i-1},v_{i+1}) \text{ is } p_i\text{-linked for each } 1\leq i\leq k.
\end{equation}
To prove this, we must show that
\[
v_{i+1} = u_{i-1}\th_{p_i}\th_{v_{i+1}} \AND u_{i-1} = v_{i+1}\th_{p_i}\th_{u_{i-1}}.
\]
First we note that \eqref{eq:ecf} and Lemma \ref{lem:Thc} give
\[
f = e\th_{p_1}\cdots\th_{p_k}\th_f \AND e = f\th_{p_k}\cdots\th_{p_1}\th_e.
\]
Combining this with Lemma \ref{lem:P}, we obtain
\begin{align*}
u_{i-1}\th_{p_i}\th_{v_{i+1}} &= e\th_{p_1}\cdots\th_{p_{i-1}} \cdot \th_{p_i} \cdot \th_{f\th_{p_k}\cdots\th_{p_{i+1}}} \\
&= e\th_{p_1}\cdots\th_{p_{i-1}} \cdot \th_{p_i} \cdot \th_{p_{i+1}}\cdots\th_{p_k}\th_f\th_{p_k}\cdots\th_{p_{i+1}} = f\th_{p_k}\cdots\th_{p_{i+1}} = v_{i+1}.
\end{align*}
The proof that $u_{i-1} = v_{i+1}\th_{p_i}\th_{u_{i-1}}$ is analogous, and \eqref{eq:uvlinked} is proved. 
Now, the linked pairs in~\eqref{eq:uvlinked} lead to the paths
\begin{align*}
\lam(u_{i-1},p_i,v_{i+1}) &= (u_{i-1},u_{i-1}\th_{p_i},v_{i+1}) &\text{and}&& \rho(u_{i-1},p_i,v_{i+1}) &= (u_{i-1},v_{i+1}\th_{p_i},v_{i+1})\\
&= (u_{i-1},u_i,v_{i+1}) &&&&=(u_{i-1},v_i,v_{i+1}),
\end{align*}
as in Definition \ref{defn:CP_P}.  Since $(\lam(u_{i-1},p_i,v_{i+1}),\rho(u_{i-1},p_i,v_{i+1}))$ is a pair of the form \ref{Om3}, we have
\begin{equation}\label{eq:uvv}
\ldb u_{i-1},v_i,v_{i+1}\rdb  = \ldb u_{i-1},u_i,v_{i+1}\rdb  \qquad\text{for each $1\leq i\leq k$.}
\end{equation}
In other words, each of the small rectangles at the bottom of Figure \ref{fig:clinked} commutes, and then it follows that the large rectangle commutes.  Formally, we repeatedly use \eqref{eq:uvv} in the indicated places, as follows:
\begin{align*}
\rho(e,\c,f) = \ldb e,v_1,\ldots,v_k,f\rdb  
&= \ldb \ul{u_0,v_1,v_2},v_3,v_4,\ldots,v_k,v_{k+1}\rdb  \\
&= \ldb u_0,\ul{u_1,v_2,v_3},v_4,\ldots,v_k,v_{k+1}\rdb  \\
&= \ldb u_0,u_1,\ul{u_2,v_3,v_4},\ldots,v_k,v_{k+1}\rdb  \\
& \hspace{1.3mm} \vdots\\
&= \ldb u_0,u_1,u_2,u_3,u_4,\ldots,u_k,v_{k+1}\rdb  = \ldb e,u_1,\ldots,u_k,f\rdb  = \lam(e,\c,f),
\end{align*}
and the proof is complete.
\epf

For a projection algebra $P$, we write $\bF(P)=(P,\ol\C,\nu)$ for the chained projection groupoid from Proposition \ref{prop:PCnu}.  The assignment $P\mt\bF(P)$ can be thought of as an object map ${\PA\to\CPG}$.  We can extend this to morphisms as well.  Indeed, fix a projection algebra morphism $\phi:P\to P'$.  In what follows we use the standard abbreviations for the constructions associated to $P$, and use dashes to distinguish those for $P'$; e.g.~$\P=\P(P)$ and $\C'=\C(P')$.  We first define a mapping
\begin{equation}\label{eq:3comp}
\varphi:\P\to\ol\C' \BY \p\varphi = \nu'([\p]\C(\phi)),\qquad \text{i.e.}\qquad (p_1,\ldots,p_k)\varphi = \ldb p_1\phi,\ldots,p_k\phi\rdb.
\end{equation}
Note that $\varphi$ is the composite of three ordered $*$-functors, namely the quotient map $\P\to\C$, followed by ${\C(\phi):\C\to\C'}$, and then the evaluation $\nu':\C'\to\ol\C'$.
Hence $\varphi$ is itself an ordered $*$-functor. It is straightforward to verify that $\Xi\sub\ker(\varphi)$, meaning that $\s\varphi=\t\varphi$ for all $(\s,\t)\in\Xi$.  Therefore there is a well-defined ordered groupoid functor
\begin{align*}
\bF(\phi):\ol\C\to\ol\C' \GIVENBY &\ldb\p\rdb\bF(\phi) = \p\varphi = \nu'([\p]\C(\phi)),\\[2mm]
\text{i.e.}\qquad &\ldb p_1,\ldots,p_k\rdb\bF(\phi) = \ldb p_1\phi,\ldots,p_k\phi\rdb.
\end{align*}

\begin{prop}\label{prop:Cfunctor}
$\bF$ is a functor $\PA\to\CPG$.
\end{prop}

\pf
We first verify that $\bF(\phi)$, as above, is a chained projection functor from ${\bF(P)=(P,\ol\C,\nu)}$ to $\bF(P')=(P',\ol\C',\nu')$.  The restriction of $\bF(\phi)$ to $P=v\ol\C$ is $\phi$, which is a projection algebra morphism by assumption.  To show that $\bF(\phi)$ preserves evaluation maps we need to check that the following diagram commutes:
\[
\begin{tikzcd}[sep=scriptsize]
\C(P) \arrow{rr}{\C(\phi)} \arrow[swap]{dd}{\nu} & ~ & \C(P')\arrow{dd}{\nu'} \\%
~&~&~\\
\ol\C(P) \arrow{rr}{\bF(\phi)}& ~ & \ol\C(P'),
\end{tikzcd}
\]
and this is routine.  It is also routine to check that $\bF(\phi\circ\phi')=\bF(\phi)\circ\bF(\phi')$ for composable morphisms $\phi$ and $\phi'$, and that $\bF(\id_P)=\id_{\bF(P)}$ for all $P$.
\epf

\subsection{The chain semigroup}\label{subsect:C_P}

As a result of Proposition \ref{prop:PCnu} and Theorem \ref{thm:iso}, we have a regular $*$-semigroup $\bS(P,\ol\C,\nu)$, which we will denote by $\PG(P)$, and call the \emph{chain semigroup of $P$}.  The choice of notation, which harkens back to the free idempotent-generated semigroups $\IG(E)$, will be justified by the results of Section \ref{sect:F}.  In line with Subsection \ref{subsect:iso}, we denote the product in $\PG(P) = \bS(P,\ol\C,\nu)$ by~$\pr$.  To give an explicit description of $\pr$, consider an arbitrary pair of reduced chains $\c=\ldb p_1,\ldots,p_k\rdb $ and $\d=\ldb q_1,\ldots,q_l\rdb $, and write $p=\br(\c)=p_k$ and $q=\bd(\d)=q_1$.  As in~\eqref{eq:pr}, and remembering that the evaluation map is $\nu$, we have
\[
\c\pr\d = \c\rest_{p'}\circ\ldb p',q'\rdb \circ{}_{q'}\corest\d \WHERE p'=q\th_p \AND q'=p\th_q.
\]
Using \eqref{eq:rest}, we have $\c\rest_{p'} = \ldb p_1',\ldots,p_k'\rdb $ and ${}_{q'}\corest\d = \ldb q_1',\ldots,q_l'\rdb $, where
\[
p_i' = p'\th_{p_k}\cdots\th_{p_i} \AND q_j' = q'\th_{q_1}\cdots\th_{q_j} \qquad\text{for $1\leq i\leq k$ and $1\leq j\leq l$,}
\]
and where $p'=p_k'$ and $q'=q_1'$.  Again, these restrictions are well defined because of Lemma~\ref{lem:OmP}.  It follows that
\[
\c\pr\d = \c\rest_{p'}\circ\ldb p',q'\rdb \circ{}_{q'}\corest\d = \ldb p_1',\ldots,p_k'\rdb \circ\ldb p_k',q_1'\rdb \circ\ldb q_1',\ldots,q_l'\rdb  = \ldb p_1',\ldots,p_k',q_1',\ldots,q_l'\rdb 
\]
is simply the \emph{concatenation} of $\c\rest_{p'}$ and ${}_{q'}\corest\d$, which we denote by $\c\rest_{p'} \op {}_{q'}\corest\d$.
As special cases we have
\begin{equation}\label{eq:spmult}
\text{$\c\pr\d = \c\op\d$ if $\br(\c) \F \bd(\d)$ \AND $\c\pr\d = \c\circ\d$ if $\br(\c) = \bd(\d)$.}
\end{equation}

\begin{defn}\label{defn:C_P}
The \emph{chain semigroup} $\PG(P)$ of a projection algebra $P$, is the regular $*$-semigroup defined as follows.
\begin{enumerate}[label=\textup{\textsf{(CP\arabic*)}},leftmargin=13mm]
\item \label{CP1} The elements of $\PG(P)$ are the reduced ($P$-)chains, $\ldb p_1,\ldots,p_k\rdb $, as in Definition \ref{defn:olC}.
\item \label{CP2} The product $\pr$ in $\PG(P)$ is defined, for $\c,\d\in\PG(P)$ with $p=\br(\c)$ and $q=\bd(\d)$, by
\[
\c\pr\d = \c\rest_{p'} \op {}_{q'}\corest\d, \WHERE p'=q\th_p \ANd q'=p\th_q,
\]
and where $\op$ denotes concatenation, as above.
\item \label{CP3} The involution in $\PG(P)$ is given by $\ldb p_1,\ldots,p_k\rdb ^*=\ldb p_k,\ldots,p_1\rdb $.
\item \label{CP4} The projections of $\PG(P)$ have the form $\ldb p\rdb = p$, for $p\in P$, and consequently ${\bPP(\PG(P))= P}$.  Moreover, these projection algebras have the same operations, as
\[
p\pr q\pr p = q\th_p \qquad\text{for all $p,q\in P$.}
\]
\item \label{CP5} The idempotents of $\PG(P)$ have the form $p\pr q=\ldb p,q\rdb $, for $(p,q)\in{\F}$.
\end{enumerate}
\end{defn}

We have proved the following.

\begin{thm}\label{thm:C_P}
If $P$ is a projection algebra, then its chain semigroup $\PG(P)$ is a projection-generated (equivalently, idempotent-generated) regular $*$-semigroup whose projection algebra is precisely $P$.  \epfres
\end{thm}

We have now introduced the main concept of our study, the chain semigroup associated to a projection algebra, and the remaining sections of the paper investigate these semigroups from a number of different angles:
\bit
\item their incarnation as free objects in the category of regular $*$-semigroups (Section \ref{sect:F}),
\item their idempotent/biordered structure (Section \ref{sect:E}),
\item their presentations by generators and defining relations (Section \ref{sect:pres}),
\item their relationship to free (regular) idempotent-generated semigroups (Section \ref{sect:pres}), and
\item their topological structure (Section \ref{sect:top}).
\eit
In Sections \ref{sect:AGa} and \ref{sect:TLn} we present a number of natural examples.  
%

\section{Freeness of the chain semigroup}\label{sect:F}

In the previous section we showed how to construct the chain semigroup $\PG(P)=\bS(P,\ol\C,\nu)$ from a projection algebra $P$.  Now we will explain how $\PG(P)$ is rightfully thought of as `the free regular $*$-semigroup with projection algebra $P$'.  In categorical language, this is to say that the chain semigroups are the objects in the image of a left adjoint to the forgetful functor $\bPP:\RSS\to\PA$ from \eqref{eq:Pfunctor}.  (The precise meanings of these terms are given below.)  The forgetful functor in question maps a regular $*$-semigroup $S$ to its projection algebra $\bPP(S)$.  It follows from Proposition~\ref{prop:Cfunctor} (and the isomorphism $\RSS\cong\CPG$) that the assignment $P\mt\PG(P)$ is the object part of a functor $\PA\to\RSS$.  Our main goal here is to prove the following result (where again definitions are given below).

\begin{thm}\label{thm:CP}
The functor $\PA\to\RSS:P\mt\PG(P)$ is a left adjoint to the forgetful functor $\RSS\to\PA:S\mt \bPP(S)$, and $\PA$ is coreflective in $\RSS$.
\end{thm}

In fact, it will be more convenient to prove the following groupoid version of Theorem \ref{thm:CP}; the two theorems are equivalent via the isomorphism $\RSS\cong\CPG$.

\begin{thm}\label{thm:PCnu}
The functor $\PA\to\CPG:P\mt(P,\ol\C,\nu)$ is a left adjoint to the forgetful functor $\CPG\to\PA:(P,\G,\ve)\mt P$, and $\PA$ is coreflective in $\CPG$.
\end{thm}

We now give the (standard) definitions of the terms appearing in the above results; for more details see \cite{Awodey2010,MacLane1998}.

\begin{defn}\label{defn:nt}
Consider two categories $\bC$ and $\bD$, and a pair of functors ${\bF,\GG:\bC\to\bD}$.  A \emph{natural transformation} $\eta:\bF\to \GG$ is a family $\eta = (\eta_C)_{C\in v\bC}$, where each $\eta_C:\bF(C)\to \GG(C)$ is a morphism in $\bD$, and such that the following condition holds:
\bit
\item For every pair of objects $C,C'\in v\bC$, and for every morphism $\phi:C\to C'$ in $\bC$, the following diagram commutes:
\[
\begin{tikzcd}[sep=scriptsize]
\bF(C) \arrow{rr}{\bF(\phi)} \arrow[swap]{dd}{\eta_C} & ~ & \bF(C') \arrow{dd}{\eta_{C'}} \\%
~&~&~\\
\GG(C) \arrow{rr}{\GG(\phi)}& ~ & \GG(C').
\end{tikzcd}
\]
\eit
We call $\eta$ a \emph{natural isomorphism} if each $\eta_C$ is an isomorphism (in $\bD$), in which case we write~${\bF\cong\GG}$.
\end{defn}

\begin{defn}\label{defn:ad}
Consider two categories $\bC$ and $\bD$.  An \emph{adjunction} $\bC\to\bD$ is a triple $(\bF,\U,\eta)$, where $\bF:\bC\to\bD$ and $\U:\bD\to\bC$ are functors, and $\eta$ is a natural transformation $\id_\bC\to\U\bF$, such that the following condition holds:
\bit
\item For every pair of objects $C\in v\bC$ and $D\in v\bD$, and for every morphism $\phi:C\to\U(D)$ in~$\bC$, there exists a unique morphism $\ol \phi:\bF(C)\to D$ in $\bD$ such that the following diagram commutes:
\[
\begin{tikzcd}
C \arrow{rrdd}{\phi} \arrow[swap]{dd}{\eta_C}  \\%
~&~&~\\
\U(\bF(C)) \arrow{rr}{\U(\ol\phi)}& ~ & \U(D).
\end{tikzcd}
\]
\eit
In this set-up, $\bF$ and $\U$ are called the \emph{left} and \emph{right adjoints}, respectively, and $\eta$ is the \emph{unit} of the adjunction.  
The \emph{$\U$-free objects in $\bD$} are the objects in the image of $\bF$, i.e.~those of the form $\bF(C)$ for $C\in v\bC$.
\end{defn}

We are particularly interested in special adjunctions where we actually have the \emph{equality} ${\U\bF = \id_{\bC}}$.  Note that $\id=(\id_C)_{C\in v\bC}$ is clearly a natural transformation $\id_\bC\to\id_\bC$ for any category~$\bC$.  Lemma \ref{lem:coref} below concerns this situation, and speaks of so-called \emph{coreflective} (sub)categories.

\begin{defn}
A \emph{coreflective subcategory} of a category $\bD$ is a full subcategory $\bB$ whose inclusion functor $\bB\to\bD$ has a right adjoint.  We say a category is \emph{coreflective in $\bD$} if it is isomorphic to a coreflective subcategory of $\bD$.  That is, $\bC$ is coreflective in $\bD$ if there is a full embedding $\bC\to\bD$ that has a right adjoint.  
\end{defn}

In the above, a subcategory $\bB$ of $\bD$ is \emph{full} if it contains every morphism of $\bD$ between objects of $\bB$.  A \emph{full embedding} is a functor $\bC\to\bD$ that is injective on objects and morphisms, and whose image is full in $\bD$.


\newpage

\begin{lemma}\label{lem:coref}
Suppose $\bC$ and $\bD$ are categories, and $\bF:\bC\to\bD$ and $\U:\bD\to\bC$ are functors with $\U\bF=\id_{\bC}$, for which the following condition holds:
\bit
\item For every pair of objects $C\in v\bC$ and $D\in v\bD$, and for every morphism $\phi:C\to\U(D)$ in $\bC$, there exists a unique morphism $\ol \phi:\bF(C)\to D$ in $\bD$ such that $\phi=\U(\ol\phi)$.
\eit
Then
\ben
\item \label{coref1} $(\bF,\U,\id)$ is an adjunction,
\item \label{coref2} $\bC\cong\bF(\bC)$ is coreflective in $\bD$.
\een
\end{lemma}

\pf
\firstpfitem{\ref{coref1}}  This is a direct translation of Definition \ref{defn:ad} in the special case that $\U\bF=\id_\bC$ and $\eta=\id$.

\pfitem{\ref{coref2}}  This follows from the fact that all components of the unit $\id = (\id_C)$ are isomorphisms; see for example \cite[Theorem IV.3.1]{MacLane1998}.
%
\epf

\pf[\bf Proof of Theorem \ref{thm:PCnu}.]
We denote the functors in question by
\[
\bF:\PA\to\CPG:P\mt(P,\ol\C,\nu) \AND \U:\CPG\to\PA:(P,\G,\ve)\mt P.
\]
We prove the theorem by applying Lemma \ref{lem:coref}.  It is clear that $\U\bF=\id_\PA$, so we are left to verify the following condition:
\bit
\item For every projection algebra $P$, every chained projection groupoid $(P',\G,\ve)$, and every projection algebra morphism ${\phi:P\to\U(P',\G,\ve)=P'}$, there exists a unique chained projection functor $\ol \phi:\bF(P)=(P,\ol\C,\nu)\to (P',\G,\ve)$ such that $\phi=\U(\ol\phi)$.
\eit
To establish the existence of $\ol\phi$, we first define
\[
\varphi:\P\to\G \BY \p\varphi = \ve([\p]\C(\phi)), \qquad \text{i.e. } (p_1,\ldots,p_k)\varphi=\ve[p_1\phi,\ldots,p_k\phi] \text{ for  } \p\in\P.
\]
where here ${\P=\P(P)}$ is the path category of $P$.  
As with \eqref{eq:3comp}, $\varphi$ is the composite of three ordered $*$-morphisms, and is hence itself an ordered $*$-morphism.  We next show that $\Xi\sub\ker(\varphi)$, i.e.~that
\begin{equation}\label{eq:Xikervarphi}
\s\varphi=\t\varphi \qquad\text{for all $(\s,\t)\in\Xi$.}
\end{equation}
This is clear when $(\s,\t)$ has the form \ref{Om1} or \ref{Om2}.  Now suppose $(\s,\t)$ has the form~\ref{Om3}, so that
\[
\s = \lam(e,p,f) = (e,e\th_p,f) \ANd \t = \rho(e,p,f) = (e,f\th_p,f) \qquad\text{for some $p$-linked pair $(e,f)$.}
\]
We show that:
\ben
\item \label{ephifphi1} $(e\phi,f\phi)$ is $b$-linked in $\G$ for the morphism $b = p\phi\in\G$, and
\item \label{ephifphi2} $\lam(e\phi,b,f\phi) = \s\varphi $ and $\rho(e\phi,b,f\phi) = \t\varphi$.	
\een
Since $\lam(e\phi,b,f\phi) = \rho(e\phi,b,f\phi)$ in $\G$, this will complete the proof that $\s\varphi=\t\varphi$.  
For \ref{ephifphi1} we need to show that
\[
f\phi = (e\phi) \Th_b \th_{f\phi} \AND e\phi = (f\phi) \Th_{b^{-1}} \th_{e\phi}.
\]
Since $p\phi\in P'$ is a projection, we have $\Th_b=\Th_{p\phi}=\th_{p\phi}$ by \cite[equation (6.3)]{EPA2024}.  Combined with the fact that $\phi$ is a projection algebra morphism, it follows that
\[
(e\phi) \Th_b \th_{f\phi} = (e\phi) \th_{p\phi} \th_{f\phi} = (e\th_p\th_f)\phi = f\phi.
\]
An analogous calculation (keeping in mind $b^{-1}=(p\phi)^{-1}=p\phi$) gives $e\phi = (f\phi) \Th_{b^{-1}} \th_{e\phi}$, and completes the proof of \ref{ephifphi1}.  
%
For \ref{ephifphi2} we first note that
\begin{align*}
\lam(e\phi,p\phi,f\phi) &= \ve[e\phi,e_1] \circ {}_{e_1}\corest(p\phi) \circ \ve[f_1,f\phi] \\
&= \ve[e\phi,e_1] \circ e_1 \circ \ve[f_1,f\phi] \\
&= \ve[e\phi,e_1] \circ \ve[f_1,f\phi] ,
\end{align*}
where $e_1 = (e\phi)\th_{p\phi} = (e\th_p)\phi$ and $f_1 = (e\phi)\Th_{p\phi} = (e\phi)\th_{p\phi} = (e\th_p)\phi \ ( = e_1)$.  Thus, continuing from above we have
\[
\lam(e\phi,p\phi,f\phi)  = \ve[e\phi,e_1] \circ \ve[f_1,f\phi] = \ve[e\phi,(e\th_p)\phi] \circ \ve[(e\th_p)\phi,f\phi] = \ve[e\phi,(e\th_p)\phi,f\phi] = \s\varphi.
\]
Analogously, $\rho(e\phi,p\phi,f\phi) = \t\varphi$, completing the proof of \ref{ephifphi2}.

Now that we have proved \eqref{eq:Xikervarphi}, it follows that there is a well-defined ordered groupoid functor
\[
\ol\phi:\ol\C(=\P/{\ssim})\to\G \GIVENBY \ldb\p\rdb\ol\phi = \p\varphi {}= \ve([\p]\C(\phi))  \qquad\text{for $\p\in\P$.}
\]
We next check that $\ol\phi$ is in fact a chained projection functor $(P,\ol\C,\nu)\to(P',\G,\ve)$, for which we need to show that:
\ben\addtocounter{enumi}{2}
\item \label{phibar1} the object map $v\ol\phi=\ol\phi|_P$ is a projection algebra morphism $P\to P'$, and
\item \label{phibar2} $\ol\phi$ respects the evaluation maps, in the sense that the following diagram commutes:
\[
\begin{tikzcd}[sep=scriptsize]
\C(P) \arrow{rr}{\C(v\ol\phi)} \arrow[swap]{dd}{\nu} & ~ & \C(P')\arrow{dd}{\ve} \\%
~&~&~\\
\ol\C \arrow{rr}{\ol\phi}& ~ & \G.
\end{tikzcd}
\]
\een
For \ref{phibar1}, we note that 
$p\ol\phi = p\varphi = \ve[p\phi] = p\phi$ for all $p\in P$, so that $v\ol\phi=\phi$ is a projection algebra morphism by assumption.  For \ref{phibar2}, if $\p\in\P$ then
\[
\ve([\p] \C(\phi)) 
= \p\varphi = \ldb \p\rdb\ol\phi = (\nu[\p])\ol\phi.
\]
So $\ol\phi$ is indeed a chained projection functor, and it follows from the proof of \ref{phibar1} above that $\U(\ol\phi) = \ol\phi|_P = \phi$.  

We have now established the existence of $\ol\phi$.  For uniqueness, suppose the chained projection functor ${\psi:(P,\ol\C,\nu)\to(P',\G,\ve)}$ also satisfies $\U(\psi)=\phi$.  
Then for any ${\p\in\P}$ we have
\begin{align*}
\ldb\p\rdb\psi = (\nu[\p])\psi &= \ve([\p]\C(v\psi)) &&\text{since $\psi$ respects evaluation maps}\\
&= \ve([\p]\C(\phi)) &&\text{since $v\psi=\U(\psi)=\phi$}\\
&= \ldb\p\rdb\ol\phi,
\end{align*}
and so $\psi=\ol\phi$.
\epf

Now that we have proved Theorem \ref{thm:CP} (via Theorem \ref{thm:PCnu} and the isomorphism $\RSS\cong\CPG$), it follows that the chain semigroups are the $\bPP$-free objects in $\RSS$.  
Henceforth, we call~$\PG(P)$ the \emph{free (projection-generated) regular $*$-semigroup over the projection algebra $P$}.  
As noted earlier, the notation is inspired by an analogy with the free (idempotent-generated) semigroup over a biordered set~$E$, which is denoted $\IG(E)$.  The relationship between $\PG(P)$ and $\IG(E)$ will be a key topic for the remainder of the paper, starting in Section~\ref{sect:E}.

We conclude the current section by drawing out two purely semigroup-theoretical results.  

\begin{thm}\label{thm:free}
If $P$ is a projection algebra, then
\ben
\item \label{free1} $\PG(P)$ is a regular $*$-semigroup with projection algebra $P$,
\item \label{free2} for any regular $*$-semigroup $S$, and any projection algebra morphism $\phi:P\to \bPP(S)$, there is a unique $*$-semigroup homomorphism~${\ol\phi:\PG(P)\to S}$ such that the following diagram commutes (where both vertical maps are inclusions):
\een
\[
\begin{tikzcd}[sep=small]
P \arrow{rr}{\ \ \phi} \arrow[swap,hookrightarrow]{dd} & ~ & \bPP(S) \arrow[hookrightarrow]{dd} \\%
~&~&~\\
\PG(P) \arrow{rr}{\ol\phi\ \ }& ~ & S,
\end{tikzcd}
\]
\ben
\item[] and moreover any $*$-semigroup homomorphism $\PG(P)\to S$ has the form $\ol\phi$ for some projection algebra morphism $\phi:P\to \bPP(S)$.  
\een
\end{thm}

\pf
\firstpfitem{\ref{free1}}  This is contained in Theorem \ref{thm:C_P}.

\pfitem{\ref{free2}}  The existence and uniqueness of $\ol\phi$ follows from the proof of Theorem \ref{thm:PCnu} and the isomorphism $\RSS\cong\CPG$.  Given any $*$-semigroup morphism $\psi:\PG(P)\to S$, we have $\psi=\ol{v\psi}$.
\epf

The previous theorem has the following consequence:

\begin{thm}\label{thm:free2}
Let $P$ and $P'$ be projection algebras.  Then for any projection algebra morphism $\phi:P\to P'$, there is a unique $*$-semigroup homomorphism~${\ol\phi:\PG(P)\to \PG(P')}$ such that the following diagram commutes (where both vertical maps are inclusions):
\[
\begin{tikzcd}[sep=small]
P \arrow{rr}{\phi} \arrow[swap,hookrightarrow]{dd} & ~ & P' \arrow[hookrightarrow]{dd} \\%
~&~&~\\
\PG(P) \arrow{rr}{\ol\phi}& ~ & \PG(P'),
\end{tikzcd}
\]
and moreover any $*$-semigroup homomorphism $\PG(P)\to \PG(P')$ has the form $\ol\phi$ for some projection algebra morphism $\phi:P\to P'$.  \epfres
\end{thm}

\begin{rem}\label{rem:coref}
Let us pause briefly to illustrate the significance of coreflectivity in Theorems~\ref{thm:CP} and \ref{thm:PCnu},
by a comparison with the category $\Sgp$ of semigroups and the set-based free objects within it.
The latter are the semigroups $X^+$, consisting of all words over a set $X$ under the operation of concatenation.
The assignment $X\mapsto X^+$ can be viewed as a functor ${\bF':\Set\rightarrow \Sgp}$.
It is a left adjoint to the forgetful functor $\U':\Sgp\rightarrow \Set$, which maps any semigroup to its underlying set.
However, $\U'\bF'$ is not naturally isomorphic, let alone equal, to the identity on~$\Sgp$: indeed,
$\U'\bF'(X)$ is the underlying set of the free semigroup $X^+$.
On the level of morphisms this is manifested by the fact that not every morphism $X^+\rightarrow Y^+$ arises from a mapping $X\rightarrow Y$.
\end{rem}

\begin{rem}\label{rem:GE}
Lawson in \cite[Theorem 2.2.4]{Lawson1998} states that groups form a \emph{reflective} subcategory of inverse semigroups.  The \emph{left} adjoint to the inclusion maps an inverse semigroup to its maximum group image (the quotient by the least congruence that identifies all idempotents).
On the other hand one can check from the definitions that semilattices (commutative idempotent semigroups) form a coreflective subcategory of inverse semigroups.  
Similarly, the category of (regular) biordered sets is coreflective in the category of (regular) semigroups; see \cite[Theorem~3.40]{Nambooripad_book} and \cite[Theorem 6.10]{Nambooripad1979}.  Coreflectivity of $\PA$ in $\RSS$ can be viewed as a `regular $*$-analogue' of these last two facts.
\end{rem}

\begin{rem}\label{rem:initial}
There is another way to view the semigroups $\PG(P)$ as free objects.  Namely, for any (fixed) projection algebra $P$, there is a category $\RSS(P)$ with:
\bit
\item objects all regular $*$-semigroups with projection algebra $P$, and
\item morphisms all $*$-morphisms that map projections identically.
\eit
We see then that $\PG(P)$ is an initial object in this category, meaning that for every object $S$ of $\RSS(P)$ there is exactly one morphism $\PG(P)\to S$ (in this category).  This follows by applying Theorem \ref{thm:free}\ref{free2} to the identity morphism $\phi=\id_P:P\to P=\bPP(S)$.  Note that the image of the morphism $\PG(P)\to S$ is the projection-generated subsemigroup of $S$, which of course is equal to the idempotent-generated subsemigroup of $S$.
\end{rem}

\section{Projection algebras and biordered sets}\label{sect:E}

We have just seen that the chain semigroups are the $\bPP$-free objects in $\RSS$, where ${\bPP:\RSS\to\PA}$ is the forgetful functor that maps a regular $*$-semigroup $S$ to its underlying projection algebra~$\bPP(S)$.  Another forgetful functor $\bEE:\RSS\to\RSBS$ has been considered (implicitly) in the literature~\cite{NP1985}, where $\RSBS$ denotes the category of \emph{regular $*$-biordered sets}; this will be defined formally below.
It is then natural to ask whether this forgetful functor has an adjoint, and if so what the $\bEE$-free objects are.  Perhaps more fundamentally, we would like to understand the relationship between the categories $\PA$ and $\RSBS$.  It turns out that these categories are in fact \emph{equivalent}, as we show in Theorem \ref{thm:equiv} below.  This then has the consequence that $\bEE$-free objects do indeed exist, but that they are the same as the $\bPP$-free objects.

To establish this equivalence, we need functors
\[
\bE:\PA\to\RSBS \AND \bP:\RSBS\to\PA
\]
with natural isomorphisms $\bP\circ\bE \cong \id_\PA$ and $\bE\circ\bP \cong \id_\RSBS$.  These functors are constructed in Subsections~\ref{subsect:EP} and \ref{subsect:PE}, and the category equivalence is established in Subsection \ref{subsect:equiv}.  A salient part of the argument is Proposition \ref{prop:E}, which shows that regular $*$-semigroups with the same projection algebra have isomorphic $*$-biordered sets.

\subsection{Preliminaries on biordered sets}\label{subsect:Eprelim}

We begin with the necessary definitions; for more background, and proofs, see \cite[Chapter 1]{Nambooripad1979} and \cite[Section 2]{NP1985}.  

The set $E=\bEE(S)=\set{e\in S}{e^2=e}$ of all idempotents of a semigroup $S$ can be given the structure of a partial algebra called a \emph{biordered set}.  This structure can be conveniently described using Easdown's arrow notation \cite{Easdown1985}.  For $e,f\in E$ we write 
\[
e \larr f \iff e=ef \text{ ($e$ is a left zero for $f$)}  \AND e \rarr f \iff e=fe   \text{ ($e$ is a right zero for $f$)}.
\]
These relations were originally denoted $\omega^l={\larr}$ and $\omega^r={\rarr}$ in \cite{Nambooripad1979}.
Both are pre-orders (reflexive and transitive relations), and if $e\larr f$, $e\rarr f$, $e\arrl f$ or $e\arrr f$, then $ef$ and $fe$ are both idempotents (at least one of which is equal to $e$ or $f$).  It follows that we can define a partial operation $\cdot$ on $E$ with domain
\[
\BP(E) = {\larr}\cup{\rarr}\cup{\arrl}\cup{\arrr} = \set{(e,f)\in E}{\{e,f\}\cap\{ef,fe\}\not=\es},
\]
and with $e\cdot f = ef$ for $(e,f)\in\BP(E)$.  The pairs in $\BP(E)$ are called \emph{basic pairs}.  The partial algebra $(E,{\cdot})$ is the \emph{biordered set} of $S$, or \emph{boset} for short.  Since the (partial) product in $E$ is just a restriction of the (total) product in $S$, we denote it simply by juxtaposition.  In what follows, we will use various combinations of arrows, specifically
\[
{\llarr} = {\larr}\cap{\arrl} \COMMA {\rrarr} = {\rarr}\cap{\arrr} \AND {\lrarr} = {\larr}\cap{\rarr} ,
\]
so for example $e \llarr f \iff [e\larr f$ and $e\arrl f] \iff [e=ef$ and $f=fe]$.  Note that $\lrarr$ is a partial order.

There is an axiomatic definition of abstract bosets (with no reference to any over-semigroup), but we will not need to give that here, as it is known that any abstract boset is the boset of idempotents of some semigroup \cite{Easdown1985}.

We denote by $\BS$ the category of bosets.  Morphisms in $\BS$ are called \emph{bimorphisms},\footnote{The term `bimorphism' comes from Nambooripad's original work \cite{Nambooripad1979}, as a contraction of `biordered set morphism', and should not be confused with other uses of the term, for example to mean a bijective morphism of graphs.} and are simply morphisms of partial algebras, i.e.~functions $\phi:E\to E'$ such that for all $e,f\in E$,
\begin{equation}\label{eq:bimorphism}
(e,f) \in \BP(E) \Implies (e\phi,f\phi)\in\BP(E') \ANd (e f)\phi = (e\phi)  (f\phi).
\end{equation}
An isomorphism in $\BS$ will be called a \emph{bisomorphism}; these are the bijections $\phi:E\to E'$ for which $\phi$ and $\phi^{-1}$ are both bimorphisms.  It is easy to see that a bijection $\phi:E\to E'$ is a bisomorphism if and only if
\[
\text{$\phi$ is a bimorphism} \AND 
\BP(E') = \BP(E)\phi = \set{(e\phi,f\phi)}{(e,f)\in\BP(E)}.
\]

A boset is called \emph{regular} if it is the boset of a regular semigroup.  Such regular bosets can be abstractly characterised as the bosets $E$ 
whose \emph{sandwich sets} $\SS(e,f)$ are non-empty for all $e,f\in E$.  These sandwich sets are defined as follows.  We begin by defining the set (of `mixed' common lower bounds):
\[
\MM(e,f) = \set{g\in E}{e \arrl g \rarr f} = \set{g\in E}{ge=g=fg} = \set{g\in E}{fge=g}.
\]
This set has a pre-order $\pre$ defined by $h\pre g \iff [eh \rarr eg$ and $hf \larr gf]$.  The sandwich set is then
\[
\SS(e,f) =  \set{g\in \MM(e,f)}{h\pre g\text{ for all }h\in \MM(e,f)},
\]
the set of all $\pre$-maximum elements in $\MM(e,f)$.  It turns out (see \cite[Theorem 1.1]{Nambooripad1979}) that if $E$ is a regular boset, and if~$S$ is \emph{any} regular semigroup with $E=\bEE(S)$, then
\begin{equation}\label{eq:SSef}
\SS(e,f) = \set{g\in E}{egf=ef \text{ and } fge=g \text{ in } S} \qquad\text{for any $e,f\in E$}.
\end{equation}
Note here that the products $egf$ and $ef$ are taken in $S$, but need not be defined in $E$ itself; these products need not even be idempotents.  

We write $\RBS$ for the category of regular bosets.  Morphisms in $\RBS$ are called \emph{regular bimorphisms}, and are the bimorphisms $\phi:E\to E'$ (as above) that map sandwich sets into sandwich sets, meaning that
\[
\SS(e,f)\phi \sub \SS(e\phi,f\phi) \qquad\text{for all $e,f\in E$.}
\]
Regular bosets are isomorphic in $\RBS$ if and only if they are isomorphic in $\BS$.  This is because sandwich sets are defined directly from the (partial) `multiplication tables' of bosets, and are hence preserved by bisomorphisms.

\newpage

Following \cite{NP1985}, a \emph{$*$-boset} is a partial algebra $E=(E,{\cdot},{}^*)$, where $(E,{\cdot})$ is a boset, and ${}^*$ is a unary operation $E\to E:e\mt e^*$ satisfying the following, for all $e,f\in E$:
\begin{enumerate}[label=\textup{\textsf{(SB\arabic*)}},leftmargin=13mm]
\item \label{SB1} $(e^*)^*=e$,
\item \label{SB2} $(e,f)\in\BP(E) \implies (e^*,f^*)\in\BP(E)$ and $(e f)^* = f^* e^*$.
\end{enumerate}
We say $E=(E,{\cdot},{}^*)$ is a \emph{regular $*$-boset} if $(E,\cdot)$ is regular, and we additionally have: 
\begin{enumerate}[label=\textup{\textsf{(SB\arabic*)}},leftmargin=13mm]\addtocounter{enumi}{2}
\item \label{SB3} For all $e\in E$, there exist elements $s=s(e)$ and $t=t(e)$ of $E$ such that 
\[
\begin{tikzpicture}[scale=1.5]
\node (e) at (0,1) {$e$};
\node (f) at (1,1) {$s$};
\node (g) at (0,0) {$t$};
\node (h) at (1,0) {$e^*$,};
\draw[<->] (e)--(f);
\draw[<->] (g)--(h);
\draw[>-<] (e)--(g);
\draw[>-<] (f)--(h);
\end{tikzpicture}
\]
and such that $e^*(gs)=(tg)e^*$ for all $g\lrarr e$.
\end{enumerate}
(We note that these were called `special $*$-biordered sets' in \cite{NP1985}, where regular $*$-semigroups were also called `special $*$-semigroups'.  Condition \ref{SB3} is known as \emph{$\tau$-commutativity}, as it can be stated in terms of commutative diagrams involving natural maps between down-sets of idempotents in the poset $(E,{\lrarr})$.)

If $S$ is a regular $*$-semigroup, then the boset $E=\bEE(S)$ becomes a regular $*$-boset whose unary operation is the restriction of the involution of $S$.  The elements $s,t\in E$ in \ref{SB3} are $s=ee^*$ and $t=e^*e$; for any $g\lrarr e$, the products $e^*(gs)$ and $(tg)e^*$ both evaluate to $e^*ge^*$.  
Conversely, we have the following:

\begin{thm}[{see \cite[Corollary 2.7]{NP1985}}]\label{thm:NS}
Any regular $*$-boset is the $*$-boset of a regular $*$-semigroup.  \epfres
\end{thm}

We denote by $\RSBS$ the category of regular $*$-bosets.  Morphisms in $\RSBS$ are called \emph{regular $*$-bimorphisms}, and are the regular bimorphisms $\phi:E\to E'$ (as above) that respect the involutions, meaning that
\[
(e^*)\phi=(e\phi)^* \qquad\text{for all $e\in E$.}
\]
As above, regular $*$-bosets $E$ and $E'$ are isomorphic in $\RSBS$ if and only if there is a bisomorphism $E\to E'$ that also respects the involutions.  

The assignment $S\mt \bEE(S)$ is the object part of a (forgetful) functor
\begin{equation}\label{eq:Efunctor}
\bEE:\RSS\to\RSBS.
\end{equation}
A $*$-morphism $\phi:S\to S'$ in $\RSS$ is sent to its restriction $\bEE(\phi) = \phi|_{\bEE(S)}:\bEE(S)\to \bEE(S')$, which is a regular $*$-bimorphism in $\RSBS$.

\subsection[From projection algebras to $*$-biordered sets]{\boldmath From projection algebras to $*$-biordered sets}\label{subsect:EP}

Theorem \ref{thm:equiv} below establishes a category equivalence $\PA\leftrightarrow\RSBS$.  For this we need functors in both directions between the two categories.  We can immediately obtain a functor 
\begin{equation}\label{eq:bEfunctor}
\bE:\PA\to\RSBS
\end{equation}
by composing the functors
\[
\PA\to\RSS:P\mt\PG(P) \AND \RSS\to\RSBS:S\mt \bEE(S) 
\]
from Theorem \ref{thm:CP} and \eqref{eq:Efunctor}.  
Note that the functors $\RSS\to\RSBS$ and $\PA\to\RSBS$ in~\eqref{eq:Efunctor} and~\eqref{eq:bEfunctor} are both denoted by $\bE$.
As these take different kinds of arguments (regular $*$-semigroups or projection algebras, and their morphisms), it will always be clear from context which one is meant.

The functor $\bE:\PA\to\RSBS$ in \eqref{eq:bEfunctor} maps a projection algebra $P$ to the regular $*$-boset of~$\PG(P)$.  Using~\ref{CP5}, we have
\begin{equation}\label{eq:EP}
\EP = \bEE(\PG(P)) = \set{\ldb p,q\rdb}{(p,q)\in{\F}} 
\end{equation}
as a set.  We will describe the biordered structure of $\EP$ in Proposition \ref{prop:EPGP} below; the involution is of course given by $\ldb p,q\rdb^*=\ldb q,p\rdb$.

The functor $\bEE:\PA\to\RSBS$ maps a projection algebra morphism $\phi:P\to P'$ to the regular $*$-bimorphism
\[
\bE(\phi) = \bEE(\ol\phi): \EP\to \EPd,
\]
where $\ol\phi:\PG(P)\to\PG(P')$ is the $*$-morphism from Theorem \ref{thm:free2}.  So $\bEE(\phi)$ is the restriction of~$\ol\phi$ to $\bEE(P)=\bE(\PG(P))$; explicitly, we have $\ldb p,q\rdb\bEE(\phi) = \ldb p\phi,q\phi\rdb$ for $(p,q)\in{\F}$.


For the next two proofs, we note that the product of idempotents $e=\ldb p,q\rdb$ and $f=\ldb r,s\rdb$ in $\PG(P)$ is given by
\begin{equation}\label{eq:p'q'r's'}
e \pr f = \ldb p',q',r',s'\rdb, \WHERE 
p' = r\th_q\th_p \COMMa
q' = r\th_q \COMMa
r' = q\th_r \ANd
s' = q\th_r\th_s.
\end{equation}
In particular, if $e\pr f$ happens to be an idempotent, then
\begin{equation}\label{eq:p's'}
e\pr f = \ldb p',s'\rdb = \ldb r\th_q\th_p, q\th_r\th_s \rdb.
\end{equation}

\begin{prop}\label{prop:EPGP}
If $P$ is a projection algebra, and if $e=\ldb p,q\rdb$ and $f=\ldb r,s\rdb$, then 
\ben
\item $e \larr f  \textup{ (i.e.~$e=e\pr f$)} \iff \br(e)\leq \br(f) \iff q\leq s$, in which case $ f \pr e= \ldb p\th_s\th_r,s\th_p\th_q \rdb$,
\item $e\rarr f \textup{ (i.e.~$e=f\pr e$)} \iff \bd(e)\leq\bd(f)  \iff p\leq r$, in which case $ e \pr f = \ldb r\th_q\th_p, q\th_r\th_s \rdb$.
\een
\end{prop}

\pf
For the first part (the second is dual), it is enough to show that $e=e\pr f \iff q\leq s$, as the expression for $f\pr e$ will then follow from \eqref{eq:p's'}.  Write $e\pr f = \ldb p',q',r',s'\rdb$, as in \eqref{eq:p'q'r's'}.
If $e=e\pr f$, then $q=s'=q\th_r\th_s\leq s$.  
Conversely, suppose $q\leq s$, so that $q=q\th_s$.  Combining this with $\th_s=\th_s\th_r\th_s$ (which holds by \ref{PA5}), we calculate
\[
q = q\th_s = q\th_s\th_r\th_s = q\th_r\th_s = s'.
\]
From $q\leq s \F r$ it follows from \ref{PA2} that $q\leqF r$, and so $q=r\th_q = q'$.  But also
\[
p' = r\th_q\th_p = q\th_p = p,
\]
and so $e\pr f = \ldb p',q',r',s'\rdb = \ldb p,q,r',q\rdb = \ldb p,q\rdb = e$, as required.
\epf

The functor $\bE:\PA\to\RSBS$ constructs $\EP$ from a projection algebra $P$ by first passing through $\PG(P)$.  Perhaps surprisingly, it turns out that we could have passed through \emph{any} regular $*$-semigroup $S$ with projection algebra $P=\bPP(S)$, and we would obtain the same $*$-boset up to isomorphism, $\bEE(S)\cong\EP$:

\begin{prop}\label{prop:E}
If $S$ is a regular $*$-semigroup with projection algebra $P=\bPP(S)$, then the map
\[
\EP \to \bEE(S) : \ldb p,q\rdb \mt pq
\]
is an isomorphism of $*$-bosets.
\end{prop}

\pf
As in Remark \ref{rem:initial}, there exists a $*$-morphism $\Phi(=\ol\id_P):\PG(P)\to S$ with $p\Phi=p$ for all $p\in P$.  Applying the functor $\bEE$ from \eqref{eq:Efunctor}, we obtain the regular $*$-bimorphism 
\[
\phi = \bEE(\Phi) = \Phi|_{\EP}:\EP\to \bEE(S),  
\]
which is the map from the statement.
As explained in Subsection \ref{subsect:Eprelim}, and since we already know that~$\phi$ is a (regular) $*$-bimorphism, we can complete the proof that this is an isomorphism by checking that:
\ben
\item \label{propE1} $\phi$ is a bijection, and
\item \label{propE2} $\BP(\EP)\phi = \BP(\bEE(S))$.
\een

\pfitem{\ref{propE1}}  This follows by applying \ref{RS7} in both $S$ and $\PG(P)$, and keeping in mind $p\pr q=\ldb p,q\rdb$.


\pfitem{\ref{propE2}}  Since $\phi$ is a bijection, this amounts to showing that
\[
(e,f)\in\BP(\EP) \iff (e\phi,f\phi)\in\BP(\bEE(S)) \qquad\text{for all $e,f\in\EP$.}
\]
The forward implication is clear, as $\phi$ is a bimorphism.  Conversely, suppose $(e\phi,f\phi)\in\BP(\bEE(S))$, where $e,f\in\EP$, and write $e=\ldb p,q\rdb$ and $f=\ldb r,s\rdb$.  By symmetry we may assume that ${e\phi\larr f\phi}$, i.e.~$e\phi=(e\phi)(f\phi)$; we complete the proof by showing that $e\larr f$, i.e.~${e=e\pr f}$.  To do so, first write $e\pr f = \ldb p',q',r',s'\rdb$, as in~\eqref{eq:p'q'r's'}.
By Theorem \ref{thm:iso}, the $*$-morphism $\Phi:\PG(P)\to S$ is also a chained projection functor ${\bG(\PG(P)) = (P,\ol\C,\nu) \to \bG(S)=(P,\G,\ve)}$.  Since $\Phi$ is the identity on $P=v\ol\C=v\G$, it follows that $\Phi$ is a $v$-functor $\ol\C\to\G$, and so
\[
s' = \br(e\pr f) = \br((e\pr f)\Phi) = \br((e\Phi)(f\Phi)) = \br((e\phi)(f\phi)) = \br(e\phi) = \br(e) = q.
\]
Consequently, $q=s'=q\th_r\th_s\leq s$, and we then obtain $e\larr f$ from Proposition \ref{prop:EPGP}.
\epf

\begin{rem}
Combining Theorem \ref{thm:NS} and Proposition \ref{prop:E}, it follows that every regular $*$-boset is isomorphic to $\EP$ for some projection algebra $P$.  
\end{rem}

\begin{rem}
Another approach to constructing a boset $\bE(P)$ from a projection algebra $P$ would be to take the underlying set
\[
\bE(P) = {\F} = \set{(p,q)\in P\times P}{p=q\th_p\text{ and }q=p\th_q},
\]
and define the $\larr$ and $\rarr$ pre-orders, and basic products, as in Proposition \ref{prop:EPGP}.  One would then need to check that the boset axioms are satisfied.  Taking this approach, Proposition \ref{prop:E} would then state that $(p,q)\mt pq$ is an isomorphism $\bE(P)\to\bE(S)$ for any regular $*$-semigroup $S$ with projection algebra $P=\bP(S)$.
\end{rem}

\subsection[From $*$-biordered sets to projection algebras]{\boldmath From $*$-biordered sets to projection algebras}\label{subsect:PE}

Now that we have constructed a functor $\bE:\PA\to\RSBS$, we wish to construct a functor $\bP:\RSBS\to\PA$ in the reverse direction.  (Again, this functor $\RSBS\to\PA$ has the same name as the forgetful functor $\bPP:\RSS\to\PA$ considered earlier.) This is somewhat more involved than the functor $\bE$, which was simply the composition of two previously existing functors.

Consider a regular $*$-boset $E=(E,{\cdot},{}^*)$.  The underlying set of the projection algebra $\PE$ is simply the set of fixed points under the involution,
\[
\PE = \set{p\in E}{p=p^*},
\]
the elements of which are called the \emph{projections} of $E$.  To define the operations in $\PE$ we need the following special case of \cite[Proposition 2.5]{NP1985}.  We give a simple adaptation of the proof for completeness.  

\begin{lemma}\label{lem:epq}
If $E$ is a regular $*$-boset, and if $p,q\in \PE$, then there exists a unique element $e=e(p,q)$ of the sandwich set $\SS(p,q)$ for which $pe,eq\in\PE$.  Moreover, we have 
\[
e = qp \AND pe = pqp
\]
in any regular $*$-semigroup $S$ with $*$-boset $E=\bEE(S)$.
\end{lemma}

\pf
Fix an arbitrary regular $*$-semigroup $S$ with $E=\bEE(S)$.  
We obtain $qp\in \SS(p,q)$ from~\eqref{eq:SSef} and \ref{RS2}, and $p(qp),(qp)q\in\PE$ from \ref{RS5}.

For uniqueness, suppose $e\in \SS(p,q)$ is such that $pe,eq\in \PE$.  We first claim that
\begin{equation}\label{eq:pqppe}
pqp=pe \AND qpq=eq.
\end{equation}
We prove the first, and the second is analogous. Define the projections $s=pqp$ and $t=pe$.  Using \ref{RS2} and \eqref{eq:SSef} in the indicated places, we calculate
\[
s = \ul{pq}p = \ul{pe}qp = tqp \implies s=ts \AND t = p\ul{e} = \ul{pq}ep = \ul{pqp}qep = sqep \implies t = st.
\]
Since $s$ and $t$ are projections it then follows that $s=s^*=(ts)^* = s^*t^* = st = t$, as claimed.  Combining \ref{RS2}, \eqref{eq:SSef} and \eqref{eq:pqppe}, we obtain
\[
e = \ul{e}e = e\ul{pe} = \ul{ep}qp = \ul{eq}p = \ul{qpqp} = qp.  \qedhere
\]
\epf

From now on we fix the notation $e(p,q)$ from Lemma \ref{lem:epq}.  It is important to note that while the products $qp$ and $pqp$ in this lemma are taken in the semigroup $S$, the products $pe$ and~$eq$ exist in the boset $E$ itself.

\begin{lemma}\label{lem:ephi}
If $\phi:E\to E'$ is a regular $*$-bimorphism, then $e(p,q)\phi = e(p\phi,q\phi)$ for all $p,q\in\PE$.
\end{lemma}

\pf
With $e=e(p,q)$, Lemma \ref{lem:epq} gives $e\in\SS(p,q)$ and $pe,eq\in \PE$.  We deduce ${e\phi\in\SS(p\phi,q\phi)}$ from regularity of $\phi$, and $(p\phi)(e\phi),(e\phi)(q\phi)\in \PEd$ because $\phi$ is a $*$-bimorphism.  It then follows from uniqueness in Lemma \ref{lem:epq} that $e\phi = e(p\phi,q\phi)$.
\epf


\begin{defn}\label{defn:PE}
For a regular $*$-boset $E$, we define $\PE$ to be the projection algebra with underlying set $\PE=\set{p\in E}{p=p^*}$, and operations given by
\[
q\th_p = p e(p,q) \qquad\text{for $p,q\in\PE$.}
\]
This is well defined by Lemma \ref{lem:epq}, which also tells us that $\PE=\bPP(S)$ is the projection algebra of any regular $*$-semigroup $S$ with $*$-boset $E=\bEE(S)$.

For a regular $*$-bimorphism ${\phi:E\to E'}$, we define $\bP(\phi)$ to be the restriction
\begin{equation}\label{eq:Pphi}
\bP(\phi) = \phi|_{\PE}:\PE\to\PEd.
\end{equation}
(Note that $\phi$ maps projections to projections because it is a $*$-bimorphism.)
\end{defn}

\begin{prop}\label{prop:bPfunctor}
$\bP$ is a functor $\RSBS\to\PA$.
\end{prop}

\pf
To show that $\bP(\phi):\PE\to \PEd$ as in \eqref{eq:Pphi} is a projection algebra morphism, let ${p,q\in\PE}$.  We then use Lemma \ref{lem:ephi} to calculate
\[
(q\th_p)\phi = (pe(p,q))\phi = (p\phi)(e(p,q)\phi) = (p\phi)e(p\phi,q\phi) = (q\phi)\th_{p\phi}.
\]
It is again clear that the laws $\bP(\phi\circ\phi')=\bP(\phi)\circ\bP(\phi')$ and $\bP(\id_E)=\id_{\bP(E)}$ hold.
\epf

\subsection{A category equivalence, and more on freeness}\label{subsect:equiv}

We now establish the promised category equivalence.

\begin{thm}\label{thm:equiv}
We have
\[
\bP\circ\bE = \id_\PA \AND \bE\circ\bP \cong \id_\RSBS.
\]
Consequently, the functors $\bP$ and $\bE$ furnish an equivalence of the categories $\PA$ and $\RSBS$.
\end{thm}

\pf
To show that $\bP\circ\bE=\id_\PA$ we need to show that
\bit
\item $\bP(\bE(P)) = P$ for any projection algebra $P$, and
\item $\bP(\bE(\phi)) = \phi$ for any projection algebra morphism $\phi:P\to P'$.
\eit
Since $\bEE(\PG(P)) = \EP$, it follows from Definition \ref{defn:PE} that $\bP(\EP) = \bPP(\PG(P)) = P$.  For the statement concerning $\phi$, we follow the definitions to compute
\[
\bP(\bE(\phi)) = \bE(\phi)|_{\bP(\EP)} = \bE(\phi)|_P = \phi.
\]

To show that $\bE\circ\bP\cong\id_\RSBS$, we will construct a natural isomorphism $\eta:\bE\circ\bP\to\id_\RSBS$.  Towards this, we claim that for any regular $*$-boset $E$, the map
\[
\eta_E: \bE(\bP(E)) \to E: \ldb p,q\rdb \mt e(q,p)
\]
is a $*$-bisomorphism.  To see this, let $S$ be a regular $*$-semigroup with $*$-boset $E=\bEE(S)$, and write $P=\PE=\bPP(S)$.  Since $pq=e(q,p)$ in $S$ by Lemma \ref{lem:epq}, it follows that $\eta_E$ is precisely the isomorphism from Proposition \ref{prop:E}.

We now show that $\eta=(\eta_E)$ is a natural isomorphism ${\bE\circ\bP\to\id_\RSBS}$.  Since each $\eta_E$ is an isomorphism, it remains to check that for any regular $*$-bimorphism $\phi:E\to E'$ the following diagram commutes:
\[
\begin{tikzcd}
\bE\circ\bP(E) \arrow{rr}{\bE\circ\bP(\phi)} \arrow[swap]{dd}{\eta_E} & ~ & \bE\circ\bP(E') \arrow{dd}{\eta_{E'}} \\%
~&~&~\\
E \arrow{rr}{\phi}& ~ & E'.
\end{tikzcd}
\]
But the elements of $\bE\circ\bP(E)$ have the form $\ldb p,q\rdb$ for $\F$-related projections $p,q\in\PE$,
and we use Lemma \ref{lem:ephi} to calculate
\[
\ldb p,q\rdb \mtlab{\eta_E} e(q,p) \mtlab{\phi} e(q\phi,p\phi)
\AND
\ldb p,q\rdb \mtlab{\bE\circ\bP(\phi)} \ldb p\phi,q\phi\rdb \mtlab{\eta_{E'}} e(q\phi,p\phi).  \qedhere
\]
\epf

Combining Theorems \ref{thm:CP} and \ref{thm:equiv}, we obtain the following:

\begin{thm}\label{thm:CPE}
The functor $\RSBS\to\RSS:E\mt\PG(\PE)$ is a left adjoint to the forgetful functor $\RSS\to\RSBS:S\mt \bEE(S)$, and $\RSBS$ is coreflective in~$\RSS$.  \epfres
\end{thm}

It follows that the chain semigroups are the free objects in $\RSS$ with respect to both forgetful functors
\[
\bPP:\RSS\to\PA:S\mt \bPP(S) \AND \bEE:\RSS\to\RSBS:S\mt \bEE(S).
\]
The regular $*$-semigroups $\PG(\PE)$ are defined in terms of the regular $*$-boset~$E$.  Analogously to $\IG(E)$ and $\RIG(E)$, we denote this by
\[
\IG^*(E) = \PG(\PE).
\]
Note that it is possible for regular $*$-bosets $E$ and $E'$ to have different underlying sets, but have identical
projection algebras $\PE = \PEd$. 
This would occur if $E$ and $E'$ were isomorphic, with  different underlying sets, but with the same set of projections.
In this case we would of course have $\IG^*(E) = \IG^*(E')$. 
This means that the assignment $E \mapsto \IG^*(E)$ is not injective.

One could get around this `problem' by instead defining $\IG^*(E)$ to be a copy of $\PG(\PE)$ in which we identify each $e\in E$ with $\ldb ee^*,e^*e\rdb\in\PG(\PE)$.
This would then reflect the situation of projection algebras, where $\PG(P)=\PG(P') \iff P=P'$, as follows from Theorem~\ref{thm:C_P}.

\section{Presentations}\label{sect:pres}

In this section we establish a number of presentations for the free (projection-generated) regular $*$-semigroup $\PG(P)$ over an arbitrary projection algebra $P$.  The first presentation (Theorem~\ref{thm:pres}) involves the generating set $P$.  The second and third (Theorems \ref{thm:presE} and \ref{thm:presE2}) are both in terms of the generating set $E=\EP$, and these highlight the connections between $\PG(P)$ and the free (idempotent-generated) semigroup $\IG(E)$ and the free regular (idempotent-generated) semigroup~$\RIG(E)$.

\subsection{Preliminaries on presentations}\label{ss:prepre}

We begin by establishing the notation we will use for presentations; for more details, see for example \cite[Section 1.6]{Howie1995}.  We also prove a general technical lemma that will be used in this section and later in the paper.

A \emph{congruence} on a semigroup $S$ is an equivalence relation $\si$ on $S$ that is compatible with multiplication, in the sense that $a\mr\si a'$ and $b\mr\si b'$ together imply $ab\mr\si a'b'$.  The quotient $S/\si$ is then a semigroup under the induced operation on $\si$-classes.  Given a semigroup homomorphism $\phi:S\to T$, the \emph{kernel} $\ker(\phi)=\set{(a,b)\in S\times S}{a\phi=b\phi}$ is a congruence on $S$, and the fundamental homomorphism theorem for semigroups states that $S/{\ker(\phi)} \cong \im(\phi)$.

For a set $X$, we denote by $X^+$ the free semigroup over $X$, which consists of all non-empty words over $X$, under concatenation.  For a set $R\sub X^+\times X^+$ of pairs of words, we write $R^\sharp$ for the congruence on $X^+$ generated by~$R$, i.e.~the least congruence on $X^+$ containing $R$.  We often write $[w]_R$ for the $R^\sharp$-class of $w\in X^+$.  We say a semigroup $S$ has \emph{presentation} $\pres XR$ if $S\cong X^+/R^\sharp$, i.e.~if there is a surjective semigroup homomorphism $X^+\to S$ with kernel $R^\sharp$.  At times we identify $\pres XR$ with the semigroup $X^+/R^\sharp$ itself.  The elements of $X$ and $R$ are called generators and (defining) relations, respectively.  A relation $(u,v)\in R$ is typically displayed as an equality: $u=v$.

\begin{lemma}\label{lem:pres}
Consider semigroups $S=\pres XR$ and $T=\pres YQ$, for which the following hold:
\ben
\item \label{lempres1} $X\sub Y$,
\item \label{lempres2} $R\sub Q^\sharp$,
\item \label{lempres3} every $y\in Y$ is $Q^\sharp$-equivalent to a word over $X$,
\item \label{lempres4} there is a morphism $\phi:Y^+\to S$ with $Q\sub\ker(\phi)$, and $x\phi=[x]_R$ for all $x\in X$.
\een
Then $S\cong T$.
\end{lemma}

\pf
By \ref{lempres1} we have a well defined morphism $\psi:X^+\to T:x\mt[x]_Q$.
By \ref{lempres2} and \ref{lempres4}, $\psi$ and~$\phi$ induce morphisms
\[
\Psi:S\to T: [x]_R\mt[x]_Q
\AND
\Phi:T\to S:[y]_Q\mt y\phi.
\]
By \ref{lempres3}, $T$ is generated by $\set{[x]_Q}{x\in X}$, and of course $S$ is generated by $\set{[x]_R}{x\in X}$.  Since
\[
[x]_R\Psi = [x]_Q \AND [x]_Q\Phi = x\phi = [x]_R \qquad\text{for all $x\in X$, by \ref{lempres4},}
\]
it follows that $\Psi$ and $\Phi$ 
are mutually inverse isomorphisms of~$S$ and $T$.
\epf

\subsection[Presentation over $P$]{\boldmath Presentation over $P$}\label{subsect:presP}

Here is the main result of this section:

\begin{thm}\label{thm:pres}
For any projection algebra $P$, the free regular $*$-semigroup $\PG(P)$ has presentation
\[
\PG(P) \cong \pres{X_P}{R_P},
\]
where $X_P = \set{x_p}{p\in P}$ is an alphabet in one-one correspondence with $P$, and where $R_P$ is the set of relations
\begin{align}
\tag*{\textsf{(R1)}} \label{R1} x_p^2 &= x_p &&\text{for all $p\in P$,}\\
\tag*{\textsf{(R2)}} \label{R2} (x_px_q)^2 &= x_px_q &&\text{for all $p,q\in P$,}\\
\tag*{\textsf{(R3)}} \label{R3} x_px_qx_p &= x_{q\th_p} &&\text{for all $p,q\in P$.}
\end{align}
\end{thm}

\begin{rem}
If $P=\bP(S)$ is the projection algebra of a regular $*$-semigroup $S$, then recall that $q\th_p=pqp$ for $p,q\in P$, where the product $pqp$ is taken in $S$.  Thus, relations of type \ref{R3} have the form $x_px_qx_p = x_{pqp}$ in this case.
\end{rem}

To prove Theorem \ref{thm:pres}, we require some technical lemmas.  But first, it is worth observing that relations \ref{R1}--\ref{R3} closely resemble projection algebra axioms \ref{P2}, \ref{P4} and \ref{P5}, i.e.~those that are stated purely in terms of the $\th$ maps.

For the rest of this subsection, we fix $P$, $X_P$ and $R_P$ as in Theorem \ref{thm:pres}.
We write ${\sim}=R_P^\sharp$ for the congruence on $X_P^+$ generated by relations \ref{R1}--\ref{R3}, and use $\sim_1$ to indicate equivalence by one or more applications of \ref{R1}, and similarly for $\sim_2$ and~$\sim_3$.

\begin{lemma}\label{lem:pres1}
If $p,q\in P$ are such that $q\leqF p$, then $x_px_q \sim x_{p'}x_q$ for some $p'\in P$ with $q\F p'\leq p$.
\end{lemma}

\pf
Let $p' = q\th_p\leq p$.  Combining $q\leqF p$ with \ref{PA1}, it follows that $q=p\th_q\F q\th_p=p'$.  We also have
\[
x_px_q \sim_2 x_px_qx_px_q \sim_3 x_{q\th_p}x_q = x_{p'}x_q.  \qedhere
\]
\epf

\begin{lemma}\label{lem:pres2}
For any $p_1,\ldots,p_k\in P$ we have $x_{p_1}\cdots x_{p_k} \sim x_{p_1'}\cdots x_{p_k'}$ for some $p_1',\ldots,p_k'\in P$ with $p_1'\F\cdots\F p_k'$ and $p_i'\leq p_i$ for all $i$.
\end{lemma}

\pf
We prove the lemma by induction on $k$.  The $k=1$ case being trivial, we assume $k\geq2$.  With $p_1''=p_2\th_{p_1}\leq p_1$ and $p_2''=p_1\th_{p_2}\leq p_2$ we have
\[
x_{p_1}x_{p_2} \sim_2 x_{p_1}x_{p_2}x_{p_1}x_{p_2}x_{p_1}x_{p_2} \sim_3 x_{p_2\th_{p_1}}x_{p_1\th_{p_2}} = x_{p_1''}x_{p_2''},
\]
and \ref{PA1} gives $p_1''\F p_2''$.  By induction, we have
\[
x_{p_2''}x_{p_3}\cdots x_{p_k} \sim x_{p_2'}x_{p_3'}\cdots x_{p_k'}
\]
for some $p_2',\ldots,p_k'\in P$ with $p_2'\F\cdots\F p_k'$, $p_2'\leq p_2''$ and $p_i'\leq p_i$ for $i=3,\ldots,k$.  Note that also $p_2'\leq p_2''\leq p_2$.  
Since $p_2'\leq p_2''\F p_1''$, \ref{PA2} gives $p_2'\leqF p_1''$.  It then follows from Lemma \ref{lem:pres1} that $x_{p_1''}x_{p_2'} \sim x_{p_1'}x_{p_2'}$ for some $p_1'\in P$ with $p_2'\F p_1' \leq p_1''$, and again we observe that $p_1'\leq p_1''\leq p_1$.  Putting everything together we have
\[
x_{p_1}x_{p_2}x_{p_3}\cdots x_{p_k} \sim x_{p_1''}x_{p_2''}x_{p_3}\cdots x_{p_k} \sim x_{p_1''}x_{p_2'}x_{p_3'}\cdots x_{p_k'} \sim x_{p_1'}x_{p_2'}x_{p_3'}\cdots x_{p_k'},
\]
with all conditions met.
\epf

Given a $P$-path $\p=(p_1,\ldots,p_k)\in\P=\P(P)$, we define the word
\[
w_\p = x_{p_1}\cdots x_{p_k} \in X_P^+.
\]
It follows from Lemma \ref{lem:pres2} that every word over $X_P$ is $\sim$-equivalent to some $w_\p$.  Using \ref{R1}, it is easy to see that
\begin{equation}\label{eq:xpq}
w_\p w_\q \sim w_{\p\circ\q} \qquad\text{for any $\p,\q\in\P$ with $\br(\p)=\bd(\q)$.}
\end{equation}
The next result refers to the congruence ${\ssim}=\Xi^\sharp$ on $\P$ from Definition \ref{defn:olC}.

\begin{lemma}\label{lem:pres3}
For any $\p,\q\in\P$, we have $\p\ssim\q \implies w_\p \sim w_\q$.
\end{lemma}

\pf
It suffices to assume that $\p$ and $\q$ differ by a single application of \ref{Om1}--\ref{Om3}, i.e.~that
\[
\p = \p'\circ\s\circ\p'' \AND \q = \p'\circ\t\circ\p'' \qquad\text{for some $\p',\p''\in\P$ and $(\s,\t)\in\Xi\cup\Xi^{-1}$.}
\]
Since $w_\p \sim w_{\p'}w_\s w_{\p''}$ and $w_\q \sim w_{\p'}w_\t w_{\p''}$ by \eqref{eq:xpq}, it is in fact enough to prove that
\[
w_\s\sim w_\t \qquad\text{for all $(\s,\t)\in\Xi$.}
\]
We consider the three forms the pair $(\s,\t)\in\Xi$ can take.

\pfitem{\ref{Om1}}  This follows immediately from \ref{R1}.

\pfitem{\ref{Om2}}  If $\s=(p,q,p)$ and $\t=(p)$ for some $(p,q)\in{\F}$, then $w_\s = x_px_qx_p \sim_3 x_{q\th_p} = x_p = w_\t$.

\pfitem{\ref{Om3}}  Finally, suppose $\s=\lam(e,p,f)=(e,e\th_p,f)$ and $\t=\rho(e,p,f)=(e,f\th_p,f)$ for some $p\in P$, and some $p$-linked pair $(e,f)$.  Then
\[
w_\s = x_ex_{e\th_p}x_f \sim_3 x_ex_px_ex_px_f \sim_2 x_ex_px_f \sim_2 x_ex_px_fx_px_f \sim_3 x_ex_{f\th_p}x_f = w_\t.  \qedhere
\]
\epf

We can now tie together the loose ends.

\pf[{\bf Proof of Theorem \ref{thm:pres}.}]
Define the homomorphism
\[
\Psi:X_P^+\to\PG(P) \BY x_p\Psi = p =\ldb p\rdb  \qquad\text{for $p\in P$.}
\]
To see that $\Psi$ is surjective, let $\c\in\PG(P)$, so that $\c=\ldb \p\rdb $ for some $\p=(p_1,\ldots,p_k)\in\P$.  
We claim that
\[
\ldb p_1,\ldots,p_i\rdb = p_1\pr\cdots\pr p_i \qquad\text{for all $1\leq i\leq k$.}
\]
Indeed, the $i=1$ case is clear, and if $2\leq i\leq k$ then
\begin{align*}
p_1\pr\cdots\pr p_i &= \ldb p_1,\ldots,p_{i-1}\rdb \pr p_i &&\text{by induction}\\
&= \ldb p_1,\ldots,p_{i-1}\rdb\rest_{p_{i-1}'} \circ\ldb p_{i-1}',p_i'\rdb \circ {}_{p_i'}\corest p_i &&\text{where $p_{i-1}'=p_i\th_{p_{i-1}}$ and $p_i'=p_{i-1}\th_{p_i}$}\\
&= \ldb p_1,\ldots,p_{i-1}\rdb \circ\ldb p_{i-1},p_i\rdb \circ p_i &&\text{since $p_{i-1}'=p_{i-1}$ and $p_i'=p_i$, as $p_{i-1}\F p_i$}\\
&= \ldb p_1,\ldots,p_{i-1},p_i\rdb,
\end{align*}
proving the claim. It then follows that
\begin{equation}\label{eq:xpPsi}
\c = \ldb p_1,\ldots,p_k\rdb = p_1\pr\cdots\pr p_k = (x_{p_1}\cdots x_{p_k})\Psi = w_\p\Psi,
\end{equation}
completing the proof that $\Psi$ is surjective.

Next, we note that $R_P\sub\ker(\Psi)$, meaning that $u\Psi=v\Psi$  (in $\PG(P)$) for all $(u,v)\in R_P$.  Indeed, this is clear when $(u,v)$ has type \ref{R1}, and follows from \ref{RS2} or \ref{CP4} for type \ref{R2} or~\ref{R3}.

It remains to show that $\ker(\Psi)\sub R_P^\sharp$.  To do so, fix some $(u,v)\in\ker(\Psi)$, so that $u,v\in X_P^+$ and $u\Psi=v\Psi$; we must show that $u\sim v$ (recall that we write $\sim$ for $R_P^\sharp$).  By Lemma \ref{lem:pres2}, we have $u\sim w_\p$ and $v\sim w_\q$ for some~$\p,\q\in\P$.  Using \eqref{eq:xpPsi}, and remembering that ${\sim}\sub\ker(\Psi)$, we have
\[
\ldb \p\rdb  = w_\p\Psi = u\Psi = v\Psi = w_\q\Psi = \ldb \q\rdb ,
\]
meaning that $\p\ssim\q$.  But then $w_\p\sim w_\q$ by Lemma \ref{lem:pres3}, so $u\sim w_\p\sim w_\q\sim v$, as required.
\epf

\subsection[Presentation as a quotient of $\IG(E)$]{\boldmath Presentation as a quotient of $\IG(E)$}\label{subsect:presE}

Consider a boset $E$, and recall that a product $ef$ is defined in $E$ precisely when $(e,f)\in\BP(E)$ is a basic pair, i.e.~when $\{ef,fe\}\cap\{e,f\}\not=\es$.  The \emph{free (idempotent-generated) semigroup over~$E$} has presentation
\begin{equation}\label{eq:IGE}
\IG(E) = \pres{X_E}{x_ex_f=x_{ef}\ \text{for all $(e,f)\in\BP(E)$}},
\end{equation}
where here $X_E = \set{x_e}{e\in E}$ is an alphabet in one-one correspondence with $E$.  The boset of~$\IG(E)$ is isomorphic to $E$ \cite{Easdown1985}, and consists of all equivalence classes of letters $x_e$ ($e\in E$).

Now consider a projection algebra $P$, and let $E=\EP=\bEE(\PG(P))$ be the (regular $*$-) boset of $\PG(P)$.  Since $\PG(P)$ is a semigroup with biordered set $E=\bE(P)$, general theory \cite[Theorem~3.3]{Easdown1985} tells us that the mapping $x_e\mt e$ ($e\in E$) induces a surmorphism ${\IG(E)\to\PG(P)}$.  Thus, we know in advance that there exists a presentation for $\PG(P)$ extending the above presentation for~$\IG(E)$ by means of additional relations.  Theorem \ref{thm:presE} below gives an explicit such presentation, with additional relations $x_px_q=x_{pq}$ for projections $p,q\in P(\sub E)$. 
These additional relations can also be viewed as a generating set for the kernel of the surmomorphism $\IG(E)\to\PG(P)$.  
Note that the product $pq$ might not exist in the boset $E$, but it certainly exists in the semigroup $\PG(P)$, and is an idempotent, and hence a well-defined element of $E$; in fact,~$pq$ is the element $e(q,p)$ from Lemma \ref{lem:epq}.  
For simplicity, we will denote the product in $\PG(P)$ by juxtaposition instead of $\pr$ throughout this subsection.

\begin{thm}\label{thm:presE}
For any projection algebra $P$, the free regular $*$-semigroup $\PG(P)$ has presentation
\[
\PG(P) \cong \pres{X_E}{R_E},
\]
where $X_E = \set{x_e}{e\in E}$ is an alphabet in one-one correspondence with $E=\EP$, and where~$R_E$ is the set of relations
\begin{align}
\tag*{\textsf{(R1)$'$}} \label{R1'} x_ex_f &= x_{ef} &&\text{for all $(e,f)\in\BP(E)$,}\\
\tag*{\textsf{(R2)$'$}} \label{R2'} x_px_q &= x_{pq}  &&\text{for all $p,q\in P$.}
\end{align}
\end{thm}

\pf
By Theorem \ref{thm:pres} we have $\PG(P)\cong\pres{X_P}{R_P}$.  Thus, we can prove the current theorem by applying Lemma \ref{lem:pres} with $S=\pres{X_P}{R_P}$ and ${T=\pres{X_E}{R_E}}$.  To do so, we must show that:
\ben
\item \label{presE1} $X_P\sub X_E$,
\item \label{presE2} $R_P\sub R_E^\sharp$,
\item \label{presE3} every $x_e$ ($e\in E$) is $\sim'$-equivalent to a word over $X_P$, where ${\sim'}=R_E^\sharp$, and
\item \label{presE4} there is a morphism $\phi:X_E^+\to \pres{X_P}{R_P}$ such that $R_E\sub\ker(\phi)$, and $x_p\phi=[x_p]$ for all $p\in P$, where we write $[w]$ for the $R_P^\sharp$-class of $w\in X_P^+$.
\een
Item \ref{presE1} is clear.  For \ref{presE2}, we check the relations from $R_P$ in turn.  We use~$\sim'_1$ and~$\sim'_2$ to denote equivalence via \ref{R1'} or \ref{R2'}.

\pfitem{\ref{R1}}  This is contained in \ref{R1'} (and in \ref{R2'}).  

\pfitem{\ref{R2}}  If $p,q\in P$, then $x_px_q \sim'_2 x_{pq} \sim'_1 x_{pq}^2 \sim'_2(x_px_q)^2$.

\pfitem{\ref{R3}}  If $p,q\in P$, then $pq\in E$, and $(p,pq)$ is a basic pair in $E$.  It follows that \ref{R1'} contains the relation $x_{pq}x_p = x_{pqp}$.  But then $x_px_qx_p \sim'_2 x_{pq}x_p \sim'_1 x_{pqp} = x_{q\th_p}$.

\aftercases For \ref{presE3}, fix some $e\in E$.  Since $e=ee^*e^*e$, with $ee^*,e^*e\in P$, we have
\[
x_e = x_{ee^*e^*e} \sim'_2 x_{ee^*}x_{e^*e} \in X_P^+.
\]
For \ref{presE4}, we first define a morphism
\[
\psi:X_E^+\to\PG(P):x_e\mt e.
\]
To see that $R_E\sub\ker(\psi)$, fix some $e,f\in E$.  If $ef\in E$, then $\psi$ maps both $x_ex_f$ and $x_{ef}$ to $ef$, and it follows that relations \ref{R1'} and \ref{R2'} are both preserved.  Composing $\psi:X_E^+\to\PG(P)$ with the isomorphism $\PG(P)\to\pres{X_P}{R_P}:p\mt[x_p]$ gives $\phi:X_E^+\to\pres{X_P}{R_P}$ with the required properties.
\epf

\begin{rem}
Relations \ref{R2'} reflect the fact that the product of two projections is always an idempotent in a regular $*$-semigroup.  In fact, any idempotent is the product of two $\F$-related projections by \ref{RS3}.
Consequently, one  might wonder if we could replace \ref{R2'} by the subset consisting of the relations
$x_px_q = x_{pq}$ for $(p,q)\in{\F}$,
and still define $\PG(P)$.
It turns out that this is not possible in general.  For example, consider the three element semilattice $S=\{0,e,f\}$ with $ef=0$, and note that $P=E=S$ and ${\F}=\De_P$ in this case.    
As \ref{R2'} is a copy of the multiplication table of~$S$ (which is the case for any semilattice), it follows that $\PG(P)\cong S$. 
Since ${\F}=\De_P$, the relations in \ref{R2'} arising from friendly pairs are just $x_p^2=x_p$ for $p\in E$,  which are already contained in~\ref{R1'}. So if we reduce the presentation in this way,
we actually arrive at the free idempotent-generated semigroup $\IG(E)$. 
As  shown in \cite[Example 2]{BMM2009}, $\IG(E)$ is infinite (and non-regular), and so certainly not isomorphic to $ \PG(P)$.
In the next subsection we will see that  if instead of $\IG(E)$ we start with a presentation for $\RIG(E)$, then the relations from \ref{R2'} arising from friendly pairs are indeed sufficient to define $\PG(P)$.
\end{rem}

\begin{rem}\label{rem:IGE_P}
Consider a projection algebra $P$, and its associated boset $E=\bE(P)$.  We have now given presentations for $\PG(P)$ in terms of (copies of) the generating sets $P$ (Theorem~\ref{thm:pres}) and~$E$ (Theorem \ref{thm:presE}).  On the other hand, $\IG(E)$ is defined in terms of a presentation with generating set $E$ (see \eqref{eq:IGE}), and one might wonder if there is a presentation utilising the generating set~$P$.  This is not the case, however, as~$P$ need not generate $\IG(E)$ in general.  
We will give a concrete instance of this in Example \ref{eg:2x2}.
\end{rem}

\subsection[Presentation as a quotient of $\RIG(E)$]{\boldmath Presentation as a quotient of $\RIG(E)$}
\label{subsect:RIG}

Consider again a projection algebra $P$, and its associated boset ${E=\EP=\bEE(\PG(P))}$.  Since~$E$ is regular (i.e.~the boset of a regular semigroup, namely $\PG(P)$), we also have the \emph{free regular (idempotent-generated) semigroup} $\RIG(E)$.  This was defined by Nambooripad in \cite{Nambooripad1979} using his groupoid machinery, and in \cite{Pastijn1980} by means of the presentation
\begin{align*}
\RIG(E) = \la X_E : {} &x_ex_f=x_{ef}\ \ \ \hspace{0.65cm} \text{for all $(e,f)\in\BP(E)$,}\\
& x_ex_gx_f = x_ex_f\ \ \ \text{for all $e,f\in E$ and $g\in \SS(e,f)$} \ra.
\end{align*}
Here $\SS(e,f)$ is the sandwich set of the idempotents $e,f\in E$, defined in Subsection \ref{subsect:Eprelim}.  Note that the characterisation of $\SS(e,f)$ in \eqref{eq:SSef} applies to $S=\PG(P)$, as $\PG(P)$ is a regular semigroup with boset $E$.
The above presentation shows that $\RIG(E)$ is a quotient of $\IG(E)$.  In turn, our next result shows that $\PG(P)$ is a quotient of $\RIG(E)$, with the additional relations generating the kernel of the canonical surmorphism $\RIG(E)\to\PG(P)$.

\begin{thm}\label{thm:presE2}
For any projection algebra $P$, the free regular $*$-semigroup $\PG(P)$ has presentation
\[
\PG(P) \cong \pres{X_E}{R_E'},
\]
where $X_E = \set{x_e}{e\in E}$ is an alphabet in one-one correspondence with $E=\EP$, and where~$R_E'$ is the set of relations
\begin{align}
\tag*{\textsf{(R1)$''$}} \label{R1''} x_ex_f &= x_{ef} &&\text{for all $(e,f)\in\BP(E)$,}\\
\tag*{\textsf{(R2)$''$}} \label{R2''} x_ex_f &= x_ex_gx_f  &&\text{for all $e,f\in E$ and $g\in S(e,f)$,}\\
\tag*{\textsf{(R3)$''$}} \label{R3''} x_px_q &= x_{pq}  &&\text{for all $(p,q)\in {\F}$.}
\end{align}
\end{thm}

\pf
We begin with the presentation $\pres{X_E}{R_E} = \pres{X_E}{\text{\ref{R1'},\ \ref{R2'}}}$ from Theorem \ref{thm:presE}, via the mapping
\[
\psi:X_E^+\to\PG(P):x_e\mt e.
\]
Since $\psi$ maps $x_ex_gx_f$ and $x_ex_f$ both to $egf=ef$ for $g\in \SS(e,f)$, we can add relations \ref{R2''} to the presentation.  Noting that \ref{R1'} and \ref{R1''} are the same sets of relations, the presentation has now become
\[
\pres{X_E}{\text{\ref{R1''},\ \ref{R2''},\ \ref{R2'}}}.
\]
Since \ref{R3''} is contained in \ref{R2'}, we can complete the proof by showing that each relation in~\ref{R2'} is implied by those in $R_E'$.  To do so, let $p,q\in P$ be arbitrary, and let $p'=q\th_p$ and $q'=p\th_q$, so that $p'\F q'$ and $p'q'=pq$.  We then calculate (again writing $\sim_1''$ for equivalence by \ref{R1''}, and so on)
\begin{align*}
x_px_q &\sim''_2 x_px_{qp}x_q &&\text{as $qp\in \SS(p,q)$}\\
&\sim''_1 x_px_{qp}x_{qp}x_q \\
&\sim''_1 x_{p\cdot qp}x_{qp\cdot q} &&\text{as $(qp,p)$ and $(q,qp)$ are basic pairs}\\
&= x_{p'}x_{q'} \sim''_3 x_{p'q'} = x_{pq}. &&\qedhere
\end{align*}
\epf

\section{Adjacency semigroups and bridging path semigroups}\label{sect:AGa}

We now turn to explicit examples of free projection-generated regular $*$-semigroups, starting with those arising from adjacency semigroups, as introduced in Subsection \ref{subsect:AGa}.

Let $\Ga=(P,E)$ be a symmetric, reflexive digraph, and let $A_\Ga$ be its adjacency semigroup.
We keep the notation of Subsection \ref{subsect:AGa}, including the projection algebra $P_0=P\cup\{0\}$, whose operations were given in \eqref{eq:thpAGa}.  We now consider the structure of the free regular $*$-semigroup~$\PGPz$.

The path category $\P=\P(P_0)$ consists of tuples of the form $(0,\ldots,0)$, and $(p_1,\ldots,p_k)$ where each $(p_i,p_{i+1})\in E(\sub{\F})$; 
the latter are simply the paths in $\Ga$ in the usual graph-theoretical sense, with repeated vertices allowed as in Subsection \ref{subsect:CP}.

As in Remark \ref{rem:conf}, any non-zero path $\p\in\P$ is $\approx$-equivalent to a unique \emph{reduced} path $\ol\p=(p_1,\ldots,p_k)$, where each $p_i$ is distinct from $p_{i+1}$ (if $i\leq k-1$) and from $p_{i+2}$ (if $i\leq k-2$).  By identifying a non-zero chain ($\approx$-class of a path) with the unique reduced path it contains, we may identify the chain groupoid $\C=\C(P_0)$ with the set of reduced paths in $\Ga$, along with $0=[0]$.

Using \eqref{eq:thpAGa}, we see that $(e,f)$ is $p$-linked if and only if $e=f=0$ or $(e,p)$ and $(p,f)$ are both edges of $\Ga$.  In these respective cases, we have
\[
\lam(e,p,f)=\rho(e,p,f) = (0,0,0) \OR \lam(e,p,f)=\rho(e,p,f) = (e,p,f).
\]
It follows that the congruences $\approx$ and $\ssim$ are equal, and so $\ol\C=\C$, and $\ldb\p\rdb=[\p]$ for all $\p\in\P$.

The free regular $*$-semigroup $\PGPz$ can therefore be viewed as follows:
\bit
\item The elements are $0$ and the reduced paths in $\Ga$.
\item The product of reduced paths $\p=(p_1,\ldots,p_k)$ and $\q=(q_1,\ldots,q_l)$ is given by
\[
\p\pr\q = \begin{cases}
\ol{\p\op\q} &\text{if $(p_k,q_1)\in E$ is an edge of $\Ga$}\\
0 &\text{otherwise.}
\end{cases}
\]
Here $\op$ denotes concatenation, so $\p\op\q=(p_1,\ldots,p_k,q_1,\ldots,q_l)$.
\item The involution is given by reversal of paths.
\item The non-zero projections are the empty paths $p=(p)$, which are in one-one correspondence with the vertices of $\Ga$.
\item The remaining non-zero idempotents are the non-loop edges $p\pr q=(p,q)$ of $\Ga$.
\eit
The authors have not seen these specific semigroups in the literature, but we note that they are closely related to the \emph{graph inverse semigroups} introduced in \cite{AH1975}, in which $\p\pr\q$ is only non-zero when $p_k=q_1$; there are some other differences as well, including the fact that the digraphs of \cite{AH1975} are not assumed to be symmetric or reflexive.
In our case $\p$ and $\q$ can also be composed if there is an edge $p_k\to q_1$, which `bridges' the end of $\p$ and the start of $\q$.  We therefore call $\PG(\bPP(A_\Ga))$ the \emph{bridging path semigroup of $\Ga$}, and denote it by $B_\Ga$.  
Note that bridging path semigroups can be defined starting from an arbitrary digraph, i.e.~without assuming symmetry and reflexivity a priori. The properties of such a semigroup would then depend on the properties of the graph~$\Gamma$, and this may be an interesting direction for study. In particular, one can verify that $B_\Gamma$ is a regular $*$-semigroup with projections $P_0$ and the involution given by reversal of paths if and only if $\Gamma$ is symmetric and reflexive.

In the special case that $\Ga$ is the complete digraph,  the adjacency semigroup~$A_\Ga$ is simply the square band $B_P=P\times P$ with a zero adjoined.  This band $B_P$ is a regular $*$-semigroup, with operations $(p,q)(r,s) = (p,s)$ and $(p,q)^*=(q,p)$.  Every element of $B_P$ is an idempotent; the projections are $p=(p,p)$; and the projection algebra $\bPP(B_P)=P$ has a trivial structure, in the sense that the $\th_p$ operations are all constant maps.  
The above analysis shows that the semigroup~$\PG(P)$ consists of all reduced tuples, with product $\p\pr\q=\ol{\p\op\q}$.  
In the next two examples we give some more explicit details concerning the square bands $B_P$, and the corresponding free regular $*$-semigroups $\PG(P)$, in the cases that $|P|\leq3$.

\begin{eg}\label{eg:2x2}
If $|P|\leq2$ then $\PG(P)=B_P$ is finite.  By contrast, the corresponding free idempotent-generated semigroup $\IG(E)$ is infinite when $|P|=2$; this is folklore and can be deduced from \cite[Theorem~5]{GR2012}.  Here we write $E=B_P=\bEE(B_P)\cong\bE(P)$.  
More specifically, it follows from from \cite[Theorem~5]{GR2012} that $\IG(E)$ is isomorphic to the $2\times2$ Rees matrix semigroup $S$ over an infinite cyclic group $H=\la a\ra$, with respect to the sandwich matrix~$\big(\begin{smallmatrix} 1&1\\ 1& a\end{smallmatrix}\big)$.  In particular, this lets us prove the claim from Remark \ref{rem:IGE_P}, namely that $\IG(E)$ is not generated by (its canonical copy of) $P$.

To see this, write $P=\{p,q\}$.  The above-mentioned isomorphism $\IG(E)\to S$ maps (the equivalence classes of) the letters $x_p$ and $x_q$ to the idempotents $\ol p=(1,1,1)$ and $\ol q=(2,a^{-1},2)$, respectively.  (The other two idempotents of $S$ are $(1,1,2)$ and $(2,1,1)$.)  It is easy to see that every element of $\la\ol p,\ol q\ra$ has the form $(i,a^{-m},j)$ for some $i,j\in\{1,2\}$ and $m\geq0$, so indeed~$\la\ol p,\ol q\ra\not=S$.
\end{eg}

\begin{eg}\label{eg:3x3}
If $|P|\geq3$ then $\PG(P)$ is infinite, and in particular $\PG(P)\not=B_P$.
It is instructive to consider the case where $P=\{p,q,r\}$ has size $3$.  For $s,t\in P$, let $H_{s,t}$ be the set of all reduced paths from $s$ to $t$, so that $\PG(P)=\bigsqcup_{s,t\in P}H_{s,t}$.  Note for example that
\[
H_{p,p} = \{p\} \cup \set{(p,q,r,p)^k,(p,r,q,p)^k}{k=1,2,\ldots}.
\]
It is easy to check that $H_{p,p}$ is isomorphic to the infinite cyclic group $H=\la a\ra$.  The identity of $H_{p,p}$ is $p$, and $(p,q,r,p)^k$ and $(p,r,q,p)^k$ are inverses of each other.  An analogous argument shows that each $H_{s,t}\cong H$.  It follows from general structure theory (see \cite[Chapter 3]{Howie1995}) that $\PG(P)$ is a $3\times3$ Rees matrix semigroup over $H$, with respect to the sandwich matrix
$\biggl(\begin{smallmatrix} 1&1&1\\ 1&1&a\\ 1& a^{-1} &1\end{smallmatrix}\biggr)$.
By way of comparison, $\IG(\bEE(B_P))=\RIG(\bEE(B_P))$ is a $3\times3$ Rees matrix semigroup over the free group of rank $4$, again following from \cite[Theorem~5]{GR2012}.
\end{eg}

The structure of an arbitrary bridging path semigroup $B_\Ga = \PG(\bP(A_\Ga))$ can be similarly described in terms of Rees $0$-matrix semigroups over free groups.  The description requires an analysis of the maximal subgroups of free regular $*$-semigroups, which is the subject of the forthcoming paper~\cite{Paper4}.

Our final example in this section is an extension $P'$ of the projection algebra $P$ of the $3\times3$ rectangular band from Example \ref{eg:3x3} by a single extra projection.  It was introduced in \cite{EPA2024} (where it had an identity adjoined), and was originally discovered by Michael Kinyon.
We note that Kinyon's projection algebra is not the projection algebra of an adjacency semigroup.  It illustrates the fact that it is possible for a projection algebra $P$ with $\PG(P)$ infinite to be contained in a projection algebra $P'$ with $\PG(P')$ finite.

\begin{eg}\label{eg:Kinyon}
Let $P'=\{p,q,r,e\}$ be the projection algebra with operations
\[
\th_p = \trans{p&q&r&e\\p&p&p&p} \COMMA
\th_q = \trans{p&q&r&e\\q&q&q&q} \COMMA
\th_r = \trans{p&q&r&e\\r&r&r&r} \AND
\th_e = \trans{p&q&r&e\\p&q&q&e} .
\]
Note that $p,q,r$ are all $\F$-related, but $e$ is $\F$-related only to itself.  As in Example \ref{eg:3x3}, it follows that 
\[
\C=\C(P') = \{e\} \cup \bigsqcup_{s,t\in\{p,q,r\}} H_{s,t}.
\]
However, $\ol\C$ is a proper quotient of $\C$ here, as there is a non-trivial linked pair.  Specifically,~$(p,r)$ is $e$-linked, and we have
\[
\lam(p,e,r) = (p,p\th_e,r) = (p,p,r)  \AND \rho(p,e,r) = (p,r\th_e,r) = (p,q,r).
\]
Consequently, $\ldb p,q,r\rdb = \ldb p,p,r\rdb = \ldb p,r\rdb$ in $\PG(P')$.  It follows from this that $\ldb s,t,u\rdb=\ldb s,u\rdb$ for distinct $s,t,u\in\{p,q,r\}$.  For example,
\[
\ldb p,r,q\rdb = \ldb p,q,r,q\rdb = \ldb p,q\rdb \AND \ldb q,p,r\rdb = \ldb q,p,q,r\rdb = \ldb q,r\rdb.
\]
The other three cases are obtained by inverting the three already considered.  This all shows that $\PG(P')$ contains exactly ten elements:
\bit
\item the projections $e,p,q,r$, and
\item the remaining idempotents $s\pr t = \ldb s,t\rdb$, for distinct $s,t\in\{p,q,r\}$.
\eit
The entire multiplication table of $\PG(P')$ can be obtained from the fact that $\PG(P')\sm\{e\} = \la p,q,r\ra$ is a $3\times3$ rectangular band, together with the rules
\[
p\pr e = p \COMMA q\pr e = q \AND r\pr e = \ldb r,q\rdb.
\]
\end{eg}

\section{Temperley--Lieb monoids}\label{sect:TLn}

In Section \ref{sect:AGa} we saw examples of naturally occurring (albeit very small) regular $*$-semigroups~$S$ that were isomorphic to their associated free regular $*$-semigroup $\PG(\bPP(S))$.
In this section, we present a decidedly non-trivial family of semigroups displaying this phenomenon, namely the Temperley--Lieb monoids, introduced in Subsection \ref{subsect:DM}.

\begin{thm}\label{thm:TLnPGP}
Let $P=\bPP(\TL_n)$ be the projection algebra of a finite Temperley--Lieb monoid~$\TL_n$.  Then~$\PG(P)\cong\TL_n$.  
\end{thm}

\pf
We begin with the presentations $\TL_n\cong\pres{X_T}{R_T}$ and $\PG(P)\cong\pres{X_P}{R_P}$ from Theorems \ref{thm:TLn} and \ref{thm:pres}, and apply Lemma \ref{lem:pres} to show that they are equivalent.  For this, we first need to convert the monoid presentation $\pres{X_T}{R_T}$ into a semigroup presentation $\pres {X_T'}{R_T'}$ with $X_T'\sub X_P$, which we do as follows.
First, we identify each generator $t_i\in X_T$ with $x_{\tau_i}\in X_P$, noting that they both represent the element $\tau_i$ of $\TL_n$; see \eqref{eq:taui}.
Then we define $X_T'=\{e\}\cup X_T$, where here $e=x_1$ represents the identity $1$ of $\TL_n$.
Finally we add new relations that ensure~$e$ acts as the identity. The resulting set $R_T'$ of relations is:
\begin{align}
\tag*{\textsf{(T1)}} \label{U1} t_i^2 &= t_i &&\text{for all $i$}\\
\tag*{\textsf{(T2)}} \label{U2} t_it_j &= t_jt_i &&\text{if $|i-j|>1$}\\
\tag*{\textsf{(T3)}} \label{U3} t_it_jt_i &= t_i &&\text{if $|i-j|=1$}\\
\tag*{\textsf{(T4)}} \label{U4} e^2 &= e \\
\tag*{\textsf{(T5)}} \label{U5} et_i=t_ie &= t_i &&\text{for all $i$.}
\end{align}
It is clear that $\TL_n\cong\pres {X_T'}{R_T'}$.  We now work towards applying Lemma \ref{lem:pres}, with ${S=\pres {X_T'}{R_T'}}$ and $T=\pres{X_P}{R_P}$.  For this we need to show that:
\ben
\item \label{presTL1} $X_T'\sub X_P$,
\item \label{presTL2} $R_T'\sub R_P^\sharp$,
\item \label{presTL3} every $x_p$ ($p\in P$) is $\sim$-equivalent to a word over $X_T'$, where again ${\sim}=R_P^\sharp$, and
\item \label{presTL4} there is a morphism $\phi:X_P^+\to \pres {X_T'}{R_T'}$ such that $R_P\sub\ker(\phi)$, and $x\phi=[x]_{R_T'}$ for all $x\in X_T'$.
\een
We have already noted that \ref{presTL1} holds.  For \ref{presTL2}, we consider the relations from $R$ in turn.  Again we write $\sim_1$ to denote equivalence via \ref{R1}, and so on.

\pfitem{\ref{U1} and \ref{U4}}  These are contained in \ref{R1}.  

\pfitem{\ref{U2}}  Fix $i,j$ with $|i-j|>1$, and define the projection $p = \tau_i\tau_j\tau_i = \tau_j\tau_i\tau_j (=\tau_i\tau_j)$.  Then
\[
t_it_j = x_{\tau_i}x_{\tau_j} \sim_2 x_{\tau_i}x_{\tau_j}x_{\tau_i}x_{\tau_j}x_{\tau_i}x_{\tau_j} \sim_3 x_{\tau_i\tau_j\tau_i}x_{\tau_j\tau_i\tau_j} =x_px_p \sim_1 x_p \ANDSIm t_jt_i \sim x_p.
\]

\pfitem{\ref{U3}}  We have $t_it_jt_i = x_{\tau_i}x_{\tau_j}x_{\tau_i} \sim_3 x_{\tau_i\tau_j\tau_i} = x_{\tau_i} = t_i$.

\pfitem{\ref{U5}}  Writing $p=\tau_i$ for simplicity, we have
\[
et_i = x_1x_p \sim_2 x_1x_px_1x_p \sim_3 x_{1p1}x_p = x_px_p \sim_1 x_p = t_i \ANDSIM t_ie\sim t_i.
\]

For \ref{presTL3} we first note that $x_1=e\in X_T'$.  For $p\not=1$ we have $p = \tau_{i_1}\cdots\tau_{i_k}$ (in~$\TL_n$) for some $i_1,\ldots,i_k$, and then
\[
p = p^*p = \tau_{i_k}\cdots\tau_{i_1}\tau_{i_1}\cdots\tau_{i_k} = \tau_{i_k}\cdots\tau_{i_2}\tau_{i_1}\tau_{i_2}\cdots\tau_{i_k} = \tau_{i_1} \th_{\tau_{i_2}}\cdots\th_{\tau_{i_k}}.
\]
It follows that $x_p = x_{\tau_{i_1} \th_{\tau_{i_2}}\cdots\th_{\tau_{i_k}}} \sim_3 x_{\tau_{i_k}}\cdots x_{\tau_{i_2}}x_{\tau_{i_1}}x_{\tau_{i_2}}\cdots x_{\tau_{i_k}} = t_{i_k}\cdots t_{i_2} t_{i_1}t_{i_2}\cdots t_{i_k}$.

Finally, for \ref{presTL4} we begin by defining a morphism
\[
\psi:X_P^+\to\TL_n:x_p\mt p.
\]
We then have $R_P\sub\ker(\psi)$ since $P=\bPP(\TL_n)$.  We then obtain the desired $\phi:X_P^+\to\pres {X_T'}{R_T'}$ by composing $\psi:X_P^+\to\TL_n$ with the canonical isomorphism $\TL_n\to\pres {X_T'}{R_T'}$.
\epf

It is natural to wonder about the relationships between other diagram monoids and their associated free projection-generated regular $*$-semigroups.  These relationships tend to be more complicated than the case of $\TL_n$, as glimpsed in Example \ref{eg:M4} below, and treated in detail for the partition monoid~$\PP_n$ in the forthcoming paper~\cite{Paper4}.

\section{Topological interpretation}\label{sect:top}

In this section we provide an alternative, topological interpretation of the groupoids $\C=\C(P)$ and $\ol\C=\ol\C(P)$ associated to a projection algebra $P$.  Specifically, we can view $\C$ as the fundamental groupoid of a natural graph $G_P$ built from the $\F$-relation of $P$, and $\ol\C$ as the fundamental groupoid of either of two $2$-complexes $K_P$ and $K_P'$.  These  complexes both have $G_P$ as their $1$-skeleton, and their $2$-cells are induced by linked pairs of projections.  The complex $K_P$ is defined in terms of all linked pairs, and $K_P'$ in terms of a special subset that generates the others; in particular~$K_P'$ is simplicial.  Maximal subgroups of the semigroup $\PG(P)$ coincide up to isomorphism with fundamental groups of either of these two complexes.

\subsection{Special and degenerate linked pairs}\label{subsect:degen}

Consider a $p$-linked pair $(e,f)$ in a projection algebra $P$.  The projections $e,f,e\th_p,f\th_p$ (which are used to define the paths $\lam(e,p,f)$ and $\rho(e,p,f)$) need not all be distinct, and an important special case occurs when $e=e\th_p$ or $f=f\th_p$ (i.e.~when $e\leq p$ or $f\leq p$).  In this case, we say $(e,f)$ is a \emph{special $p$-linked pair}.  


Let $\Xi'$ be the subset of $\Xi$ consisting of all pairs $(\s,\t)\in\P\times\P$ of the form \ref{Om1} and \ref{Om2}, as well as:
\begin{enumerate}[label=\textup{\textsf{($\mathsf{\Om}$\arabic*)$'$}},leftmargin=10mm]\addtocounter{enumi}{2}
\item \label{Om3'} $\s=\lam(e,p,f)$ and $\t=\rho(e,p,f)$, for some $p\in P$, and some special $p$-linked pair $(e,f)$,
\end{enumerate}
and let ${\ssim'}=\Xi'^\sharp$ be the congruence on $\P$ generated by $\Xi'$.

\begin{prop}\label{prop:ssim'}
We have ${\ssim'}={\ssim}$.
\end{prop}

\pf
We just need to show that $\s\ssim'\t$ whenever $(\s,\t)$ is a pair of type \ref{Om3}.  So suppose $\s=\lam(e,p,f)$ and $\t=\rho(e,p,f)$ for some $p$-linked pair $(e,f)$, and let $e'=e\th_p$ and $f'=f\th_p$.  Then one can show that $(e,f')$ and $(e',f)$ are special $p$-linked pairs (with $f'\leq p$ in the first, and $e'\leq p$ in the second), and we have
\begin{align}
\nonumber \lam(e,p,f') &= (e,e',f') , &  \lam(e',p,f) &= (e',e',f) \ssim' (e',f), \\
\label{eq:e'f'} \rho(e,p,f') &= (e,f',f') \ssim' (e,f') , &    \rho(e',p,f) &= (e',f',f).
\end{align}
Since $\Xi'$ contains the pairs $(\lam(e,p,f'),\rho(e,p,f'))$ and $(\lam(e',p,f),\rho(e',p,f))$, it follows that $(e,e',f')\ssim'(e,f')$ and $(e',f)\ssim'(e',f',f)$.  But then
\[
\s = \lam(e,p,f) = (e,e',f) \ssim' (e,e',f',f) \ssim' (e,f',f) = \rho(e,p,f) = \t.  \qedhere
\]
\epf

Another case in which we do not need to include a pair $(\s,\t)$ of type \ref{Om3} is when we already have $\s\approx\t$.  

\begin{defn}\label{defn:degen}
We say a $p$-linked pair $(e,f)$ is \emph{degenerate} if $\lam(e,p,f) \approx \rho(e,p,f)$.
\end{defn}

The next result characterises such degenerate pairs at the level of the projection algebra structure.

\begin{prop}\label{prop:degen}
A $p$-linked pair $(e,f)$ is degenerate if and only if $e\th_p=f\th_p$ or $e,f\leq p$.
\end{prop}

\pf
Throughout the proof we write $e'=e\th_p$ and $f'=f\th_p$, and also $\lam=\lam(e,p,f)=(e,e',f)$ and $\rho=\rho(e,p,f)=(e,f',f)$.

\pfitem{($\Leftarrow$)}  If $e'=f'$, then $\lam=\rho$.  If $e,f\leq p$, then $e=e'$ and $f=f'$, so that
\[
\lam=(e,e,f) \approx (e,f) \approx (e,f,f) = \rho.
\]

\pfitem{($\Rightarrow$)}  Suppose $\lam\approx\rho$.  If $\lam$ and $\rho$ are reduced (in the sense of Remark \ref{rem:conf}), then we must have $\lam=\rho$, which implies $e'=f'$.  So now suppose $\lam$ and $\rho$ are not reduced.  In particular, $e\not=f$; cf.~\ref{Om2}.  For $\lam$ not to be reduced, it must then be the case that $e'\in\{e,f\}$, and similarly $f'\in\{e,f\}$.  If $e'=f'$ then we are done.  Otherwise, $\{e,f\} = \{e',f'\} = \{e\th_p,f\th_p\}$, and it follows that $e,f\leq p$.
%
%
\epf

It follows from Propositions \ref{prop:ssim'} and \ref{prop:degen} that the congruence $\ssim$ on $\P$ is generated by the set $\Xi''$ of all pairs $(\s,\t)\in\P\times\P$ of the form \ref{Om1} and \ref{Om2}, as well as:
\begin{enumerate}[label=\textup{\textsf{($\mathsf{\Om}$\arabic*)$''$}},leftmargin=12mm]\addtocounter{enumi}{2}
\item \label{Om3'} $\s=\lam(e,p,f)$ and $\t=\rho(e,p,f)$, for some $p\in P$, and some non-degenerate special $p$-linked pair $(e,f)$.
\end{enumerate}



\begin{figure}[h]
\begin{center}
\begin{tikzpicture}[scale=0.6]
\tikzstyle{vertex}=[circle,draw=black, fill=white, inner sep = 0.07cm]
\node (e) at (-2,0){$e$};
\node (e') at (0,2){$\phantom{'}e'$};
\node (f') at (0,-2){$f'$};
\node (f) at (2,0){$f$};
\draw[->-=0.6] (e)--(e');
\draw[->-=0.6] (e')--(f);
\draw[->-=0.6] (e)--(f');
\draw[->-=0.6] (f')--(f);
\node () at (0,-3){(a)};
%
%
\begin{scope}[shift={(6,0)}]
\node (e') at (0,2){$e=e'$};
\node (f') at (0,-2){$f'$};
\node (f) at (2,0){$f$};
\node () at (1,-3){(b)};
%
\draw[->-=0.6] (e')--(f');
\draw[->-=0.6] (e')--(f);
\draw[->-=0.6] (f')--(f);
\end{scope}
\begin{scope}[shift={(14,0)}]
\node (e) at (-2,0){$e$};
\node (e') at (0,2){$e'$};
\node (f') at (0,-2){$f=f'$};
%
\draw[->-=0.6] (e)--(e');
\draw[->-=0.6] (e)--(f');
\draw[->-=0.6] (e')--(f');
\node () at (-1,-3){(c)};
\end{scope}
\begin{scope}[shift={(20,0)}]
\node (e) at (-2,0){$e$};
\node (e') at (0,2){$\phantom{'}e'$};
\node (f') at (0,-2){$f'$};
\node (f) at (2,0){$f$};
\draw[->-=0.6] (e)--(e');
\draw[->-=0.6] (e')--(f);
\draw[->-=0.6] (e)--(f');
\draw[->-=0.6] (f')--(f);
\draw[->-=0.6] (e')--(f');
\node () at (0,-3){(d)};
\end{scope}
\end{tikzpicture}\\[7mm]
\begin{tikzpicture}[scale=0.6]
\tikzstyle{vertex}=[circle,draw=black, fill=white, inner sep = 0.07cm]
\node (e) at (-3,0){$e$};
\node (e') at (0,0){$e'$};
\node (f) at (3,0){$f$};
\node () at (0,-.6){$\rotatebox{90}{=}$};
\node () at (0,-1.2){$f'$};
\draw[->-=0.6] (e)--(e');
\draw[->-=0.6] (e')--(f);
\node () at (0,-2.25){(e)};
\begin{scope}[shift={(10,0)}]
\node (e') at (-2,0){$e$};
\node () at (-1.95,-.6){$\rotatebox{90}{=}$};
\node () at (-2,-1.2){$\phantom{'}e'$};
\node (f') at (2,0){$f$};
\node () at (2.05,-.6){$\rotatebox{90}{=}$};
\node () at (2,-1.2){$\phantom{'}f'$};
%
\draw[->-=0.6] (e')--(f');
\node () at (0,-2.25){(f)};
\end{scope}
\end{tikzpicture}
\caption{Diamonds, triangles and degeneracy of linked pairs of projections; see Remark \ref{rem:nondeg}.}
\label{fig:LP}
\end{center}
\end{figure}
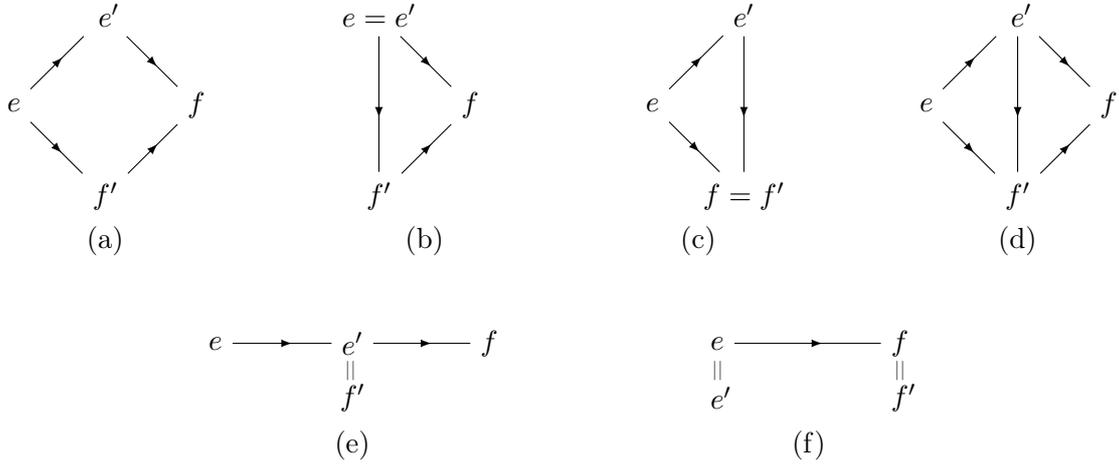

\begin{rem}\label{rem:nondeg}
Consider a non-degenerate $p$-linked pair $(e,f)$, and write $e'=e\th_p$ and $f'=f\th_p$, and also $\lam=\lam(e,p,f)=(e,e',f)$ and $\rho=\rho(e,p,f)=(e,f',f)$.  By Proposition \ref{prop:degen} there are three possibilities:
\ben
\item \label{nondeg1} $\{e,f,e',f'\}$ has size $4$,
\item \label{nondeg2} $\{e,f,e',f'\}$ has size $3$ and $e\leq p$ (i.e.~$e=e'$),
\item \label{nondeg3} $\{e,f,e',f'\}$ has size $3$ and $f\leq p$ (i.e.~$f=f'$),
\een
with $(e,f)$ being special in the second and third.
In case \ref{nondeg1}, identification of the paths $\lam$ and~$\rho$ (cf.~\ref{Om3}) amounts to commutativity of the diamond in Figure \ref{fig:LP}(a), in the groupoid $\ol\C$.
In case \ref{nondeg2}, we have $\lam\approx(e,f)$, so equating $\lam$ and $\rho$ amounts to  commutativity of the triangle in Figure \ref{fig:LP}(b).  Case~\ref{nondeg3} corresponds to Figure \ref{fig:LP}(c).  Considering again case \ref{nondeg1}, the proof of Proposition \ref{prop:ssim'} showed that $(e,f')$ and $(e',f)$ are special $p$-linked pairs; since $e',f'\leq p$, these are of types \ref{nondeg3} and~\ref{nondeg2}, respectively.  This means that the two triangles commute in Figure~\ref{fig:LP}(d); commutativity of these triangles of course implies commutativity of the outer diamond in the same diagram, as per the final line of the proof of Proposition \ref{prop:ssim'}.

On the other hand, a degenerate $p$-linked pair $(e,f)$ leads to one of diagrams (e) or (f) in Figure~\ref{fig:LP}, in the cases $e'=f'$ and $e,f\leq p$, respectively.  These diagrams already commute in~$\C$, and we note that (e) and (f) picture the generic case in which all projections displayed are distinct (it is possible to have even more collapse). 
\end{rem}

\subsection{Graphs and complexes}\label{subsect:top}

We are now in a position to give the promised topological interpretation of the chain and reduced chain groupoids $\C=\C(P)$ and $\ol\C=\ol\C(P)$.
We begin by defining a graph $G_P$ with:
\bit
\item vertex set $P$ (a given projection algebra), and
\item an (undirected) edge $\{p,q\}$ for every pair of distinct $\F$-related projections $(p,q)\in{\F}\sm\De_P$.
\eit
We then define two $2$-complexes, $K_P$ and $K_P'$.  These both have $G_P$ as their $1$-skeleton.
\bit
\item The complex $K_P$ has a $2$-cell with boundary $(e,e',f,f',e)$ for every $p$-linked pair $(e,f)$, where as usual we write $e'=e\th_p$ and $f'=f\th_p$.
\item The complex $K_P'$ is a sub-complex of $K_P$, and contains only the $2$-cells corresponding to non-degenerate special $p$-linked pairs $(e,f)$.  Such a cell has boundary $(e,f,f',e)$ or $(e,e',f,e)$ when the pair is of type \ref{nondeg2} or \ref{nondeg3}, as enumerated in Remark \ref{rem:nondeg}.
\eit
In particular, all cells in $K_P'$ are triangles.

The following result is essentially the observation that our categories $\C=\P/\Om^\sharp$ and ${\ol\C=\P/\Xi^\sharp=\P/\Xi''^\sharp}$ are constructed in precisely the same way as the fundamental groupoids of the above graph and complexes.  See for example \cite[Section 6]{Higgins2005}.

\begin{thm}\label{thm:fundgpd}
For any projection algebra $P$ we have
\[
\C(P) \cong \pi_1(G_P) \AND \ol\C \cong \pi_1(K_P) \cong \pi_1(K_P'),
\]
as (unordered) groupoids.  \epfres
\end{thm}

This has an immediate important consequence concerning the subgroups of the free regular $*$-semigroups $\PG(P)$.
It follows from general semigroup theory that in an arbitrary semigroup~$S$, maximal subgroups 
are precisely the $\H$-classes of idempotents; see \cite[Exercise 2.3.1]{CPbook}.
The relation $\H$ is one of the five Green's equivalences (see \cite[Section 2.1]{Howie1995}), but we do not require its actual definition, only that the $\H$-class of
an idempotent $e$ has the form
\[
H_e=\set{ s\in S}{se=es=s \text{ and } st=ts=e\text{ for some } t\in S}.
\]
If $e$ and $f$ are $\R$-equivalent idempotents, which here can be simply taken to mean $e\rrarr f$ using the notation introduced in Subsection \ref{subsect:Eprelim}, then $H_e\cong H_f$ \cite[Proposition 2.3.6]{Howie1995}.
Specialising to the case where $S$ is a regular $*$-semigroup, in which every idempotent $e$ is $\R$-related to the projection $ee^*$ \cite[Theorem 2.2]{NS1978}, we see that 
every maximal subgroup of $S$ is isomorphic to one of the form $H_p$ with $p\in \bP(S)$, and that
$H_p=\{ s\in S\colon ss^*= s^*s=p\}$.
Finally, in the special case that $S=\PG(P)$ for a projection algebra $P$, we see from \eqref{eq:spmult} that the multiplication in $S$ restricted to $H_p$ is precisely the same as the composition in the reduced chain groupoid $\ol\C=\ol\C(P)$ restricted to $\ol\C(p,p)$. Combining this with Theorem~\ref{thm:fundgpd}, we obtain the following as a corollary:

\begin{thm}
\label{thm:fundgp}
Let $P$ be a projection algebra. Every maximal subgroup of $\PG(P)$ is isomorphic to a fundamental group of $K_P$, or (equivalently) of $K_P'$. Specifically, the $\H$-class of any projection $p\in P$ is isomorphic to $\pi_1(K_P,p)\cong\pi_1(K_P',p)$.  \epfres
\end{thm}

We conclude this section and the paper by recasting the examples from Sections \ref{sect:AGa} and \ref{sect:TLn} in the terminology of this section, and adding two new ones.  As always, the reader should have at the back of their mind the category isomorphism between regular $*$-semigroups and chained projection groupoids (Subsection \ref{subsect:iso}), and that under this isomorphism the free regular $*$-semigroup $\PG(P)$ and the reduced chain groupoid $\ol\C(P)$ correspond to each other (Section \ref{sect:CP}). The exposition will freely move between these two structures.

\begin{eg}
Let $\Ga=(P,E)$ be a symmetric, reflexive digraph, $A_\Ga$ the adjacency semigroup,~$P_0$ its projection algebra, and $B_\Ga=\PG(P_0)$ the bridging path semigroup, as in Subsection \ref{subsect:AGa} and Section \ref{sect:AGa}.  The graph $G_{P_0}$ is the simple, undirected reduct of $\Ga$ with $0$ added as a new, isolated vertex.  As we saw in Section \ref{sect:AGa}, all linked pairs in~$P_0$ are degenerate, and so $K_{P_0}=K_{P_0}'=G_{P_0}$.  In particular, $\C=\ol\C$ can be realised as the fundamental groupoid of the graph $G_{P_0}$,
and Theorem \ref{thm:fundgp} implies that all maximal subgroups of $B_\Ga=\PG(P_0)$ are free.  This gives another way to see why $\PG(P)$ is finite in Example \ref{eg:2x2} and infinite in Example \ref{eg:3x3}.
\end{eg}

\begin{eg}
Let $P'=\{p,q,r,e\}$ be the projection algebra from Example \ref{eg:Kinyon}.  The graph~$G_{P'}$ consists of a triangle with vertices $p,q,r$, and an isolated vertex $e$.  The complex $K_{P'}=K_{P'}'$ has a $2$-cell attached to the triangle in $G_{P'}$, coming from the $e$-linked pair $(p,r)$.
\end{eg}

\begin{eg}\label{eg:TLn}
In Theorem \ref{thm:TLnPGP}, we showed that when $P=\bPP(\TL_n)$ is the projection algebra of the Temperley--Lieb monoid $\TL_n$, the free regular $*$-semigroup $\PG(P)$ is isomorphic to $\TL_n$.  
This has a topological consequence.  
The monoid $\TL_n$ is known to have no non-trivial subgroups; equivalently, it is $\H$-trivial \cite{Wilcox2007}.
Hence, all the fundamental groups of the complexes $K_P$ and~$K_P'$ are trivial, i.e.~their connected components are simply connected.
It does not seem obvious, a priori, that this ought to be the case.
\end{eg}

In the final two examples we use the topological viewpoint to gain some insight into the free regular $*$-semigroups $\PG(\bP(\M_n))$ associated to the Motzkin monoids $\M_n$, as defined in Subsection \ref{subsect:DM}, and the relationship to the idempotent-generated subsemigroups $\la E(\M_n)\ra$, which were studied in \cite{DEG2017}.

\begin{figure}[h]
\begin{center}
\scalebox{0.8}{
\begin{tikzpicture}[scale=1]
\begin{scope}
\nc\prcol{\ColA!40}
\nc\idcol{\ColA!40}
\nc\xx{1}
\nc\yy{0}
\foreach \x in {0,...,\yy} {\fill[\prcol] (\x,\yy-\x)--(\x,\yy-\x+1)--(\x+1,\yy-\x+1)--(\x+1,\yy-\x)--(\x,\yy-\x);}
\foreach \x in {0,...,\xx} {\draw (\x,0)--(\x,\xx) (0,\x)--(\xx,\x); }
\node () at (.5,.5) {\mot{1,2,3}{}};
\end{scope}
%
\begin{scope}[shift={(2,-2)}]
\nc\prcol{\ColB!40}
\nc\idcol{\ColB!40}
\nc\xx{3}
\nc\yy{2}
\foreach \x in {0,...,\yy} {\fill[\prcol] (\x,\yy-\x)--(\x,\yy-\x+1)--(\x+1,\yy-\x+1)--(\x+1,\yy-\x)--(\x,\yy-\x);}
\foreach \x in {0,...,\xx} {\draw (\x,0)--(\x,\xx) (0,\x)--(\xx,\x); }
\node () at (0.5,2.5) {\mot{1,2}{}};
\node () at (1.5,2.5) {\mott{1/1,2/3}{}{}};
\node () at (2.5,2.5) {\mott{1/2,2/3}{}{}};
\node () at (0.5,1.5) {\mott{1/1,3/2}{}{}};
\node () at (1.5,1.5) {\mot{1,3}{}};
\node () at (2.5,1.5) {\mott{1/2,3/3}{}{}};
\node () at (0.5,0.5) {\mott{2/1,3/2}{}{}};
\node () at (1.5,0.5) {\mott{2/1,3/3}{}{}};
\node () at (2.5,0.5) {\mot{2,3}{}};
\end{scope}
%
\begin{scope}[shift={(6,-4)}]
\nc\prcol{\ColC!40}
\nc\idcol{\ColC!20}
\nc\xx{5}
\nc\yy{4}
\foreach \x in {0,...,\yy} {\fill[\prcol] (\x,\yy-\x)--(\x,\yy-\x+1)--(\x+1,\yy-\x+1)--(\x+1,\yy-\x)--(\x,\yy-\x);}
\foreach \x/\y in {1/2,2/3,2/1,3/2,3/0,4/1} {\fill[\idcol] (\x,\y)--(\x,\y+1)--(\x+1,\y+1)--(\x+1,\y)--(\x,\y);}
\foreach \x in {0,...,\xx} {\draw (\x,0)--(\x,\xx) (0,\x)--(\xx,\x); }
\node () at (1.5,3.5) {\mott{1/1}{}{}};
\node () at (2.5,3.5) {\mott{1/1}{}{2/3}};
\node () at (3.5,3.5) {\mott{1/3}{}{1/2}};
\node () at (4.5,3.5) {\mott{1/3}{}{}};
\node () at (0.5,3.5) {\mott{1/2}{}{}};
\node () at (1.5,2.5) {\mott{1/1}{2/3}{}};
\node () at (2.5,2.5) {\mott{1/1}{2/3}{2/3}};
\node () at (3.5,2.5) {\mott{1/3}{2/3}{1/2}};
\node () at (4.5,2.5) {\mott{1/3}{2/3}{}};
\node () at (0.5,2.5) {\mott{1/2}{2/3}{}};
\node () at (1.5,1.5) {\mott{3/1}{1/2}{}};
\node () at (2.5,1.5) {\mott{3/1}{1/2}{2/3}};
\node () at (3.5,1.5) {\mott{3/3}{1/2}{1/2}};
\node () at (4.5,1.5) {\mott{3/3}{1/2}{}};
\node () at (0.5,1.5) {\mott{3/2}{1/2}{}};
\node () at (1.5,0.5) {\mott{3/1}{}{}};
\node () at (2.5,0.5) {\mott{3/1}{}{2/3}};
\node () at (3.5,0.5) {\mott{3/3}{}{1/2}};
\node () at (4.5,0.5) {\mott{3/3}{}{}};
\node () at (0.5,0.5) {\mott{3/2}{}{}};
\node () at (1.5,4.5) {\mott{2/1}{}{}};
\node () at (2.5,4.5) {\mott{2/1}{}{2/3}};
\node () at (3.5,4.5) {\mott{2/3}{}{1/2}};
\node () at (4.5,4.5) {\mott{2/3}{}{}};
\node () at (0.5,4.5) {\mott{2/2}{}{}};
\end{scope}
%
\begin{scope}[shift={(12,-3)}]
\nc\prcol{\ColD!40}
\nc\idcol{\ColD!20}
\nc\xx{4}
\nc\yy{3}
\foreach \x in {0,1,2,3} \foreach \y in {0,1,2,3} {\fill[\idcol] (\x,\y)--(\x,\y+1)--(\x+1,\y+1)--(\x+1,\y)--(\x,\y);}
\foreach \x in {0,...,\yy} {\fill[\prcol] (\x,\yy-\x)--(\x,\yy-\x+1)--(\x+1,\yy-\x+1)--(\x+1,\yy-\x)--(\x,\yy-\x);}
\foreach \x in {0,...,\xx} {\draw (\x,0)--(\x,\xx) (0,\x)--(\xx,\x); }
\node () at (0.5,3.5) {\mott{}{1/2}{1/2}};
\node () at (1.5,3.5) {\mott{}{1/2}{1/3}};
\node () at (2.5,3.5) {\mott{}{1/2}{2/3}};
\node () at (3.5,3.5) {\mott{}{1/2}{}};
\node () at (0.5,2.5) {\mott{}{1/3}{1/2}};
\node () at (1.5,2.5) {\mott{}{1/3}{1/3}};
\node () at (2.5,2.5) {\mott{}{1/3}{2/3}};
\node () at (3.5,2.5) {\mott{}{1/3}{}};
\node () at (0.5,1.5) {\mott{}{2/3}{1/2}};
\node () at (1.5,1.5) {\mott{}{2/3}{1/3}};
\node () at (2.5,1.5) {\mott{}{2/3}{2/3}};
\node () at (3.5,1.5) {\mott{}{2/3}{}};
\node () at (0.5,0.5) {\mott{}{}{1/2}};
\node () at (1.5,0.5) {\mott{}{}{1/3}};
\node () at (2.5,0.5) {\mott{}{}{2/3}};
\node () at (3.5,0.5) {\mott{}{}{}};
\end{scope}
\end{tikzpicture}
}
\caption{The elements of the Motzkin monoid~$\M_3$; see Example \ref{eg:M3}.}
\label{fig:egg-box_M3}
\end{center}
\end{figure}
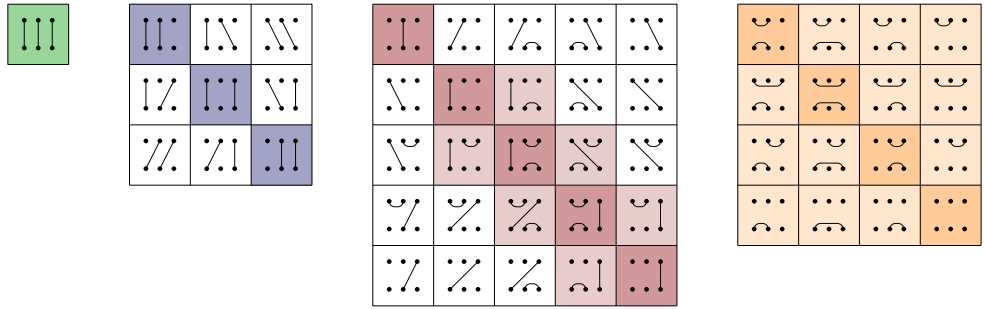

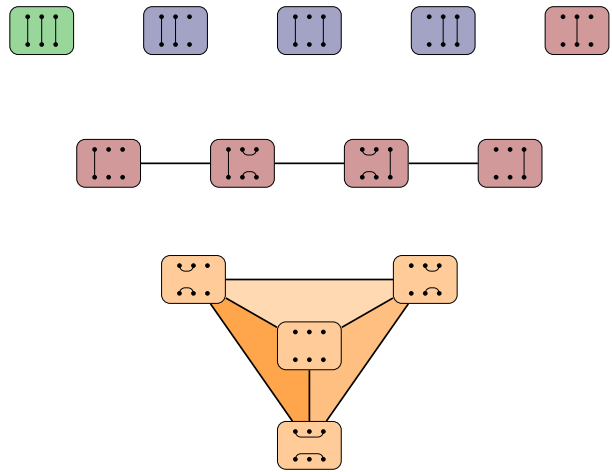
\begin{figure}[h]
\begin{center}
\scalebox{0.8}{
\begin{tikzpicture}[scale=1.1]
\nc\col{gray}
\node[draw,rounded corners,rectangle,fill=\ColA!40] () at (0,0) {\mot{1,2,3}{}};
\node[draw,rounded corners,rectangle,fill=\ColB!40] () at (2,0) {\mot{1,2}{}};
\node[draw,rounded corners,rectangle,fill=\ColB!40] () at (4,0) {\mot{1,3}{}};
\node[draw,rounded corners,rectangle,fill=\ColB!40] () at (6,0) {\mot{2,3}{}};
\node[draw,rounded corners,rectangle,fill=\ColC!40] () at (8,0) {\mot{2}{}};
\begin{scope}[shift={(1,-2)}]
\node[draw,rounded corners,rectangle,fill=\ColC!40] (a) at (0,0) {\mot{1}{}};
\node[draw,rounded corners,rectangle,fill=\ColC!40] (b) at (2,0) {\mot{1}{2/3}};
\node[draw,rounded corners,rectangle,fill=\ColC!40] (c) at (4,0) {\mot{3}{1/2}};
\node[draw,rounded corners,rectangle,fill=\ColC!40] (d) at (6,0) {\mot{3}{}};
\end{scope}
%
%
\begin{scope}[shift={(4,-4.75)}]
\fill[\ColD!70] (150:2)--(270:1.5)--(0:0);
\fill[\ColD!50] (30:2)--(270:1.5)--(0:0);
\fill[\ColD!30] (150:2)--(30:2)--(0:0);
\node[draw,rounded corners,rectangle,fill=\ColD!40] (e) at (0:0) {\mot{}{}};
\node[draw,rounded corners,rectangle,fill=\ColD!40] (f) at (150:2) {\mot{}{1/2}};
\node[draw,rounded corners,rectangle,fill=\ColD!40] (g) at (30:2) {\mot{}{2/3}};
\node[draw,rounded corners,rectangle,fill=\ColD!40] (h) at (270:1.5) {\mot{}{1/3}};
\end{scope}
\draw [thick]
(a)--(b)--(c)--(d)
(e)--(f)--(g)--(h)--(e)--(g)
(f)--(h)
;
\end{tikzpicture}
}
\caption{The complex $K_P'$, where $P=\bPP(\M_3)$ is the projection algebra of the Motzkin monoid~$\M_3$; see Example \ref{eg:M3}.}
\label{fig:M3}
\end{center}
\end{figure}

\begin{eg}\label{eg:M3}
Consider the Motzkin monoid $\M_3$, and let $P=\bPP(\M_3)$ be the projection algebra of this monoid.  The elements of $\M_3$ are shown in Figure~\ref{fig:egg-box_M3} in a so-called \emph{egg-box diagram} (see \cite{CPbook} for more details), in which the idempotents are shaded; the projections are indicated by darker shading.  The complex $K_P'$ is shown in Figure~\ref{fig:M3}.  The connected component of~$K_P'$ at the bottom of the figure contains three triangular $2$-cells, which are indicated by shading.  (The outer triangle of this component is not the boundary of a $2$-cell.)  To see for example that the `upper' triangle is a $2$-cell, denote its vertices by $e = \mot{}{1/2}$, $f = \mot{}{2/3}$ and $g = \mot{}{}$.  Then with~$p=\mot{1,2}{}$, one can check that $(e,f)$ is a special $p$-linked pair, with $e\th_p=e$ and $f\th_p = g$.  

It follows from this that the connected components of $K_P'$ are simply connected, and hence that $\PG(P)$ is $\H$-trivial.  It follows that $\PG(P)$ is isomorphic to its image under the natural morphism $\ol\id_P:\PG(P)\to\M_3$ from Theorem \ref{thm:free}.  This image is the idempotent-generated subsemigroup $\la \bEE(\M_3)\ra$.
\end{eg}

\nc\ColA{OliveGreen}
\nc\ColB{MidnightBlue}
\nc\ColC{RawSienna}
\nc\ColD{orange}

One may wonder whether the same holds for larger Motzkin monoids, but this is not the case, as our final example shows.

\begin{figure}[h]
\begin{center}
\scalebox{0.7}{
\begin{tikzpicture}[scale=2]
\tikzstyle{vertex}=[draw,rounded corners,rectangle,fill=\col!40]
\nc\col{\ColD}
\nc\xx{0.15}
%
\fill[\col!70] (3,1)--(4+\xx,1)--(5,2);
\fill[\col!50] (3,1)--(4+\xx,1)--(5,0);
\fill[\col!30] (5,0)--(4+\xx,1)--(5,2);
\fill[\col!70] (9,1)--(8-\xx,1)--(7,2);
\fill[\col!50] (9,1)--(8-\xx,1)--(7,0);
\fill[\col!30] (7,0)--(8-\xx,1)--(7,2);
\node[vertex] (a) at (1-\xx-\xx,1) {\moot{3}{}};
\node[vertex] (b) at (2-\xx,1) {\moot{3}{1/2}};
\node[vertex] (c) at (3,1)  {\moot{1}{2/3}};
\node[vertex] (d) at (4+\xx,1)  {\moot{1}{}};
\node[vertex] (e) at (5,2)  {\moot{1}{2/4}};
\node[vertex] (f) at (5,0)  {\moot{1}{3/4}};
\node[vertex] (g) at (7,2)  {\moot{4}{1/2}};
\node[vertex] (h) at (7,0)  {\moot{4}{1/3}};
\node[vertex] (i) at (8-\xx,1)  {\moot{4}{}};
\node[vertex] (j) at (9,1)  {\moot{4}{2/3}};
\node[vertex] (k) at (10+\xx,1)  {\moot{2}{3/4}};
\node[vertex] (l) at (11+\xx+\xx,1)  {\moot{2}{}};
%
%
\draw[thick] (a)--(b)--(c)--(d)--(f)--(e)--(c)--(f) (d)--(e) (l)--(k)--(j)--(i)--(g)--(h)--(j)--(g) (i)--(h) (e)--(g) (f)--(h);
\end{tikzpicture}
}
\caption{A connected component of the complex $K_P'$, where $P=\bPP(\M_4)$ is the projection algebra of the Motzkin monoid $\M_4$; see Example \ref{eg:M4}.}
\label{fig:M4}
\end{center}
\end{figure}
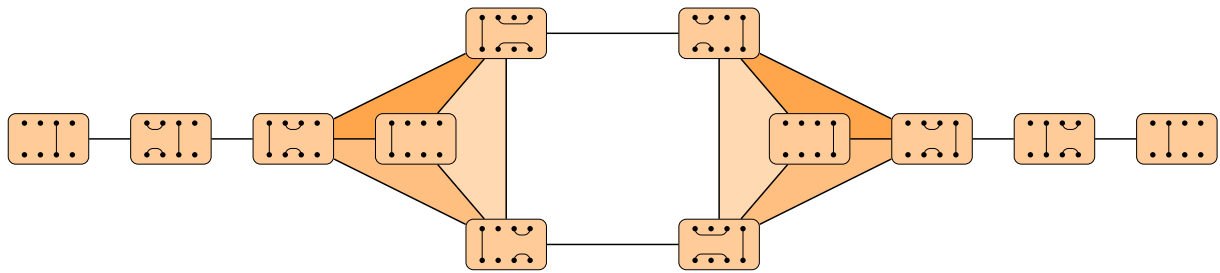

\begin{eg}\label{eg:M4}
Let $P=\bPP(\M_4)$ be the projection algebra of the Motzkin monoid $\M_4$.  The complex $K_P'$ has $35$ vertices, and eleven connected components, one of which is shown in Figure~\ref{fig:M4}; its six $2$-cells are shaded.  It is apparent that the fundamental group(oid) of this component is infinite; specifically, a loop around the central square has infinite order.  Consequently, $\PG(P)$ is infinite, and hence not isomorphic to $\la \bEE(\M_4)\ra$.
\end{eg}

Maximal subgroups of free projection-generated regular $*$-semigroups will be the main topic of our paper \cite{Paper4}.  
This will include a general theory of presentations for maximal subgroups of arbitrary $\PG(P)$, and detailed computations for $\PG(\bPP(\PP_n))$, the free regular $*$-semigroups arising from the partition monoids.  
These results will be compared and contrasted with known presentations \cite{GR2012} for maximal subgroups of the free (regular) idempotent-generated semigroups~$\IG(E)$ and~$\RIG(E)$.  This will again include explicit results for  partition monoids, which will highlight the significant difference between $\IG(\bE(\PP_n))$ and $\PG(\bP(\PP_n))$, and involve unexpected connections with \emph{twisted} partition monoids \cite{ER2022b,ER2022c,KV2023}.

\footnotesize
\def\bibspacing{-1.1pt}
\bibliography{biblio}
\bibliographystyle{abbrv}

\end{document}